\documentclass[11pt]{article}

\usepackage{latexsym}
\usepackage{amsfonts}
\usepackage{amsmath}

\usepackage[dvips]{color}

\newtheorem{thrm}{Theorem}[section]
\newtheorem{lemma}[thrm]{Lemma}
\newtheorem{prop}[thrm]{Proposition}

\newtheorem{cor}[thrm]{Corollary}
\newtheorem{remark}[thrm]{Remark}

\numberwithin{equation}{section}

\usepackage{enumerate}

\usepackage{amssymb}
\usepackage{mathtools}
\mathtoolsset{showonlyrefs}
\usepackage{mathrsfs}
\usepackage{comment}
\usepackage{hyperref}
\usepackage{xcolor}

\def\E{\mathbb{E} }
\def\P{\mathbb{P} }
\def\Q{\mathbb{Q} }
\def\R{\mathbb{R} }
\def\N{\mathbb{N} }

\makeatletter
\begin{document}
\allowdisplaybreaks

\title{\Large \bf{
Tails of extinction time and maximal displacement for
 critical branching killed L\'{e}vy process}
	\footnote{The research of this project is supported
     by the National Key R\&D Program of China (No. 2020YFA0712900).}}
\author{ \bf  Haojie Hou \hspace{1mm}\hspace{1mm}
Yan-Xia Ren\footnote{The research of this author is supported by NSFC (Grant Nos. 12071011 and 12231002) and LMEQF.
 } \hspace{1mm}\hspace{1mm} and \hspace{1mm}\hspace{1mm}
Renming Song\thanks{Research supported in part by a grant from the Simons
Foundation
(\#960480, Renming Song).}
\hspace{1mm} }
\date{}
\maketitle

\begin{abstract}
In this paper, we study asymptotic behaviors of the tails of extinction time and maximal displacement of  a critical branching killed L\'{e}vy process $(Z_t^{(0,\infty)})_{t\ge 0}$
in $\R$, in which all particles (and their descendants) are killed upon exiting $(0, \infty)$.
Let $\zeta^{(0,\infty)}$ and $M_t^{(0,\infty)}$ be the extinction time and maximal position of all the particles alive at time $t$ of this branching killed L\'{e}vy process and
define $M^{(0,\infty)}: = \sup_{t\geq 0} M_t^{(0,\infty)}$. Under the assumption that the offspring distribution belongs to
the domain of attraction of an $\alpha$-stable distribution, $\alpha\in (1, 2]$, and some moment conditions on the spatial motion, we give the decay rates of the survival probabilities
$$
\P_{y}(\zeta^{(0,\infty)}>t), \quad \P_{\sqrt{t}y}(\zeta^{(0,\infty)}>t)
$$
and the tail probabilities
$$
\P_{y}(M^{(0,\infty)}\geq x), \quad \P_{xy}(M^{(0,\infty)}\geq x).
$$
We also study the scaling limits of
$M_t^{(0,\infty)}$ and the point process $Z_t^{(0,\infty)}$ under $\P_{\sqrt{t}y}(\cdot |\zeta^{(0,\infty)}>t)$
 and $\P_y(\cdot |\zeta^{(0,\infty)}>t)$.
The scaling limits under $\P_{\sqrt{t}y}(\cdot |\zeta^{(0,\infty)}>t)$ are represented in terms of
 super killed Brownian motion.
  \end{abstract}

\medskip


\medskip


\section{Introduction and main results}

\subsection{Background and motivation}

A branching L\'{e}vy process on $\R$ is defined as follows: initially there is a particle at position $x\in R$ which moves according to a L\'{e}vy process $(\xi_t, \mathbf{P}_x)$ on $\R$. We will use $\mathbf{E}_x$ to denote expectation with respect to $\mathbf{P}_x$.
The lifetime of this particle is an exponentially distributed random variable with parameter $\beta>0$ and when it dies, this particle gives birth to a random number of offspring with law $\{p_k: k\geq 0\}$. The children of this particle  independently repeat their parent's behavior from their birthplace. The procedure goes on. We use $N(t)$ to denote the set of particles alive at time $t$ and for each $u\in N(t)$, we denote
by $X_u(t)$ the position of $u$ at time $t$. Also, for any  $u\in N(t)$ and $s\leq t$, we use $X_u(s)$ to denote the position of  $u$ or its ancestor at time $s$.  The point process $(Z_t)_{t\ge 0}$ defined by
\[
Z_t:= \sum_{u\in N(t)} \delta_{X_u(t)}
\]
is called a branching L\'evy process. We will use $\P_x$ to denote the law of this process and use $\E_x$ to denote the corresponding expectation. We will use the convention  $\P:= \P_0$ and $\E:= \E_0$.

Suppose that $m:= \sum_{k=0}^\infty kp_k\in (0,\infty)$. It is well known that the branching
L\'evy process $(Z_t)_{t\ge 0}$ will become extinct with probability one if and only if $m<1$ (subcritical) or $m=1$ and $p_1\neq 1$(critical). In this paper, we will focus on the critical case,
 that is, we always assume that $m=1$ and $p_1\neq 1$.

For any $t$, let $M_t:= \sup_{u\in N(t)} X_u(t)$ be the maximal position of all the particles alive at time $t$
and we use the convention that $M_t=-\infty$ if $N(t)=\emptyset.$ Now we define the maximal displacement and extinction time respectively by
\begin{align}\label{Def-M-zeta}
    M=\sup_{t\geq 0} M_t\quad \mbox{and}\quad \zeta:= \inf\left\{t>0: N(t)=\emptyset\right\}.
\end{align}
 Since we  always assume  $m=1$ and $p_1\neq 1$, we have $\P(M <\infty)=\P(\zeta <\infty)=1$.

Due to the homogeneity of the branching rate $\beta$ and offspring law $\{p_k: k\geq 0\}$, $\zeta$ is equal in law to the extinction time of a continuous-time Galton-Waston process with the same offspring distribution as $(Z_t)_{t\ge 0}$, so the decay rate of the survival probability $\P(\zeta>t)$ is clear.
 For example, suppose that
\begin{itemize}
	\item [{\bf(H1)}]  The offspring distribution $\{p_k: k\geq 0\}$ belongs to the domain of attraction of an $\alpha$-stable, $\alpha\in (1,2]$, distribution.
	More precisely, either there exist $\alpha\in (1,2)$ and $\kappa(\alpha)\in(0,\infty)$ such that
	\[
	\lim_{n\to\infty} n^\alpha \sum_{k=n}^\infty p_k = \kappa(\alpha),
	\]
or that (corresponding $\alpha=2$)
	\[
	\sum_{k=0}^\infty k^2 p_k<\infty.
	\]
\end{itemize}
Then it is known (see, for example, \cite{Kolmogorov38, Slack1968, Zolotarev1957}) that,  there exists a $C(\alpha)\in (0,\infty)$ such that
\begin{align}\label{Survival-prob-zeta}
	\lim_{t\to\infty} t^{\frac{1}{\alpha-1}}\P(\zeta >t) = C(\alpha).
\end{align}

The tail probability of the maximal displacement $M$ has been intensively studied in the literature.
Sawyer and Fleischman \cite{SF79} proved that under the assumption $\sum_{k=0}^\infty k^3 p_k<\infty$ and that the spatial motion $\xi$  is a standard Brownian motion, there exists a constant $\theta(2)>0$ such that
\begin{align}\label{Tail-Maximal-BBM}
	\lim_{x\to\infty} x^2\P(M\geq x) = \theta(2).
\end{align}
For corresponding results in the case of critical branching random walks with offspring distribution having finite third moment, see \cite{LS15}, and for these in the case of
critical branching L\'{e}vy processes with offspring distribution having finite third moment, see \cite{LS16, Profeta21, Profeta22}. In the case of critical branching L\'{e}vy processes
with offspring distribution belonging to the domain of attraction of an
$\alpha$-stable distribution with $\alpha\in (1, 2]$,
Hou, Jiang, Ren and Song \cite{HJRS} proved that under assumption {\bf(H1)} (although \cite{HJRS} did not deal with  the case $\alpha=2$,  the proof is actually the same as the case $\alpha\in (1,2)$, without additional assumption $\sum_{k=0}^\infty k^3 p_k<\infty$, see the argument in \cite[Theorem 1.1 below]{HJRS}) and the following assumptions on the spatial motion
\begin{itemize}
	\item [{\bf(H2)}]
	\begin{align}
		\mathbf{E}_0 (\xi_1)=0,\quad \sigma^2=\mathbf{E}_0(\xi_1^2)\in (0,\infty)
	\end{align}
\end{itemize}
and
\begin{itemize}
	\item [{\bf(H3)}]
	\begin{align}
		\mathbf{E}_0\left((	\xi_1\vee 0)^{r_0} \right)<\infty\quad \mbox{for some}\  r_0>\frac{2\alpha}{\alpha-1},
	\end{align}
\end{itemize}
where $\alpha$ is the constant in {\bf(H1)}, there exists a constant $\theta(\alpha)>0$ such that
\begin{align}\label{Tail-probability-M}
	\lim_{x\to\infty} x^{\frac{2}{\alpha-1}}\P(M\geq x)= \theta(\alpha).
\end{align}

The main concern of this paper is on critical branching L\'evy processes. If, in the critical branching L\'evy process, we kill  all particles (and their potential descendants) once they exit $(0, \infty)$, we obtain a point process $(Z_t^{(0,\infty)})_{t\ge 0}$ with
\[
Z_t^{(0,\infty)}: = \sum_{u\in N(t)} 1_{\{\inf_{s\leq t}X_u(s)>0\}}\delta_{X_u(t)}.
\]
The process $(Z_t^{(0,\infty)})_{t\ge 0}$ is called a critical branching killed L\'{e}vy process.
Let $(Z^0_t, \P_y)$ stand for a branching Markov process with spatial motion
$\xi_{t\wedge \tau_0^-}$ where
$\tau_0^-:= \inf\{t> 0: \xi_t \leq 0\}$,
branching rate $\beta$ and offspring distribution $\{p_k: k\geq 0\}$.
Then it is easy to see that for any $t, y>0$,
\begin{equation}\label{e:altint}
\big(Z_t^{(0,\infty)}, \P_y \big)\stackrel{\mathrm{d}}{=}\big(Z^0_t|_{(0, \infty)}, \P_y \big).
\end{equation}
Define
\begin{equation}\label{def-M_t}
M_t^{(0,\infty)}:= \sup_{u\in N(t): \inf_{s\leq t} X_u(s) >0 } X_u(t),\quad M^{(0,\infty)}:= \sup_{t\geq 0}M_t^{(0,\infty)}
\end{equation}
and
\begin{align}\label{Def-killed-M-zeta}
 \zeta^{(0,\infty)}:= \inf\big\{t>0:  Z_t^{(0,\infty)}((0,\infty))=0\big\}
\end{align}
with the convention $M_t^{(0,\infty )}= -\infty$ when $Z_t^{(0,\infty)}((0,\infty))=0$. When the underlying motion $\xi$ is a standard Brownian motion, Lalley and Zheng \cite{LZ15} proved that,
if $\sum_{k=0}^\infty k^3 p_k<\infty$, then
\begin{align}\label{Tail-Maximal-y}
	\lim_{x\to\infty} x^3\P_y(M^{(0,\infty)}\geq x )= \theta^{(0,\infty)}(2) y, \quad \mbox{ for all } y>0,
\end{align}
where $\theta^{(0,\infty)}(2)\in(0,\infty)$ is a constant independent of $x$ and $y$.
Comparing \eqref{Tail-Maximal-BBM} and \eqref{Tail-Maximal-y}, we see that
the tail $\P_y(M^{(0,\infty)}\geq x )$ of
critical branching killed Brownian motion decays
 to 0 in the order $x^{-3}$,
while the tail $\P_y(M\geq x )$ of
 critical branching Brownian motion decays
 to 0 in the order $x^{-2}$.
Lalley and Zheng \cite{LZ15} also showed that there exists a continuous function $(0, 1)\ni y\mapsto \theta_y^{(0,\infty)}(2)$ such that
\begin{align}\label{Tail-Maximal-xy}
	\lim_{x\to\infty} x^2 \P_{xy}(M^{(0,\infty)}\geq x) = \theta_y^{(0,\infty)}(2), \quad y\in (0, 1).
\end{align}
The argument of \cite{LZ15} relies heavily on the construction of $\P_y(M^{(0,\infty)}\geq x)$ via Weierstrass' $\mathcal{P}$-functions in the special case $p_0=p_2=\frac{1}{2}$ and a comparison argument for general offspring distributions.

 There are also some works
 on the survival probability and maximal displacement of branching killed L\'evy processes
 when $m=\sum_{k=0}^\infty kp_k>1$ and the spatial motion $\xi$ is a Brownian motion with  drift $-\mu$ where
     $\mu=\sqrt{2\beta (m-1)}$, see
 \cite{BBS15, BBS14, Kesten78, MS23,  MS20}.
  For these branching processes, $\sqrt{2\beta (m-1)}$ is the critical value of the drift  in the sense that the process will die out with probability $1$ if and only if  $\mu\geq\sqrt{2\beta (m-1)}$.
 When $p_2=1$ and $\mu=\sqrt{2\beta}$ (critical drift case),
 Berestycki, Berestycki and Schweinsberg \cite{BBS15}
 studied, among other things, the  asymptotic behavior of the position of the right-most particle
 as the position $y$ of the initial particle tends to infinity.
 In the case $\sum_{k=0}^\infty kp_k=2$ and $\mu =\sqrt{2\beta}$ (critical drift case),  the survival probability was
first studied by Kesten \cite{Kesten78}, and the result of  \cite{Kesten78} was later refined in \cite{BBS14, MS20}.
For the (all-time) maximal displacement $M^{(0,\infty)}=\max_{s\geq 0} M_s^{(0,\infty)}$
and the time when this all-time maximum is achieved
$m^{(0,\infty)}:=\mbox{arg max}_{s\geq 0} M_s^{(0,\infty)}$,
 Maillard and Schweinsberg \cite{MS23} proved the weak convergence for the conditioned law of $(M^{(0,\infty)}, m^{(0,\infty)})$ on  the event $\{\zeta^{(0,\infty)}>t\}$.

The purpose of the paper is to study the asymptotic behaviors of the tails of the survival probability and maximal displacement of
critical branching killed L\'evy processes.
More precisely, our goals are as follows:
\begin{itemize}
	\item[(i)] generalize \eqref{Tail-Maximal-y} and \eqref{Tail-Maximal-xy}  to critical 
	branching killed L\'evy processes 
	with offspring distribution satisfying {\bf(H1)} and spatial motion satisfying {\bf(H2)}--{\bf(H4)}, with {\bf(H4)} given in Section \ref{Main-results} below;
	\item[(ii)] find the exact decay rate of the survival probability $\P_y(\zeta^{(0,\infty)}>t)$;
	\item[(iii)] give probabilisitic interpretations of the limit in the generalization of \eqref{Tail-Maximal-xy}  and the limit of the survival probability when the initial position is at $\sqrt{t}y$ for fixed $y>0$;
	\item[(iv)] find scaling limits of $Z^{(0, \infty)}_t$ and $M^{(0, \infty)}_t$ under law $\P_y\big(\cdot |\zeta^{(0,\infty)}>t\big)$ and law $\P_{\sqrt{t}y}\big(\cdot |\zeta^{(0,\infty)}>t\big)$.
\end{itemize}

Our approach for proving the generalizations of \eqref{Tail-Maximal-y} and \eqref{Tail-Maximal-xy} is different from that of \cite{LZ15}. The probabilistic interpretations of the limits under $\P_{\sqrt{t}y}\big(\cdot |\zeta^{(0,\infty)}>t\big)$ are given in terms of
a particular critical superprocess. So in the next subsection, we will give a description of this superprocess and some basic facts about it.

\subsection{Critical super killed Brownian motion}\label{Intro-SBM}

Set $\R_+:=[0,\infty)$.
Let $\mathcal{M}_F(\R_+)$ and $\mathcal{M}_F((0, \infty))$ be the families of finite Borel measures on $\R_+$ and on $(0, \infty)$ respectively. We will use $\mathbf{0}$ to denote  the null measure on $\R_+$ and on $(0, \infty)$. Let $B_b(\R_+)$ and $B_b^+(\R_+)$ be the spaces of bounded Borel functions and non-negative bounded Borel functions on $\R_+$ respectively.
In this paper, whenever we are given
a function $f$ on $(0, \infty)$, we automatically extend it to $\R$ by setting $f(x)=0$ for $x\leq 0$.
The meanings of $B_b((0, \infty))$ and $B_b^+((0, \infty))$ are similar.
For any $f\in B_b(\R_+)$ and $\mu\in \mathcal{M}_F(\R_+)$,
we use $\langle f, \mu\rangle$ to denote the integral of $f$ with respect to $\mu$. For any $\alpha\in (1, 2]$, the function
\begin{align}\label{Stable-Branching-mechanism}
\varphi(\lambda) := \mathcal{C}(\alpha) \lambda^\alpha:=	\left\{\begin{array}{ll}
\frac{\beta \kappa(\alpha)\Gamma(2-\alpha)}{\alpha-1}  \lambda^\alpha,\quad &\mbox{when}\  \alpha\in (1,2),\\
\frac{\beta }{2}\left(\sum_{k=1}^\infty k(k-1)p_k\right) \lambda^2,\quad &\alpha=2,
\end{array}\right.
\end{align}
where $\kappa(\alpha)$ is given in {\bf(H1)} and $\Gamma(z):=\int_0^\infty t^{z-1}e^{-t}\mathrm{d}t$ is the Gamma function,  is a branching mechanism. Since $\varphi'(0)=0$, $\varphi$ is a critical branching mechanism.
For any $x\in \R_+$,
let $(W_t, \mathbf{P}_x)$ be a Brownian motion starting from $x$,
with variance $\sigma^2t$,  where  $\sigma^2$ is given in {\bf(H2)}.
Let $W^0_t:=W_{t\wedge \tau^{W, -}_0}$ be the process $W$ stopped at the first exit time $\tau^{W, -}_0$ of $(0, \infty)$.
Note that, when starting from $0$,
$W^0$ stays at $0$.

In this paper, for any $\mu \in \mathcal{M}_F(\R_+)$,  we will use
$X=\{(X_t)_{t\geq 0}; \mathbb P_\mu\}$
 to denote a superprocess with spatial motion $W^0$ and branching mechanism $\varphi$, that is, an $\mathcal{M}_F(\R_+)$-valued Markov process such that for any $f\in B_b^+(\R_+)$,
$$
-\log  \E_{\mu} \left(\exp\left\{ -\langle f, X_{t} \rangle\right\}\right)=\langle v_f^X(t, \cdot), \mu\rangle,
$$
where $(t, x)\mapsto v_f^X(t, x)$ is the unique locally bounded non-negative solution to
\begin{align}\label{Evolution-cumulant-semigroup}
v^X_f(t, x)&=\mathbf{E}_x\left(f(W^0_t) \right)- \mathbf{E}_y\Big(\int_0^t \varphi\left( v_f^X(t-s,W^0_s)\right)\mathrm{d}s\Big).
\end{align}
Taking $f\equiv \theta 1_{\R_+}$
in \eqref{Evolution-cumulant-semigroup}, the uniqueness of the solution implies that
\[
 -\log \E_{\delta_y} \left(\exp\left\{ -\langle \theta , X_r \rangle\right\}\right) = \left((\alpha-1)\mathcal{C}(\alpha) r^{\alpha}+\theta^{1-\alpha}\right)^{-\frac{1}{\alpha-1}}.
\]
Therefore, letting $\theta \to +\infty$ in the above equation, we obtain that
\begin{align}\label{Extinct-probability-superprocess}
	-\log \P_{\delta_y} (X_r=\mathbf{0}) &=  \lim_{\theta \to\infty} \left((\alpha-1)\mathcal{C}(\alpha)r^{\alpha}+\theta^{1-\alpha}\right)^{-\frac{1}{\alpha-1}}\nonumber\\
	& = \left((\alpha-1)\mathcal{C}(\alpha)r^{\alpha}\right)^{-\frac{1}{\alpha-1}},\quad \mbox{for all }r>0, 		y\in \R_+.
\end{align}

Next, we introduce the $\N$-measures associated to the superprocess $X$. Without loss of generality, we assume that $X$ is the coordinate process on
\[
\mathbb D:=\{ w= (w_t)_{t\geq 0}: w \text{ is an $\mathcal{M}_F(\R_+)$-valued
c\`{a}dl\`{a}g function on $\R_+$}
\}.
\]
We assume that $(\mathcal{F}_\infty, (\mathcal{F}_t)_{t\ge 0})$ is the natural filtration on $\mathbb D$, completed as usual with the $\mathcal{F}_\infty$-measurable and $\mathbb P_\mu$-negligible sets for every $\mu\in\mathcal{M}_F(\R_+)$. Let $\mathbb W^+_0$ be the family of
$\mathcal{M}_F(\R_+)$-valued c\`{a}dl\`{a}g functions on $(0, \infty)$ with $\mathbf{0}$ as a trap and with
$\lim_{t\downarrow0}w_t= \mathbf{0}$.
Note that $\mathbb W^+_0$ can be regarded as a subset of $\mathbb{D}$.

By \eqref{Extinct-probability-superprocess}, $\P_{\delta_y} (X_t=\mathbf{0}) >0$ for all $t>0$ and  $y\in \R_+$,  which implies that there exists a unique family of $\sigma$-finite measures $\{\mathbb N_y; y\in \R_+ \}$ on $\mathbb W^+_0$ such that for any $\mu\in \mathcal {M}_F(\R_+)$,
if ${\mathcal N}(\mathrm{d}w)$ is a Poisson random measure on $\mathbb W^+_0$ with intensity measure
$$
\mathbb N_\mu(\mathrm{d}w):=\int_{\R_+} \mathbb N_y(\mathrm{d}w)\mu(\mathrm{d}y),
$$
then the process defined by
$$
\widehat X_0:=\mu, \quad \widehat X_t:=\int_{\mathbb W^+_0 }w_t{\mathcal N}(\mathrm{d}w), \quad t>0,
$$
is a realization of the superprocess $X=\{(X_t)_{t\geq 0}; \mathbb P_\mu\}$.
Furthermore, for any $t>0$, $y\in \R_+$ and  $f\in B^+_b(\R_+)$,
\begin{align}\label{N-measure-equation}
 \mathbb N_y\left(1- \exp\left\{-\langle f, w_t \rangle \right\}\right)= -\log \mathbb{E}_{\delta_y}\left(\exp\left\{-\langle f, X_t\rangle\right\}\right),
\end{align}
see \cite[Theorems 8.27 and 8.28]{LZ}.
 $\{\mathbb N_y; y\in \R_+\}$
are called the $\mathbb N$-measures associated to $\{\mathbb P_{\delta_y}; y\in \R_+\}$. One can also see \cite{DyKu} for the definition of  $\{\mathbb N_y; y\in \R_+\}$.
Note that for any $y>0$, $\P_{\delta_y} (X_1((0,\infty))=0) >0$.
Thus by  \eqref{N-measure-equation}, we see that $\mathbb N_y(w_1((0,\infty))\neq 0)<\infty$.

Now we define
\[
X_t^{(0,\infty)}:= X_t\big|_{(0,\infty)},\quad t\ge 0.
\]
$(X_t^{(0,\infty)})_{t\ge 0}$ is called a critical super killed Brownian motion.
By the definition of $X^{(0,\infty)}$ we see that  for any $t, y>0$, under $\P_{\delta_y}$,
\begin{align}\label{Coupling}
	\langle f, X_t \rangle = \langle f, X_t^{(0,\infty)} \rangle ,\quad  \mbox{for any $f\in B_b^+((0, \infty))$}.
\end{align}

\subsection{Main results}\label{Main-results}

The following condition is stronger than {\bf(H3)} since it requires that $-\xi$ also satisfies {\bf(H3)}:
\begin{itemize}
	\item [{\bf(H4)}] For the $\alpha \in (1,2]$ in {\bf(H1)}, it holds that
	\begin{align}
		\mathbf{E}_0\left(|\xi_1|^{r_0} \right)<\infty\quad\mbox{for some} \quad r_0 >\frac{2 \alpha}{\alpha-1}.
	\end{align}
\end{itemize}
The assumption {\bf(H4)} will be used in the proofs of Lemmas \ref{Technical-lemma-3}, \ref{lemma5} and  \ref{lemma9}. In the proof of Lemma \ref{lemma5}, we need to apply \eqref{Tail-probability-M}
 to critical branching killed L\'evy processes with spatial motion $-\xi$.
In  the case $\alpha =2$, the assumption {\bf(H4)} is the same as that in \cite{LS15} and is weaker than that in  \cite[Remark 1.4]{Neuman-Zheng}.

For any $x\in \R$, define
\begin{align}
  \tau_x^+:= \inf\{t>0: \xi_t \ge x\}\quad \mbox{and}\quad  \tau_x^-:= \inf\{t>0: \xi_t\le x\}.
\end{align}
By Lemma \ref{Prop-Renewal-Function}  below,  under condition {\bf(H2)},  we have $\mathbf{E}_x|\xi_{\tau_0^-}|<\infty$.
Define
\begin{align}\label{Def-R}
	R(x):= x- \mathbf{E}_x\big(\xi_{\tau_0^-}\big) = -\mathbf{E}_0\big(\xi_{\tau_{-x}^-}\big),
\quad x\geq 0.
\end{align}
Note that $R(x)\ge x$ and that $R(x)$ is non-decreasing in $x$.
In Lemma \ref{Prop-Renewal-Function}, we will show that $R$ is harmonic in $(0, \infty)$ with respect to the process $\xi_{t\wedge \tau_0^-}$. When $\xi$ is a Brownian motion, $R(x)\equiv x.$

 Our main results are  as follows.

\begin{thrm}\label{thm1}
	Assume  {\bf (H1)}, {\bf(H2)} and  {\bf(H4)} hold.
	
(i) For any $y>0$,
\[
\lim_{t\to\infty} t^{\frac{1}{\alpha -1}} \P_{\sqrt{t}y}(\zeta^{(0,\infty)} >t) = \N_y\left(w_1((0,\infty))\neq 0\right),
\]
 where $\N_y$ is the $\N$-measure of the super  Brownian motion  defined in Section \ref{Intro-SBM}.

(ii)  There exists a constant 	$C^{(0,\infty)}(\alpha)\in(0,\infty)$
such that for any $y>0$,
\[
\lim_{t\to\infty} t^{\frac{1}{\alpha -1} +\frac{1}{2}} \P_{y}(\zeta^{(0,\infty)} >t) = 	C^{(0,\infty)}(\alpha)R(y).
\]
\end{thrm}

\begin{remark}
 For a critical branching L\'evy process, the tail $\P(\zeta>t)$ of extinction time decays to zero like $t^{-1/(\alpha-1)}$, see \eqref{Survival-prob-zeta}. Theorem \ref{thm1} tells that, for a critical branching killed L\'evy process starting from a single particle
 at $y>0$, the tail $\P_{y}(\zeta^{(0,\infty)} >t) $ decays to zero like $t^{-1/(\alpha-1)-1/2}$, while the tail $\P_{\sqrt{t}y}(\zeta^{(0,\infty)} >t) $ reverts back to $t^{-1/(\alpha-1)}$.
\end{remark}

 For any $t>0$, we define the following scaled version of $Z^{(0, \infty)}_t$:
\begin{equation}\label{def-Z_1}
Z_1^{(0,\infty), t}:=t^{-\frac{1}{\alpha -1}}\sum_{u\in N(t)} 1_{\{\inf_{s\leq t}X_u(s)>0\}}\delta_{X_u(t)/\sqrt{t}}.
\end{equation}
The next theorem is about the limits of $Z_1^{(0,\infty), t}$
under $\P_{\sqrt{t}y}(\cdot \big| \zeta^{(0,\infty)}>t)$ and $\P_y(\cdot \big| \zeta^{(0,\infty)}>t)$ as $t\to\infty$.
It is similar in spirit to the result that  the Dawson-Watanabe process is an appropriate scaling limit of branching Markov processes (see, for instance, \cite[Proposition 4.6]{LZ}).
This result partially answers the question in \cite[Question 1.8]{EP19} when the domain is assumed to be the half line.

\begin{thrm}\label{thm2}
		Assume  {\bf (H1)}, {\bf(H2)} and  {\bf(H4)} hold.
		
		(i) For any $y>0$,
		\[
		\lim_{t\to\infty} \P_{\sqrt{t}y} \big( Z_1^{(0,\infty), t} \in \cdot \  \big| \zeta^{(0,\infty)}>t\big) = \N_y\left( w_1|_{(0,\infty)} \in \cdot \big| w_1((0,\infty))\neq 0\right),
		\]
		where $w_1|_{(0,\infty)}$ is the restriction of the random measure $w_1$ on $(0,\infty)$.
		
		(ii) There exists a random measure $(\eta_1, \P)$ on $(0,\infty)$ such that for any $y>0$,
		\[
		\lim_{t\to\infty} \P_{y} \big( Z_1^{(0,\infty), t} \in \cdot \  \big| \zeta^{(0,\infty)}>t\big) = \P(\eta_1 \in \cdot).
		\]
\end{thrm}

\begin{remark}
Powell \cite{EP19} studied critical branching diffusion processes $Z^D_t$
 killed upon exiting
a bounded domain  $D\subset \R^d$. It was proved in \cite[Theorem 1.6]{EP19} that for any $y\in D$ and  non-negative bounded continuous function $f$ on $D$,
$\frac{1}{t} \langle f, Z_t^{D}\rangle $  under $\P_y(\cdot | Z_t^{(D)}(D)>0)$
converges weakly to an exponential random variable. In \cite[Question 1.8]{EP19}, Powell asked what happens when $D$ is unbounded. Our Theorem \ref{thm2} answers this question
in the case that
$D$ is the half-plane $(0,\infty)$.
\end{remark}

The following theorem generalizes \eqref{Tail-Maximal-y}  and \eqref{Tail-Maximal-xy}
and also provides a probabilistic interpretation for the limit of the generalization of \eqref{Tail-Maximal-xy}. When specialized to the case $\alpha=2$, the next theorem also gives an alternative proof of  \cite[Theorem 6.1]{LZ15}.
Define the maximal displacement of $X^{(0,\infty)}$ by
\begin{align}
	M^{(0,\infty), X}:= \sup_{r>0} \inf\{y\in \R: X_r^{(0,\infty)}((y,\infty))=0\}.
\end{align}

\begin{thrm}\label{thm3}
	Assume  {\bf (H1)}, {\bf(H2)} and  {\bf(H4)} hold.
	
	(i) For any $y>0$, it holds that
	\begin{equation}\label{G-Maximum-y}
	\lim_{x\to\infty} x^{\frac{2}{\alpha-1}}\P_{xy}(M^{(0,\infty)}\geq x)= -\log \P_{\delta_y}\big(M^{(0,\infty),X}<1\big).
\end{equation}
	
	(ii)  Assume further that  {\bf(H4)} holds with $r_0> 1+\frac{2\alpha}{\alpha-1}$. Then there exists a constant
  $\theta^{(0,\infty)}(\alpha)\in (0,\infty)$
  such that for any $y>0$, it holds that
	\begin{equation}\label{G-Maximum-xy}
	\lim_{x\to\infty} x^{\frac{2}{\alpha-1}+1} \P_y(M^{(0,\infty)}\geq x)= \theta^{(0,\infty)}(\alpha) R(y).
	\end{equation}
\end{thrm}

The  higher moment condition in Theorem \ref{thm3} (ii) is used in  \eqref{Technical-step}.

\begin{remark}

(1)
When $\alpha=2$, Theorem \ref{thm3} (i) is consistent with \eqref{Tail-Maximal-xy},
and Theorem \ref{thm3} (ii) is  consistent with \eqref{Tail-Maximal-y}.
In Lalley and Zheng \cite{LZ15}, the offspring distribution is assumed to have finite 3rd moment.
Our assumption  {\bf(H1)} on the offspring distribution is weaker and optimal in some sense.

In  \cite{LZ15},
the constant $\theta^{(0,\infty)}(2)$ and
the limit in \eqref{G-Maximum-y} are given in terms of Weierstrass' $\mathcal{P}$-functions.
 Our limit in \eqref{G-Maximum-y} is given in terms of superprocess and limit in \eqref{G-Maximum-xy} is given in terms of $R(y)$ defined in \eqref{Def-R}.
Using the fact that $\frac{1}{y}W_t1_{\{\min_{s\leq t} W_s>0\}}$ is a martingale under $\mathbf{P}_y$, we define
\[
\frac{\mathrm{d}\mathbf{P}_y^\uparrow}{\mathrm{d} \mathbf{P}_y}\bigg|_{\sigma(W_s, s\leq t)}:=\frac{1}{y}W_t1_{\{\min_{s\leq t} W_s>0\}} .
\]
In the case  $\sigma=1$,
it is well-known that $(W_t, \mathbf{P}_y^\uparrow)$ is a Bessel-3 process. Then combining
\eqref{Feynman-Kac-of-K} (with $z= \frac{1}{2}$) and \eqref{expression-theta}, we can
give the following  probabilistic representation for $\theta^{(0,\infty)}(\alpha)$:
\begin{align}
	\theta^{(0,\infty)}(\alpha)&= K^X\Big(\frac{1}{2}\Big) \lim_{y\to0+} \frac{1}{y}\mathbf{E}_y\Big(\exp\Big\{-  \int_0^{\tau_{1/2}^{W,+}} \psi^X\left(  K^X(W_s)\right) \mathrm{d}s \Big\};
	\tau_{1/2}^{W,+}< \tau_0^{W, -} \Big)\nonumber\\
	& = 2K^X\Big(\frac{1}{2}\Big) \lim_{y\to0+} \mathbf{E}_y^\uparrow \Big(\exp\Big\{-  \int_0^{\tau_{1/2}^{W,+}} \psi^X\left(  K^X(W_s)\right) \mathrm{d}s \Big\}\Big)\nonumber\\
	& = 2K^X\Big(\frac{1}{2}\Big) \mathbf{E}_0^\uparrow \Big(\exp\Big\{-  \int_0^{\tau_{1/2}^{W,+}} \psi^X\left(  K^X(W_s)\right) \mathrm{d}s \Big\}\Big),
\end{align}
where $\psi^X(v):=\varphi(v)/v$ and $K^X(\cdot)$ is the unique solution of \eqref{PDEin(0,1)}.

(2)   \eqref{G-Maximum-xy} says that  unlike  \eqref{Tail-probability-M} for critical branching L\'evy processes, the tail $\P_y(M^{(0,\infty)}\geq x)$ decays to zero like
$x^{-2/(\alpha-1)-1}$.
\end{remark}

The following result gives a  Yaglom-type limit for $M_t^{(0,\infty)}$.
 Define
\[
M_1^{(0,\infty),X}:= \inf\{y\in \R: X_1^{(0,\infty)}((y,\infty))=0\}.
\]

\begin{thrm}\label{thm4}
	Assume  {\bf (H1)}, {\bf(H2)} and  {\bf(H4)} hold.
	
		(i)   For any $y>0$, it holds that
		\begin{align}
			&\lim_{t\to\infty} \P_{\sqrt{t}y}\Big(\frac{M_{t}^{(0,\infty)}}{\sqrt{t}} \in \cdot \big| \zeta^{(0,\infty)}> 	t\Big)= \N_y\big( M_1^{(0,\infty), X} \in \cdot \big| w_1((0,\infty))\neq 0\big).
		\end{align}
	
	(ii) Let $\eta_1$ be the random measure in Theorem \ref{thm2}(ii) and  $M^{\eta_1}:= \inf\{y\in \R: \eta_1((y,\infty))=0\}$. Then for any $y>0$,
	\begin{align}
		& \lim_{t\to\infty} \P_{ y} \Big(\frac{M_{t}^{(0,\infty)}}{\sqrt{t}}\in \cdot \big| \zeta^{(0,\infty)} >t\Big)
		= \P(M^{\eta_1}\in \cdot).
	\end{align}
\end{thrm}

 We mention in passing here that  the proof of Theorem \ref{thm4}  does not use the conclusions of Theorem \ref{thm3}. So, we  only need {\bf(H4)}, not the enhanced  version of {\bf(H4)} in Theorem  \ref{thm3}(ii).

Theorem \ref{thm4} (ii) is similar in spirit to Lalley and Shao \cite[Theorem 3]{LS15} for branching random walks.
Let $M_n$ be the maximal position of a critical branching random walk at time $n$. \cite[Theorem 3]{LS15} says that, conditioned on survival at time $n$,
$M_n/\sqrt{n}$ converges in distribution to the maximum of the support of
a random measure $Y_1$, where $Y_1$ is the conditional limit of a super-Brownian motion $X$ such that for any nonnegative bounded continuous test function $f$,
$$
\lim_{t\to\infty}\P_{\delta_0}\big(t^{-1}\langle f(\sqrt{t}\cdot) , X_t\rangle \in \cdot| X_t \neq
\mathbf{0}
\big)=\P(\langle f, Y_1\rangle \in \cdot ),
$$
see \cite[Proposition 21]{LS15}.
  Theorem \ref{thm4} (i)
  corresponds to  Lalley and Shao \cite[Theorem 3]{LS15}, and note that since there is killing at $0$, to get the conditional limit as $t\to\infty$, the starting point needs to be at $\sqrt{t}y$.

\subsection{Proof strategies and organization of the paper}

Now we sketch the main idea of the proof of Theorem \ref{thm1}. The main ideas for the proofs of Theorem \ref{thm2} and Theorem \ref{thm3} are similar, and Theorem \ref{thm4} follows from  Theorems \ref{thm1} and \ref{thm2}. For $t>0$, $s\ge 0$ and $y>0$, let
$$
v^{(t)}_\infty(s, y):=t^{\frac{1}{\alpha-1}}\P_{\sqrt{t}y}\big(Z_{ts}^{(0,\infty)}((0,\infty))>0\big)= t^{\frac{1}{\alpha-1}}\P_{\sqrt{t}y}\big(\zeta^{(0,\infty)}>ts\big).
$$
In Section \ref{Section-F-K},  we derive an integral equation for $v^{(t)}_\infty(s, y)$.
In Section \ref{Section-2.2}, we use the Feynman-Kac formula to prove some analytical properties of $v^{(t)}_\infty(s, y)$
and show that, for any $s_0\in (0, 1)$,  $\{v^{(t)}_\infty(s, y): s\ge s_0, y>0\}_{t\ge 1}$ is tight.  Then we show that the limit
$v^X_\infty(s, y)$ exists,
is unique and can be represented via the superprocess $X$.

The remainder of this paper is organized as follows. In Section \ref{pre}, we give some preliminaries. 
The proofs of the main results are given in Section \ref{Main}. The proofs of some auxiliary results used in 
Section  \ref{Main} are given in Section \ref{aux}.

In the remainder of this paper, the notation
$f(x)\lesssim g(x)$
means
that there exists some constant $C$ independent of $x$ such that $f(x)\leq C g(x)$ holds for all $x$.

\section{Preliminaries}\label{pre}

Recall that $(\xi_t, \mathbf{P}_y)$ is a L\'{e}vy process starting from $y$, and
 for any function $f$ on $(0, \infty)$, we automatically extend it to $\R$ by setting $f(x)=0$ for all $x\le 0$.
For a  random variable $X$ and events $A, B$, we will use $\E(X; A )$ and $\E(X; A,B)$ to denote $\E(X 1_{A})$ and $\E(X1_{A}1_{B})$ respectively.

\subsection{Feynman-Kac representation}\label{Section-F-K}
Define
	\[
\phi(v):=\beta\Big(\sum_{k=0}^\infty p_k(1-v)^k -(1-v)\Big),\quad v\in[0,1].
\]
Let $L$ be a random variable with law $\{p_k\}$, then by our assumption, $\E L=m=1$. By Jensen's inequality, we have
$\phi(v)=\beta \left(\E\left( (1-v)^L\right)-(1-v)\right)\geq \beta \left((1-v)^{\E L}- (1-v)\right)=0$, which implies that $\phi$ is a non-negative function on $[0,1]$.

\begin{lemma}\label{lemma1}
	For any $f\in B_b^+((0, \infty))$,
	\[
	u_f(t,y): = \E_y\Big(\exp\Big\{- \int f(y) Z_t^{(0,\infty)}(\mathrm{d}y)\Big\}\Big),\quad t> 0,  y\in \R_+,
	\]
	solves the equation
	\begin{align}
		u_f(t,y)= \mathbf{E}_y\Big(\exp\Big\{-f(\xi_{t\land \tau_0^-}) \Big\}\Big)+ \beta\mathbf{E}_y\Big(\int_0^t \Big(\sum_{k=0}^\infty p_k u_f(t-s, \xi_{s\land \tau_0^-})^k -u_f(t-s, \xi_{s\land \tau_0^-})\Big) \mathrm{d}s\Big).
	\end{align}
Consequently, $v_f(t,y):=1-u_f(t,y)$ satisfies
	\begin{align}\label{Evolution-eq-v}
		v_f(t,y)= \mathbf{E}_y\Big(1- \exp\Big\{-f(\xi_{t\land \tau_0^-}) \Big\}\Big)- \mathbf{E}_y\Big(\int_0^t \phi(v_f(t-s, \xi_{s\land \tau_0^-}))\mathrm{d}s\Big).
	\end{align}
\end{lemma}
\textbf{Proof:}
It follows from \eqref{e:altint} that,  for any $f\in B_b^+((0, \infty))$,
\[
u_f(t,y)= \E_y \Big(\exp\Big\{- \int f(y) Z_t^0(\mathrm{d}y)\Big\}\Big).
\]
By considering the first splitting time of the branching Markov process
$Z_t^0$, we get
\begin{align}
	u_f(t,y)= e^{-\beta t}\mathbf{E}_y\Big(\exp\Big\{-f(\xi_{t\land \tau_0^-})\Big\}\Big)+\beta\mathbf{E}_y\Big(\int_0^t e^{-\beta s} \sum_{k=0}^\infty p_ku_f(t-s, \xi_{s\land \tau_0^-})^k \mathrm{d}s \Big).
\end{align}
Now the first result follows from \cite[Lemma 4.1]{E.B1.}.
Equation \eqref{Evolution-eq-v} follows from the first result and the definition of $v_f$. This completes the proof of the lemma.

\hfill$\Box$

For any $t>0$,  $s\geq 0, x, y\in \R$ and $v\in [0, t^{\frac{1}{\alpha -1}}]$,
define $f_{(t)}(\cdot):= f\left(\frac{\cdot}{\sqrt{t}}\right)$ and
\begin{align}\label{Def-V-phi-psi}
	v_f^{(t)}(s,y):= t^{\frac{1}{\alpha -1}}  v_{f_{(t)}}(ts, \sqrt{t}y),
	\quad \phi^{(t)} (v):=t^{\frac{\alpha}{\alpha -1}}\phi\left(vt^{-\frac{1}{\alpha-1}}\right), \quad \psi^{(t)}(v):= \frac{\phi^{(t)}(v)}{v}
\end{align}
and
\begin{align}\label{Def-xi-t}
	 \xi_s^{(t)}:= \frac{\xi_{st}}{\sqrt{t}}, \quad \tau_x^{(t), +}:=  \inf\{s>0: \xi_s^{(t)} \ge x\}\quad \tau_x^{(t), -}:= \inf\{s>0: \xi_s^{(t)} \le x\}.
	\end{align}
With a slight abuse of notation,
 we also use $\mathbf P_y$ to denote the law of   $\{\xi_s^{(t)}, s\geq 0\}$ with $\xi_0^{(t)}=y$. Then  $(\xi_{(tr)\land \tau_0^-}, \mathbf{P}_{\sqrt{t}y})\stackrel{\mathrm{d}}{=} (\sqrt{t}\xi_{r\land \tau_0^{(t),-}}^{(t)}, \mathbf{P}_y)$.

\begin{cor}\label{cor1}
For any $t>0, y\in \R_+$ and $0\leq w \leq r$, it holds that
	\begin{align}\label{Evolution-eq-V}
	v_f^{(t)}(r,y)=\mathbf{E}_{y}\big(v_f^{(t)}(r-w,\xi_{w\land \tau_0^{(t),-}}^{(t)} )\big) - \mathbf{E}_{y}\Big(\int_0^{w} \phi^{(t)}(v_f^{(t)}(r-s, \xi_{s\land \tau_0^{(t),-}}^{(t)} ))\mathrm{d} s\Big).
	\end{align}
\end{cor}
\textbf{Proof: }
It follows from \eqref{Evolution-eq-v} and the Markov property that for any
$y\in \R_+$ and $0\le r\le t$,
\[
	v_f(t,y)= \mathbf{E}_y(v_f(r,\xi_{(t-r)\land \tau_0^-}))- \mathbf{E}_y\Big(\int_0^{t-r} \phi(v_f(t-s, \xi_{s\land \tau_0^-}))\mathrm{d}s\Big).
\]
By the equality above with $t$ replaced by $tr$, $r$ replaced by $t(r-w)$ and $f$ replaced by $f_{(t)}$, we get
\begin{align}
	& v_f^{(t)}(r,y)= t^{\frac{1}{\alpha-1}}\mathbf{E}_{\sqrt{t}y}\big(
	v_{f_{(t)}}(t(r-w), \xi_{(tw)\land \tau_0^-})\big) - t^{\frac{1}{\alpha-1}}\mathbf{E}_{\sqrt{t}y}\Big(\int_0^{tw} \phi(v_{f_{(t)}}(tr-s, \xi_{s\land \tau_0^-}))\mathrm{d} s\Big)\nonumber\\
	& =  \mathbf{E}_{y}\big(v_f^{(t)}(r-w, \xi_{w\land \tau_0^{(t),-}}^{(t)})\big) - \mathbf{E}_{y}\Big(\int_0^{w} \phi^{(t)}(v_f^{(t)}(r-s,  \xi_{s\land \tau_0^{(t),-}}^{(t)}))\mathrm{d} s\Big),
\end{align}
where in the last equality we used the fact that  $(\xi_{(tr)\land \tau_0^-}, \mathbf{P}_{\sqrt{t}y})\stackrel{\mathrm{d}}{=} (\sqrt{t}\xi_{r\land \tau_0^{(t),-}}^{(t)}, \mathbf{P}_y)$.
This completes the proof of the corollary.
\hfill$\Box$
\bigskip

Taking $f=\theta 1_{(0,\infty)}(\cdot)$ in \eqref{Evolution-eq-V} and then letting $\theta \to +\infty$,
Corollary \ref{cor1} tells us that
\begin{align}\label{Def-V-infty}
v_\infty^{(t)}(r,y):= t^{\frac{1}{\alpha-1}}\P_{\sqrt{t}y}\big(Z_{tr}^{(0,\infty)}((0,\infty))>0\big)= t^{\frac{1}{\alpha-1}}\P_{\sqrt{t}y}\big(\zeta^{(0,\infty)}>tr\big)
\end{align}
 satisfies  the following equation:  for any $t>0, y\in \R_+$ and $0\leq w \leq r$,
\begin{align}\label{Evolution-eq-V2}
	v_\infty^{(t)}(r,y)=\mathbf{E}_{y}\big(v_\infty^{(t)}(r-w,\xi_{w\land \tau_0^{(t),-}}^{(t)} )\big) - \mathbf{E}_{y}\Big(\int_0^{w} \phi^{(t)}(v_\infty^{(t)}(r-s, \xi_{s\land \tau_0^{(t),-}}^{(t)} ))\mathrm{d} s\Big).
\end{align}

\begin{prop}\label{Feynman-Kac}
	For any $t>0, y\in \R_+$ and any $0< w< r$, it holds that
	\begin{align}
		v_\infty^{(t)}(r,y) = \mathbf{E}_y\Big(\exp\Big\{-\int_0^w \psi^{(t)}\big(v_\infty^{(t)}(r-s, \xi_s^{(t)})\big)\mathrm{d}s \Big\}v_\infty^{(t)} (r-w, \xi_s^{(t)}); \tau_0^{(t),-} >w\Big).
	\end{align}
	Also, for any $f\in B_b^+((0, \infty))$, it holds that
		\begin{align}
		v_f^{(t)}(r,y) = \mathbf{E}_y\Big(\exp\Big\{-\int_0^w \psi^{(t)}\big(v_f^{(t)}(r-s, \xi_s^{(t)})\big)\mathrm{d}s \Big\}v_f^{(t)} (r-w, \xi_w^{(t)});\tau_0^{(t),-} >w \Big).
		\end{align}
\end{prop}
\textbf{Proof:} For any fixed $t, y>0$, combining \eqref{Evolution-eq-V2} and the Feynman-Kac formula, we get
\[
v_\infty^{(t)}(r,y) = \mathbf{E}_y\Big(\exp\Big\{-\int_0^w \psi^{(t)}\big(v_\infty^{(t)}(r-s, \xi_{s\land \tau_0^{(t),-}}^{(t)})\big)\mathrm{d}s \Big\}v_\infty^{(t)} (r-w, \xi_{w\land \tau_0^{(t),-}}^{(t)}) \Big).
\]
Since $v_\infty^{(t)}(r-w, y)=0$ for all $y\leq 0$ and $0<w<r$,
we get the first result. The case for $v_f^{(t)}$ is similar.
\hfill$\Box$

\begin{remark}
	By the Markov property of $\xi_{r\land \tau_0^{(t),-}}^{(t)}$, for $y>0$ and $w\in [0,r]$, it holds that
	\begin{align}
		&\Upsilon_w:= \exp\Big\{-\int_0^w \psi^{(t)}\big(v_\infty^{(t)}(r-s, \xi_{s\land \tau_0^{(t),-}}^{(t)})\big)\mathrm{d}s \Big\}v_\infty^{(t)} (r-w, \xi_{w\land \tau_0^{(t),-}}^{(t)})\nonumber\\
		&= \mathbf{E}_y \Big(\exp\Big\{-\int_0^r \psi^{(t)}\big(v_\infty^{(t)}(r-s, \xi_{s\land \tau_0^{(t),-}}^{(t)})\big)\mathrm{d}s \Big\}v_\infty^{(t)} (0, \xi_{r\land \tau_0^{(t),-}}^{(t)}) \Big| \xi_{s\land \tau_0^{(t),-}}^{(t)}: s\leq w\Big).
	\end{align}
	Hence, $\{\Upsilon_w: w\in [0,r]\}$ is a $\mathbf{P}_y$-martingale.
	Thus, for any stopping time $T$ of the L\'{e}vy process $\xi_s^{(t)}$ and any $t>0, 0<w<r$, we have
	\begin{align}\label{Feynman-Kac-1}
		&v_\infty^{(t)}(r,y) = \mathbf{E}_y\left(\Upsilon_{w\land T} \right)\nonumber\\
		&=\mathbf{E}_y\Big(\exp\Big\{-\int_0^{w\land T} \psi^{(t)}\big(v_\infty^{(t)}(r-s, \xi_s^{(t)})\big)\mathrm{d}s \Big\}v_\infty^{(t)} (r-w\land T, \xi_{w\land T}^{(t)}); \tau_0^{(t),-} >w\land T  \Big).
	\end{align}
\end{remark}

For any $0<y<x$, define
\begin{equation}\label{v(y;x)}
v(y;x):= \P_{y} (M^{(0,\infty)}\geq x).
\end{equation}

\begin{prop}\label{Feynman-Kac-2}
For any $0<y<x$, it holds that
\begin{align}\label{FK-v}
	v(y ;x)=\mathbf{E}_y\Big(\exp\Big\{-\int_0^{\tau_x^+} \psi(v(\xi_s;x)) \mathrm{d}s \Big\}; \tau_x^+<\tau_0^-\Big)
\end{align}
where
	\[
	\psi(v):= \frac{\phi(v)}{v}=\frac{\beta\left(\sum_{k=0}^\infty p_k(1-v)^k -(1-v)\right)}{v},\quad v\in[0,1].
	\]
Consequently, 	for $0<y<z< x$, by the strong Markov property, we have
	\begin{align}\label{FK-v2}
        v(y ;x)=\mathbf{E}_y\Big(v(\xi_{\tau_z^+};x)\exp\Big\{-\int_0^{\tau_z^+} \psi(v(\xi_s;x)) \mathrm{d}s \Big\}; \tau_z^+<\tau_0^- \Big),
	\end{align}
	
\end{prop}
\textbf{Proof: }
Assume $0<y<x$.
Comparing the first branching time with  $\tau_0^-$, we get
\begin{align}
		v(y;x) =  &\int_0^\infty \beta e^{-\beta s} 	\mathbf{P}_y\left(\tau_x^+< \tau_0^-,  \tau_x^+\leq s\right)\mathrm{d} s \nonumber\\
	& + \int_0^\infty \beta e^{-\beta s}  \mathbf{E}_y\Big(
	\Big(1-\sum_{k=0}^\infty p_k\left(1-v(\xi_s;x)\right)^k \Big); s<\tau_0^-\land \tau_x^+ \Big)\mathrm{d}s\nonumber\\
	 = &\mathbf{E}_y\left(e^{-\beta \tau_x^+} ; \tau_x^+ < \tau_0^-\right)+ \int_0^\infty \beta e^{-\beta s}  \mathbf{E}_y\Big(
	\Big(1-\sum_{k=0}^\infty p_k\left(1-v(\xi_s;x)\right)^k \Big); s<\tau_0^-\land \tau_x^+ \Big)\mathrm{d}s.
\end{align}
By \cite[Lemma 4.1]{E.B1.},
the above equation is equivalent to
\begin{align}
	& v(y;x)+\beta \int_0^\infty  \mathbf{E}_y\left( v(\xi_s; x); s< \tau_0^-\land \tau_x^+ \right)\mathrm{d}s \nonumber\\
	&= \mathbf{P}_y\left( \tau_x^+<\tau_0^- \right)+ \beta\int_0^\infty \mathbf{E}_y \Big(1-\sum_{k=0}^\infty p_k\Big(1-v(\xi_s;x)\Big)^k ; s< \tau_x^+\land \tau_0^- \Big)\mathrm{d}s,
\end{align}
which is also equivalent to
\begin{align}
	v(y;x)& =  \mathbf{P}_y\left( \tau_x^+<\tau_0^- \right)-\beta  \int_0^\infty \mathbf{E}_y \Big(\sum_{k=0}^\infty p_k\left(1-v(\xi_s;x)\right)^k - (1-v(\xi_s;x)) ; s< \tau_x^+\land \tau_0^- \Big)\mathrm{d}s\nonumber\\
	& = \mathbf{P}_y\left( \tau_x^+<\tau_0^- \right) - \mathbf{E}_y\Big(\int_0^{\tau_x^+ \land \tau_0^-} \psi( v(\xi_s;x))v(\xi_s;x) \mathrm{d}s\Big).
\end{align}
Since $\psi(v)\geq 0$ for all $v\in [0,1]$, by the Feynman-Kac formula, we have
\begin{align}
	v(y ;x)&=\mathbf{E}_y\Big(\exp\Big\{-\int_0^{\tau_x^+\land \tau_0} \psi(v(\xi_s;x)) \mathrm{d}s \Big\}; \tau_x^+<\tau_0^-\Big)\nonumber\\
	&=\mathbf{E}_y\Big(\exp\Big\{-\int_0^{\tau_x^+} \psi(v(\xi_s;x)) \mathrm{d}s \Big\}; \tau_x^+<\tau_0^-\Big),
\end{align}
which gives \eqref{FK-v}.
\hfill$\Box$
\bigskip

For any $x>0$ and $y\in \R_+$, define
\begin{align}\label{def-K(x)}
	K^{(x)}(y):=x^{\frac{2}{\alpha-1}} v(xy;x) =x^{\frac{2}{\alpha-1}}\P_{xy}(M^{(0,\infty)}\geq x).
\end{align}
Then  $K^{(x)}(0)=0$ and that $K^{(x)}(y)= x^{\frac{2}{\alpha-1}}$ when $y\geq 1$.

\begin{lemma}\label{lemma6}
	For every $x>0$ and $0<y<z<1$, it holds that
	\[
	K^{(x)}(y)= \mathbf{E}_y\Big(\exp\Big\{-\int_0^{\tau_z^{(x^2),+}} \psi^{(x^2)}\big( K^{(x)}\big(\xi_s^{(x^2)} \big)  \big)\mathrm{d}s \Big\}  K^{(x)}\big(\xi_{\tau_z^{(x^2),+}}^{(x^2)} \big); \tau_z^{(x^2),+}< \tau_0^{(x^2),-} \Big).
	\]
\end{lemma}
\textbf{Proof: }
By the definition of $\tau^{(t), +}_x$ in \eqref{Def-xi-t},
$\big(x^2\tau^{(x^2), +}_z, \mathbf{P}_y\big)\stackrel{\mathrm{d}}{=}\left( \tau^+_{xz}, \mathbf{P}_{xy}\right).$
Therefore, combining
\eqref{FK-v2}
and the definition of $\xi^{(t)}$  in \eqref{Def-xi-t}, we get that
\begin{align}
	& K^{(x)}(y)= x^{\frac{2}{\alpha-1}}  \mathbf{E}_{xy}\Big(\exp\Big\{-\int_0^{\tau_{xz}^+} \psi(v(\xi_s;x)) \mathrm{d}s \Big\}
	v(\xi_{\tau^+_{xz}};x); \tau_{xz}^+<\tau_0^- \Big) \nonumber\\
	& = x^{\frac{2}{\alpha-1}}  \mathbf{E}_{y}\Big(\exp\Big\{-\int_0^{x^2\tau_{z}^{(x^2),+}} \psi(v(x\xi_{sx^{-2}}^{(x^2)};x)) \mathrm{d}s \Big\} v\big(x\xi_{ \tau_z^{(x^2),+}}^{(x^2)};x\big); \tau_{z}^{(x^2),+}<\tau_0^{(x^2),-}\Big) \nonumber\\
	& =  \mathbf{E}_{y}\Big(\exp\Big\{-\int_0^{\tau_{z}^{(x^2),+}} x^2\psi(v(x\xi_{s}^{(x^2)};x)) \mathrm{d}s \Big\} K^{(x)}\big(\xi_{ \tau_z^{(x^2),+}}^{(x^2)}\big); \tau_{z}^{(x^2),+}<\tau_0^{(x^2),-} \Big) \nonumber\\
	& = \mathbf{E}_y\Big(\exp\Big\{-\int_0^{\tau_z^{(x^2),+}} \psi^{(x^2)}\big( K^{(x)}\big(\xi_s^{(x^2)} \big)  \big)\mathrm{d}s \Big\}  K^{(x)}\big(\xi_{\tau_z^{(x^2),+}}^{(x^2)} \big); \tau_z^{(x^2),+}< \tau_0^{(x^2),-} \Big),
\end{align}
where in the last equality we used the fact that
\begin{align}\label{step_35}
	x^2\psi(v x^{-\frac{2}{\alpha-1}})= x^2 \frac{\phi (v x^{-\frac{2}{\alpha-1}})}{v x^{-\frac{2}{\alpha-1}}}= x^{\frac{2\alpha}{\alpha-1}}\frac{\phi (v x^{-\frac{2}{\alpha-1}})}{v } = \psi^{(x^2)}(v)
\end{align}
and the definition of $\psi^{(t)}$ in \eqref{Def-V-phi-psi}. This completes the proof.
\hfill$\Box$

\subsection{Some useful properties of L\'{e}vy processes}

In this subsection, we always assume that the L\'{e}vy process fulfills {\bf(H2)}.

\begin{lemma}\label{lemma12}
	If $\mathbf{E}_0\left( ((-\xi_1)\vee 0)^\lambda \right)<\infty$ for some $\lambda>2$, then
	\begin{align}
		\sup_{x>0} \mathbf{E}_x\Big(\left|\xi_{\tau_0^-}\right|^{\lambda-2} \Big)<\infty.
	\end{align}
	If $\mathbf{E}_0\left( (\xi_1\vee 0)^\lambda \right)<\infty$ for some $\lambda>2$, then
	\begin{align}
		\sup_{x>0} \mathbf{E}_{-x}\Big(\xi_{\tau_0^+}^{\lambda-2} \Big)<\infty.
	\end{align}
\end{lemma}
\textbf{Proof: } For the first result, see \cite[Lemma 2.1]{HJRS}. The second result follows by the first result with $\xi$ replaced by $-\xi$.

\hfill$\Box$

\begin{lemma}\label{Prop-Renewal-Function}	
	(i) For any $x>0$, it holds that $\mathbf{E}_{x}|\xi_{\tau_0^-}|<\infty$ and
\begin{align}\label{e:l2.8(i)}
		\lim_{x\to \infty}\frac{R(x)}{x} =1- \lim_{x\to\infty} \frac{\mathbf{E}_{x}\big(\xi_{\tau_0^-}\big)}{x}=1.
	\end{align}
Furthermore, $R(x)\lesssim x+1$.
	
	(ii)  $R$ is harmonic with respect to $\xi_{t\wedge \tau_0^-}$, that is,
	\begin{align}\label{R-X-2}
		R(x)= \mathbf{E}_x\left(R(\xi_s); \tau_0^->s \right), \quad s>0, x>0.
	\end{align}
\end{lemma}
\textbf{Proof: }
(i)  The first equality in \eqref{e:l2.8(i)} is an immediate consequence of the definition of $R(x)$, so we only prove the second equality of \eqref{e:l2.8(i)}.
We will use the decomposition introduced in \cite[p.208]{DM02}. Suppose
that $\pi^\xi$ is the L\'{e}vy measure of $\xi$. If $\pi^\xi(|x|>1)=0$, then $\mathbf{E}_0(|\xi_1|^3)<\infty$, which implies the boundness of $\mathbf{E}_x\left(|\xi_{\tau_0^-}|\right)$ according to Lemma \ref{lemma12}. Now assume that $\pi^\xi(|x|>1)>0$.  Let $\sigma_n$ be the $n$-th time that $\xi$ has a jump of magnitude larger than $1$, and put $\sigma_0=0$, then $\{\sigma_n-\sigma_{n-1}, n\geq 1\}$ are iid exponential random variables with parameter $\pi(\{ |x|>1 \})$. We can define a random walk $\hat{Z}_n$ given in \cite[p. 208]{DM02}:
\[
\hat{Z}_n = \xi_{\sigma_n},\quad n\geq 1,\quad \mbox{and}\quad \hat{Z}_0 =\xi_0.
\]
 Similar to \cite[(2.4) and (2.5)]{HJRS}, under {\bf(H2)}, $\mathbf{E}_0(\hat{Z}_1)=0$ and $\mathbf{E}_0\left(\hat{Z}_1^2\right)<\infty$. For $x\geq 0,$ define
\[
R^Z(x):= x- \mathbf{E}_x\big(\hat{Z}_{\tau_0^{Z,-}}\big),
\]
where $\tau_0^{Z,-}:=\inf\{n: \hat{Z}_n<0\}.$
It is well-known that,  under the assumption {\bf(H2)}, $\mathbf{E}_x| \hat{Z}_{\tau_0^{Z,-}}|<\infty$.
Using a martingale argument for the ladder heights process of $\hat{Z}$, we know that (for example, see \cite[(3.4) and (3.6)]{HRS})
$$
	\lim_{x\to \infty} \frac{R^Z(x)}{x}=1,
$$
which is equivalent to
\begin{align}\label{eq:8}
\lim_{x\to \infty} \frac{\mathbf{E}_x\big(\big| \hat{Z}_{\tau_0^{Z,-}}\big|\big)}{x}=0.
\end{align}
By \cite[p.209]{DM02}, for any $z>1$ and any $x>0$,
\begin{align}\label{eq:9}
	\mathbf{P}_x\big(|\xi_{\tau_0^-}| >z\big)\leq \mathbf{P}_x\big(\big|\hat{Z}_{\tau_0^{Z,-}}  \big|>z\big).
\end{align}
Combining  \eqref{eq:8} and \eqref{eq:9}, we conclude that
\begin{align}
	& \frac{1}{x}\mathbf{E}_x\big(|\xi_{\tau_0^-}|\big)\leq \frac{1}{x}+ \frac{1}{x}\int_1^\infty \mathbf{P}_x\big(|\xi_{\tau_0^-}| >z\big)\mathrm{d}z\nonumber\\
	& \leq \frac{1}{x}+ \frac{1}{x}\int_0^\infty \mathbf{P}_x\big(|\hat{Z}_{\tau_0^{Z,-}} | >z\big)\mathrm{d}z = \frac{1+ \mathbf{E}_x\big(\big| \hat{Z}_{\tau_0^{Z,-}}\big|\big)}{x}\stackrel{x\to \infty}{\longrightarrow}0.
\end{align}
The last assertion of (i) follows from \eqref{e:l2.8(i)} and the monotonicity of $R$.

(ii)
Note that
\begin{align}
	x= \mathbf{E}_x\big(\xi_{\tau_0^-\land t}\big) = \mathbf{E}_x\big(\xi_t ; \tau_0^->t \big) + \mathbf{E}_x\big(\xi_{\tau_0^-}; \tau_0^-<t \big).
\end{align}
Letting $t\to\infty$ in the above equation, using the definition of $R(x)$ and the Markov property, we have
\begin{align}\label{R-X}
	R(x)&= \lim_{t\to\infty} \mathbf{E}_x\left(\xi_t; \tau_0^->t  \right)= \lim_{t\to\infty} \mathbf{E}_x\left(\xi_{t+s}; \tau_0^->t+s  \right)\nonumber\\
&=\lim_{t\to\infty} \mathbf{E}_x\left(\mathbf{E}_{\xi_s}\left(\xi_{t}; \tau_0^->t  \right); \tau_0^->s\right)=\mathbf{E}_x\left(R(\xi_s), \tau_0^->s\right),
\end{align}
where in the last equality, we used dominated convergence theorem and the fact
$$
|\mathbf{E}_{x}\left(\xi_{t}; \tau_0^->t  \right)|=|x-\mathbf{E}_{x}(\xi_{\tau^-_0};\tau^-_0\leq t)|\leq x+R(x)\lesssim x+1, x>0.
$$
The proof is complete.

\hfill$\Box$

\begin{remark}\label{Remark}
     It follows from Lemma \ref{Prop-Renewal-Function}(i)  that,  under {\bf(H2)},
         $\mathbf{E}_y(|\xi_{\tau_0^-}|)\lesssim y+1$
     for any $y>0$. Similarly, replacing $\xi$ by $-\xi$, we see that
     \begin{align}
     	 \mathbf{E}_{-y}\big(\xi_{\tau_0^+}\big)\lesssim y+1,
     	 \quad\mbox{for all}\quad y>0,\quad \mbox{and}\quad \lim_{y\to\infty} \frac{\mathbf{E}_{-y}\big(\xi_{\tau_0^+}\big)}{y}=0.
     \end{align}
\end{remark}

Recall the definitions $\xi^{(t)}, \tau^{(t),+} $ and $\tau^{(t),-}$ in \eqref{Def-xi-t}.

\begin{lemma}\label{lemma11}
	(i) For any $y, s, t>0$,
	\begin{align}\label{e:lemma2.8(i)}
	\mathbf{P}_y\left(\tau_0^{(t),-} > s\right)\lesssim \frac{\sqrt{t} y +1}{\sqrt{st}}\quad \mbox{and}\quad \mathbf{P}_0\left(\tau_y^{(t),+} > s\right)\lesssim \frac{\sqrt{t} y  +1}{\sqrt{st}}.
	\end{align}
	(ii) For any $0<y<z$ and any $t>0$,
	\[
	\mathbf{P}_y\left(\tau_0^{(t),-} \leq \tau_{z}^{(t),+} \right) \lesssim \frac{\sqrt{t}(z-y)  + 1 }{\sqrt{t}z}.
	\]
\end{lemma}
\textbf{Proof: } (i) By definition of $\tau_0^{(t),-}$, we have
\begin{align}
	\mathbf{P}_y\big(\tau_0^{(t),-} > s\big) = 	\mathbf{P}_y\big(\inf_{\ell \leq s} \xi_\ell^{(t)}> 0 \big) =\mathbf{P}_{\sqrt{t}y}\big(\inf_{\ell\leq s t} \xi_\ell> 0 \big).
\end{align}
Now set $S_n:= \xi_n$ for $n\in \mathbb{N}$. We use the trivial upper bound $1$ in the case  $st\leq 1$.
 Now we assume that $st>1$, then
\[
\mathbf{P}_y\big(\tau_0^{(t),-} > s\big) \leq
\mathbf{P}_{\sqrt{t}y}\big(\inf_{j \leq [st]} S_j > 0 \big)
\lesssim\frac{\sqrt{t}y +1}{\sqrt{[st]}}\leq 2\frac{\sqrt{t} y  +1}{\sqrt{st}},
\]
where we used \cite[(2.7)]{AE} in the second inequality above. For the second inequality in \eqref{e:lemma2.8(i)}, noticing that
\[
\mathbf{P}_0\big(\tau_y^{(t),+} > s\big) =
\mathbf{P}_0\big(\sup_{\ell \leq st} \xi_\ell < \sqrt{t}y\big)= \mathbf{P}_{\sqrt{t}y}\big(\inf_{\ell \leq s t} (-\xi_\ell)> 0 \big),
\]
an argument similar to that used to prove the first inequality in \eqref{e:lemma2.8(i)} with $\xi$ replaced by $-\xi$ leads to the desired assertion.

(ii)  By the definitions of $\tau_0^{(t),-}$ and $\tau_{z}^{(t),+} $,
\[
\mathbf{P}_y\big(\tau_0^{(t),-} \leq \tau_{z}^{(t),+} \big) = \mathbf{P}_{\sqrt{t}y}\big(\tau_0^- \leq \tau_{\sqrt{t}z}^+ \big).
\]
Since $(\xi_s, \mathbf{P}_{\sqrt{t}y})$ is a  martingale with mean $\sqrt{t}y$,
it holds that
\begin{align}\label{step_38}
	\sqrt{t}y& = \mathbf{E}_{\sqrt{t}y}\big(\xi_{s\land \tau_0^- \land \tau_{\sqrt{t}z}^+} \big) =
	\mathbf{E}_{\sqrt{t}y}\big(\xi_s ; s< \tau_0^- \land \tau_{\sqrt{t}z}^+ \big) + \mathbf{E}_{\sqrt{t}y}\big(\xi_{ \tau_0^- \land \tau_{\sqrt{t}z}^+} ; s\geq  \tau_0^- \land \tau_{\sqrt{t}z}^+ \big).
\end{align}
According to Remark \ref{Remark}, we have
\begin{align}\label{finite-moment}
\mathbf{E}_{\sqrt{t}y}\big(|\xi_{ \tau_0^- \land \tau_{\sqrt{t}z}^+}| \big)&\leq \mathbf{E}_{\sqrt{t}y}\big(| \xi_{ \tau_0^- } |\big) + \mathbf{E}_{\sqrt{t}y}
     \big( \xi_{\tau_{\sqrt{t}z}^+} \big)\nonumber\\
&=  \mathbf{E}_{\sqrt{t}y} \big(| \xi_{ \tau_0^- } |\big)  + \sqrt{t}z + \mathbf{E}_{\sqrt{t}y-\sqrt{t}z}
\big( \xi_{\tau_0^+} \big)<\infty.
\end{align}
 Noticing that $|\xi_s|\leq \sqrt{t}z $ when $s<\tau_0^-\land\tau_{\sqrt{t}z}^+$, taking $s\to\infty$ in \eqref{step_38}, we get
$$
\sqrt{t}y=  \mathbf{E}_{\sqrt{t}y}\big(\xi_{ \tau_0^- \land \tau_{\sqrt{t}z}^+} \big),
$$
which implies
\begin{align}\label{Doob-Martingale}
 \mathbf{E}_{\sqrt{t}y} \big(\xi_{\tau_{\sqrt{t}z}^+}\big) -\sqrt{t}y = \mathbf{E}_{\sqrt{t}y}\big( \big(\xi_{\tau_{\sqrt{t}z}^+} -\xi_{\tau_0^-} \big);\tau_0^-\leq \tau_{\sqrt{t}z}^+\big).
\end{align}
By Remark \ref{Remark},  we conclude that
\begin{align}
		& \mathbf{P}_y\big(\tau_0^{(t),-} \leq \tau_{z}^{(t),+} \big)  =\mathbf{P}_{\sqrt{t}y}\big(\tau_0^- \leq \tau_{\sqrt{t}z}^+ \big)\leq \frac{1}{\sqrt{t}z}
		\mathbf{E}_{\sqrt{t}y}\big( \big(\xi_{\tau^+_{\sqrt{t}z}} -\xi_{\tau_0^-} \big);\tau_0^-\leq \tau^+_{\sqrt{t}z} \big)\nonumber\\
		& =  \frac{\mathbf{E}_{\sqrt{t}y} \big(\xi_{\tau_{\sqrt{t}z}^+}\big) -\sqrt{t}y}{\sqrt{t}z}= \frac{\sqrt{t}(z-y)  +  \mathbf{E}_{-(\sqrt{t}(z-y))}\big(\xi_{\tau_0^+}\big)}{\sqrt{t}z}\lesssim \frac{\sqrt{t}(z-y)  + 1 }{\sqrt{t}z},
\end{align}
which completes the proof of (ii).

\hfill$\Box$

Let $S_n$ be the random walk defined by $S_n=\xi_n$, $n\in \N$.  Sakhanenko \cite{Sak06} proved that (see also \cite[Lemma 3.6]{GX}) under the assumption
\begin{align}\label{eq:7}
	\mathbf{E}_0\big(|\xi_1|^{2+\delta}\big)<\infty \quad \mbox{for some }\quad \delta>0,
\end{align}
 we can find  a Brownian motion  $W_t$
 with variance $\sigma^2t$, starting from the origin,
 such that for any $\gamma\in (0, \frac{\delta}{2(2+\delta)})$ and any $t>1$, there exists a constant $N_*(\gamma)>1$ such that
\begin{align}\label{Coupling-BM}
	\mathbf{P}_0\Big(\sup_{0\leq s\leq 1} \left|S_{[ts]} -	W_{ts}
	\right|> \frac{1}{2}t^{\frac{1}{2}-\gamma} \Big) \leq \frac{N_*(\gamma)}{t^{\delta/2 -(2+\delta)\gamma}} = \frac{N_*(\gamma)}{t^{(\frac{1}{2}-\gamma)(\delta+2)-1}}.
\end{align}
By comparing $S_n$ with $\xi_t$ for $t\in [n, n+1]$, we immediately have the following result.

\begin{lemma}\label{Coupling-Levy-BM}
	Assume that  \eqref{eq:7} holds. Let  $W_s$ be the Brownian motion in \eqref{Coupling-BM},
	then for any $\gamma \in (0, \frac{\delta}{2(2+\delta)})$ and any $t>1$,
	\begin{align}
			\mathbf{P}_0\Big(\sup_{0\leq s\leq 1} \left|\xi_{ts} -	W_{ts}
			\right|> t^{\frac{1}{2}-\gamma} \Big) \lesssim  \frac{N_*(\gamma)}{t^{(\frac{1}{2}-\gamma)(\delta+2)-1}}.
	\end{align}
\end{lemma}
\textbf{Proof: }
By Doob's inequality,
\begin{align}
	& \mathbf{P}_0\Big(\sup_{0\leq s\leq 1} \left|\xi_{ts} -\xi_{[ts]} \right|>  \frac{1}{2} t^{\frac{1}{2}-\gamma} \Big)\leq \lceil t\rceil \mathbf{P}_0\Big(\sup_{0\leq s\leq 1} \left|\xi_{s}  \right|>  \frac{1}{2} t^{\frac{1}{2}-\gamma} \Big)\nonumber\\
	& \leq \lceil t\rceil \frac{2^{2+\delta}}{t^{(\frac{1}{2}-\gamma)(2+\delta)}} \mathbf{E}_0\Big(|\xi_1|^{2+\delta}\Big)\lesssim  \frac{1}{t^{(\frac{1}{2}-\gamma)(2+\delta)-1}}.
\end{align}
Combining this with \eqref{Coupling-BM}, we immediately get the desired result.
\hfill$\Box$
\bigskip

Recall that the function $R$ is defined in \eqref{Def-R}.
The following result for L\'{e}vy processes is analogous to the corresponding result for random walks proved in \cite[Theorem 2.9]{GX}. We postpone its proof to Section \ref{Proof: Technical-lemma-1}.

\begin{lemma}\label{Technical-lemma-1}
	Assume that \eqref{eq:7} holds.
	For any $y>0$ and any bounded continuous function $f$ on	$(0,\infty)$,
	 it holds that
	\begin{align}
		&\lim_{t\to\infty} \sqrt{t}\mathbf{E}_y\Big(f\Big(\frac{\xi_t}{\sigma\sqrt{t}}\Big) 1_{\left\{ \tau_0^->t \right\}} \Big) =\frac{2R(y)}{\sqrt{2\pi \sigma^2}}\int_0^\infty ze^{-\frac{z^2}{2}} f(z)\mathrm{d}z.
	\end{align}
	Consequently, for any $r>0$ and any bounded continuous function $f$ on
	$(0,\infty)$,
	\begin{align}
		&\lim_{t\to\infty} \sqrt{t}\mathbf{E}_{y/\sqrt{t}} \Big(f\big(\xi_r^{(t)}\big) 1_{\left\{ \tau_0^{(t),-}>r \right\}} \Big) = \lim_{t\to\infty} \sqrt{t}\mathbf{E}_y\Big(f\Big(\frac{\xi_{tr}}{\sqrt{t}}\Big) 1_{\left\{ \tau_0^->tr \right\}} \Big) \nonumber\\
		&=\frac{1}{\sqrt{r}}\frac{2R(y)}{\sqrt{2\pi \sigma^2}}\int_0^\infty ze^{-\frac{z^2}{2}} f(z\sigma \sqrt{r})\mathrm{d}z.
	\end{align}
\end{lemma}

Since $\xi_t^2- \sigma^2 t$ is a martingale, under the assumption $\mathbf{E}_0(|\xi_1|^4)<\infty$,
by the  optional stopping theorem and Lemma \ref{lemma12}, we have
\begin{align}\label{e:rsnew}
y^2+ \sigma^2  \mathbf{E}_y\left(\tau_0^- \land \tau_x^+\right)
	= \mathbf{E}_y\big(\big(\xi_{\tau_0^- \land \tau_x^+} \big)^2\big) \leq \mathbf{E}_y\big(\big(\xi_{\tau_0^- } \big)^2\big)+\mathbf{E}_y\big(\big(\xi_{\tau_x^+} \big)^2\big)<\infty.
\end{align}
The exit time estimates in the next result
will be used to prove \eqref{Technical-step}--\eqref{eq:6}.  The requirement for the $(4+\varepsilon)$th moment of $\xi$ is to ensure that the $(2+\varepsilon)$th moment of the overshoots $\xi_{\tau_0^-}$ and
 $\xi_{\tau_x^+}$
are finite, see \eqref{eq:3} below. Note that when {\bf(H4)} holds, the condition of Lemma \ref{Technical-lemma-3} is fulfilled because that for
$\alpha\in (1,2]$, we have
$r_0> (2\alpha)/(\alpha-1)=2+2/(\alpha -1)\geq 2+2=4$,
which implies that we can take $\varepsilon_0= r_0-4>0$ in the following lemma.

\begin{lemma} \label{Technical-lemma-3}
	Assume that  $\mathbf{E}_0(|\xi_1|^{4+\varepsilon_0})<\infty$ for some $\varepsilon_0>0$.
	
	(i) For any $y,z>0$, as $x\to +\infty$,
	\begin{align}\label{e:2.13i}
		&\lim_{x\to\infty} x \mathbf{P}_{y}\big( \tau_{xz}^{+}< \tau_0^{-}\big)=\lim_{x\to\infty} x\mathbf{P}_{yx^{-1}}\Big( \tau_z^{(x^2),+}< \tau_0^{(x^2),-}\Big)= \frac{R(y)}{z}.
	\end{align}
	
	(ii)
	For any  $y>0$, there exists $C>0$ such that for any $z>0$ and any
		$x>\max\{1, \frac{y}{z}\}$,
	\begin{align}\label{e:2.13ii}
		x\mathbf{E}_{yx^{-1}} \Big( \tau_{z}^{(x^2),+}; \tau_{z}^{(x^2),+}<\tau_0^{(x^2),-}  \Big)
		=\frac{1}{x}\mathbf{E}_y\big(\tau_{xz}^+ ; \tau_{xz}^+< \tau_0^- \big)
		\le  Cx z^2 \mathbf{P}_y\big( \tau_{xz}^+<\tau_0^-  \big)+  \frac{C}{x}.
	\end{align}
\end{lemma}
\textbf{Proof: } (i)
The first equality in \eqref{e:2.13i} follows immediately from the definition of $\xi^{(x^2)}_t$, so we only need prove that the first limit in  \eqref{e:2.13i} is equal to the right hand side of  \eqref{e:2.13i}.
According to \eqref{Doob-Martingale} and the definition of $R(y)$, we have
\begin{align}\label{eq:1}
	R(y)&= y- \mathbf{E}_y\big(\xi_{\tau_0^-}\big)= \mathbf{E}_y\big( \big(\xi_{\tau_{xz}^+} - \xi_{\tau_0^-}\big); \tau_{xz}^{+}< \tau_0^{-} \big) \nonumber\\
	&\geq xz\mathbf{P}_y\big( \tau_{xz}^{+}< \tau_0^{-}\big).
\end{align}
On the other hand, for any $\delta>0$,
\begin{align}\label{eq:2}
	R(y) &=\mathbf{E}_y\big( \big(\xi_{\tau_{xz}^+} - \xi_{\tau_0^-}\big); \tau_{xz}^{+}< \tau_0^{-}\Big)\\
&\leq (z+2\delta)x \mathbf{P}_y\big(\tau_{xz}^{+}< \tau_0^{-}, \xi_{\tau_0^-}> -\delta x, \xi_{\tau_{xz}^+} < (z+\delta)x\big)\nonumber\\
	&\quad + \mathbf{E}_y\big( \big(\xi_{\tau_{xz}^+} - \xi_{\tau_0^-}\big)
	\big(1_{\{ \xi_{\tau_0^-}\leq  -\delta x\}} +1_{\{\xi_{\tau_{xz}^+} \geq  (z+\delta)x\}}\big);\tau_{xz}^{+}< \tau_0^{-}\big) \nonumber\\
	&\leq (z+2\delta)x \mathbf{P}_y\big(\tau_{xz}^{+}< \tau_0^{-}\big)\nonumber\\
	&\quad + \sqrt{2}
	\sqrt{\mathbf{E}_y\big(\big(\xi_{\tau_{xz}^+} - \xi_{\tau_0^-} \big)^2 \big)}
	\sqrt{\mathbf{P}_y( \xi_{\tau_0^-}\leq  -\delta x) + \mathbf{P}_{y}\big(\xi_{\tau_{xz}^+} \geq  (z+\delta)x\big) }.
\end{align}
Note that, by Lemma \ref{lemma12}, for any fixed $y,z>0$,
$$
\mathbf{E}_y
\Big(\Big(\xi_{\tau_{xz}^+} - \xi_{\tau_0^-} \Big)^2 \Big)\leq
2\mathbf{E}_y
\big(\xi_{\tau_{xz}^+}^2+\xi_{\tau_0^-}^2\big)=2\Big(
\mathbf{E}_{y-xz}\big(\xi_{\tau_0^+}+xz\big)^2+\mathbf{E}_y\xi_{\tau_0^-}^2\Big)
\lesssim x^2.$$
 By Markov's inequality, and using   Lemma \ref{lemma12} again,
 we have that
\begin{align}\label{eq:3}
	&
	 \sqrt{\mathbf{E}_y\big(\big(\xi_{\tau_{xz}^+} - \xi_{\tau_0^-} \big)^2 \big)}
	 \sqrt{\mathbf{P}_y( \xi_{\tau_0^-}\leq  -\delta x) + \mathbf{P}_{y}\big(\xi_{\tau_{xz}^+} \geq  (z+\delta)x\big) }\nonumber\\
	&\lesssim \sqrt{x^2}\sqrt{ \frac{1}{(\delta x)^{2+\varepsilon_0}}\mathbf{E}_y\big(|\xi_{\tau_0^-}|^{2+\varepsilon_0}\big) + \frac{1}{(\delta x)^{2+\varepsilon_0}}
		    \mathbf{E}_{-xz+y}\big(\xi_{\tau_0^+}^{2+\varepsilon_0}\big)} \nonumber\\
	&\lesssim x^{-\varepsilon_0 /2}.
\end{align}
Putting \eqref{eq:3} into \eqref{eq:2} and applying \eqref{eq:1}, we obtain
\begin{align}
	&\frac{R(y)}{z+2\delta} \leq \liminf_{x\to\infty} x \mathbf{P}_y\left(\tau_{xz}^{+}< \tau_0^{-}\right)\leq \limsup_{x\to\infty} x \mathbf{P}_y\left(\tau_{xz}^{+}< \tau_0^{-}\right) \leq \frac{R(y)}{z}.
\end{align}
Since $\delta$ is arbitrary, we arrive at the desired result.

(ii)
The first equality in \eqref{e:2.13ii} follows immediately from the definition of $\xi^{(x^2)}_t$, so we only need prove the inequality in \eqref{e:2.13ii}.
According to \eqref{e:rsnew},
\begin{align}\label{eq:4}
	& \mathbf{E}_y\big(\big(\xi_{\tau_{xz}^+\land \tau_0^-}\big)^2\big) = \sigma^2 \mathbf{E}_y
	\left(\tau_{xz}^+ ; \tau_{xz}^+< \tau_0^-\right) + \sigma^2 \mathbf{E}_y \left( \tau_{0}^- ; \tau_{0}^-< \tau_{xz}^+\right) +y^2\nonumber\\
	&\geq \sigma^2 \mathbf{E}_y \left(\tau_{xz}^+;\tau_{xz}^+< \tau_0^- \right).
\end{align}
For the left-hand side of \eqref{eq:4}, by \eqref{eq:1} and Lemma \ref{lemma12},
\begin{align}
	&\mathbf{E}_y\big(\big(\xi_{\tau_{xz}^+\land \tau_0^-}\big)^2\big) =
	\mathbf{E}_y\big(\big(\xi_{\tau_{xz}^+}\big)^2; \tau_{xz}^+<\tau_0^- \big) +
	\mathbf{E}_y\big(\big(\xi_{\tau_0^-}\big)^2; \tau_0^-<\tau_{xz}^+ \big) \nonumber\\
	&\leq 	2\mathbf{E}_y\big(\big(\xi_{\tau_{xz}^+}-xz\big)^2; \tau_{xz}^+<\tau_0^-  \big)
	+ 2x^2z^2 \mathbf{P}_y\left( \tau_{xz}^+<\tau_0^-  \right)+
	\mathbf{E}_y\big(\big(\xi_{\tau_0^-}\big)^2\big) \nonumber\\
	&\lesssim
x^2z^2 \mathbf{P}_y\left( \tau_{xz}^+<\tau_0^-  \right)+1.
\end{align}
Plugging this upper bound back to \eqref{eq:4}, we conclude that
\begin{align}
	\frac{1}{x}\mathbf{E}_y \left(\tau_{xz}^+;\tau_{xz}^+< \tau_0^- \right)
	\lesssim x z^2 \mathbf{P}_y\left( \tau_{xz}^+<\tau_0^-  \right)+
	\frac{1}{x}.
\end{align}
This implies the result of the lemma.

\hfill$\Box$

\subsection{Preliminary estimates for the survival probability}\label{Section-2.2}

Recall that $\psi^{(t)}$ is defined in \eqref{Def-V-phi-psi}.

\begin{lemma}\label{lemma2}
	Assume that {\bf (H1)} and {\bf(H2)} hold.
	
	\noindent
	(i) For any $t,r, y>0$, it holds that
	\[
	v^{(t)}_\infty(r,y) \lesssim \frac{1}{r^{\frac{1}{\alpha-1}}}\land \frac{ \sqrt{t}y +1}{r^{\frac{1}{\alpha-1}+\frac{1}{2}} \sqrt{t}}.
	\]
	
	\noindent
	(ii) For any $r>0$ and 	$y\in \R_+$, it holds that
	\[
		    \psi^{(t)}\big( v^{(t)}_\infty(r,y)\big)
	\lesssim  \frac{1}{r},
	\]
	and that for each $K>0$, uniformly for  $v\in [0,K]$,
	\[
	\lim_{t\to\infty} \frac{\psi^{(t)}(v)}{v^{\alpha-1}}=  \mathcal{C}(\alpha),
	\]
	where $\mathcal{C}(\alpha)$ is given in \eqref{Stable-Branching-mechanism}.
\end{lemma}
\textbf{Proof: }
(i) Recall the definitions of $\zeta, \zeta^{(0,\infty)}$ and  $v_\infty^{(t)}(r,y)$ in \eqref{Def-M-zeta}, \eqref{Def-killed-M-zeta} and \eqref{Def-V-infty} respectively.
Since $\zeta^{(0,\infty)}\leq \zeta$, by  \eqref{Survival-prob-zeta},
\begin{align}\label{step_3}
	v^{(t)}_\infty(r,y) \leq t^{\frac{1}{\alpha-1}}\P_{\sqrt{t}y}\left(\zeta>tr\right)\lesssim t^{\frac{1}{\alpha-1}} \frac{1}{(tr)^{\frac{1}{\alpha-1}}}= \frac{1}{r^{\frac{1}{\alpha-1}}}.
\end{align}
On the other hand, taking $w=\frac{r}{2}$ in Proposition \ref{Feynman-Kac}, combining  \eqref{step_3} and Lemma \ref{lemma11}(i), we get
\begin{align}\label{step_4}
	& v^{(t)}_\infty(r,y) \leq \mathbf{E}_y \big(v_\infty^{(t)} \big(r/2, \xi_{r/2}^{(t)}\big);\tau_0^{(t),-} >r/2 \big)\lesssim
	 \frac{1}{r^{\frac{1}{\alpha-1}}}\mathbf{P}_y\big(\tau_0^{(t),-}>r/2\big)\nonumber\\
	&\lesssim \frac{1}{r^{\frac{1}{\alpha-1}}}\frac{\sqrt{t} y +1}{\sqrt{rt}} = \frac{ \sqrt{t}y +1}{r^{\frac{1}{\alpha-1}+\frac{1}{2}} \sqrt{t}}.
\end{align}
Now the first result follows easily from  \eqref{step_3} and \eqref{step_4}.

(ii) For $\alpha\in (1,2)$, by \cite[Lemma 3.1]{HJRS}, we have
$$
\lim_{v\to 0+}  \frac{\phi(v)}{v^\alpha} =\frac{\beta \kappa(\alpha)\Gamma(2-\alpha)}{\alpha-1},
$$
which implies
\begin{equation}\label{step_34}
\phi(v)\lesssim v^{\alpha},\quad v\in [0,1].
\end{equation}
When $\alpha=2$, \eqref{step_34} also holds  since we have
\[
\lim_{v\to 0+} \frac{\phi(v)}{v^2}= \phi''(0+)=\mathcal{C}(2).
\]
Therefore, by part (i),
\[
\psi^{(t)}\big(v^{(t)}_\infty(r,y) \big)\lesssim \frac{1}{v^{(t)}_\infty(r,y) }t^{\frac{\alpha}{\alpha-1}}\Big(v^{(t)}_\infty(r,y) t^{-\frac{1}{\alpha-1}}\Big)^\alpha= \big(v^{(t)}_\infty(r,y) \big)^{\alpha-1}\lesssim\frac{1}{r},\quad r>0, y\in \R.
\]
Since, for any $K>0$, $vt^{-\frac{1}{\alpha-1}}$ converges uniformly to $0$ for $v\in [0, K]$ as $t\to\infty$, we have, uniformly for $v\in [0,K]$,
\[
\lim_{t\to\infty} \frac{\psi^{(t)}(v)}{v^{\alpha-1}}=  \lim_{t\to\infty} \frac{\phi(vt^{-\frac{1}{\alpha-1}})}{(vt^{-\frac{1}{\alpha-1}})^{\alpha}} = \frac{\beta \kappa(\alpha)\Gamma(2-\alpha)}{\alpha-1} =\mathcal{C}(\alpha).
\]
Therefore, the assertion of (ii) is valid.
\hfill$\Box$

\begin{lemma}\label{lemma3}
	Assume that {\bf (H1)} and {\bf(H2)} hold.
	
	\noindent
	(i) For any fixed $r_0\in (0,1)$,  there exists a constant $N_1(r_0)>0$ such that for any $t>0$,
	\[
	\left|v_\infty^{(t)}(r, y+w)- v_\infty^{(t)}(r,y)\right|\leq N_1(r_0)
	\frac{1+\sqrt{t}w}{\sqrt{tw}},\quad r> 2r_0, y\in \R_+, w\in\left(0, r_0 \right).
	\]
	(ii)   For any fixed $r_0\in (0,1)$,  there exists a constant $N_2(r_0)>0$ such that for any $t>0$,
	\[
	\left|v^{(t)}_\infty(r, y)- v_\infty^{(t)}(r+q, y)\right|\leq N_2(r_0) \frac{1+\sqrt{t}q^{1/4}}{\sqrt{t}q^{1/8}},\quad r>2r_0,
	y\in \R_+, 	q\in\left(0,  r_0^4\right).
	\]
\end{lemma}
\textbf{Proof: } (i)
Since $v_\infty^{(t)}(r,0)=0$, the assertion for $y=0$ follows from Lemma \ref{lemma2}(i).
 So we assume $y>0$.
Note that for $0< y<z$ and $t>0$,
\begin{align}
	& \P_y(\zeta^{(0,\infty)}>t)= \P_y\big(\exists u\in N(t):\  \inf_{s\leq t} X_u(t)>0\big)\nonumber\\
	& = \P_z\Big(\exists u\in N(t):\  \inf_{s\leq t} X_u(t)>z-y\Big)\leq \P_z\Big(\exists u\in N(t):\
	\inf_{s\leq t} X_u(t)>0\Big)= \P_z(\zeta^{(0,\infty)}>t),
\end{align}
which implies that
\begin{align}   \label{step_7}
\left|v_\infty^{(t)}(r, y+w)- v_\infty^{(t)}(r,y)\right|=
 v_\infty^{(t)}(r, y+w)- v_\infty^{(t)}(r,y), \quad y, w>0.
\end{align}
Now for $y>0$ and $w>0$, taking  $T=\tau_{y+w}^{(t),+}$ in \eqref{Feynman-Kac-1}, we get
\begin{align}\label{step_5}
	& v_\infty^{(t)}(r,y) \nonumber\\
	& = \mathbf{E}_y\Big(\exp\Big\{-\int_0^{w\land \tau_{y+w}^{(t),+}} \psi^{(t)}\big(v_\infty^{(t)}(r-s, \xi_s^{(t)})\big)\mathrm{d}s \Big\}v_\infty^{(t)} (r-w\land \tau_{y+w}^{(t),+}, \xi_{w\land \tau_{y+w}^{(t),+}}^{(t)}); \tau_0^{(t),-} >w\land \tau_{y+w}^{(t),+}  \Big)\nonumber\\
	&\geq \mathbf{E}_y\Big(\exp\Big\{-\int_0^{w} \psi^{(t)}\big(v_\infty^{(t)}(r-s, \xi_s^{(t)})\big)\mathrm{d}s \Big\}v_\infty^{(t)} (r-\tau_{y+w}^{(t),+}, \xi_{ \tau_{y+w}^{(t),+}}^{(t)}); \tau_0^{(t),-} > \tau_{y+w}^{(t),+}, \tau_{y+w}^{(t),+}\leq w  \Big)\nonumber\\
	&\geq v_\infty^{(t)}(r,y+w) \mathbf{E}_y\Big(\exp\Big\{-\int_0^{w} \psi^{(t)}\big(v_\infty^{(t)}(r-s, \xi_s^{(t)})\big)\mathrm{d}s \Big\}; \tau_0^{(t),-} > \tau_{y+w}^{(t),+}, \tau_{y+w}^{(t),+}\leq w  \Big),
\end{align}
where in the last inequality we used the facts that $\P_y(\zeta^{(0,\infty)}>t)$ is decreasing in $t$ and  increasing in $y$, and $\xi_{\tau_{y+w}^{(t),+}}^{(t)}\geq y+w$.
 By Lemma \ref{lemma2}(ii), there exists a constant $\Gamma>0$
 such that for any $r_0\in (0, r/2)$ and
 $w\in (0, r_0)$,
\begin{align}\label{step_6}
\int_0^{w} \psi^{(t)}\big(v_\infty^{(t)}(r-s, \xi_s^{(t)})\big)\mathrm{d}s\leq \Gamma \int_0^{w} \frac{1}{r-s}\mathrm{d}s\leq \frac{\Gamma w}{r-w}\leq \frac{\Gamma}{r_0} w.
\end{align}
Combining \eqref{step_7}, \eqref{step_5} and \eqref{step_6}, we see that for any $ t>0, y>0, w\in\left(0, r_0\right)$ and $r>2r_0$,
\begin{align}
	&\left|v_\infty^{(t)}(r, y+w)- v_\infty^{(t)}(r,y)\right| \nonumber\\
	& \leq v_\infty^{(t)}(r,y+w) \Big(1-\mathbf{E}_y\Big(\exp\Big\{-\int_0^{w} \psi^{(t)}\big(v_\infty^{(t)}(r-s, \xi_s^{(t)})\big)\mathrm{d}s \Big\}; \tau_0^{(t),-} > \tau_{y+w}^{(t),+}, \tau_{y+w}^{(t),+}\leq w \Big)\Big)\nonumber\\
	&\leq v_\infty^{(t)}(r,y+w) \big(1-e^{-\frac{\Gamma}{r_0}w}\mathbf{P}_y\big(\tau_0^{(t),-} > \tau_{y+w}^{(t),+}, \tau_{y+w}^{(t),+}\leq w   \big)\big)\nonumber\\
	& \leq  v_\infty^{(t)}(r,y+w) \big(1-e^{-\frac{\Gamma}{r_0}w} + \mathbf{P}_y\big(\tau_0^{(t),-} \leq \tau_{y+w}^{(t),+}\big) + \mathbf{P}_y\big( \tau_{y+w}^{(t),+}>  w\big)\big).
\end{align}
Combining Lemma \ref{lemma2}(i), the inequality $1-e^{-x}\leq x$ for $x\geq 0$,
and Lemma \ref{lemma11},
we conclude from the above inequality that
\begin{align}
	&\left|v_\infty^{(t)}(r, y+w)- v_\infty^{(t)}(r,y)\right| \nonumber\\
	& \lesssim \frac{1}{r^{\frac{1}{\alpha-1}}} \left(\frac{\Gamma}{r_0}w +  \frac{\sqrt{t}w  + 1 }{\sqrt{t}(y+w)}+ \frac{\sqrt{t} w+1}{\sqrt{wt}} \right)\lesssim \frac{1}{r_0^{\frac{1}{\alpha-1}}}\left(w+  \frac{\sqrt{t}w  + 1 }{\sqrt{t}(y+w)} +\frac{\sqrt{t} w+1}{\sqrt{wt}} \right)\nonumber\\
	& \leq \frac{1}{r_0^{\frac{1}{\alpha-1}}}\left(\frac{\sqrt{t}w  + 1 }{\sqrt{t}(y+w)} +2\frac{\sqrt{t} w+1}{\sqrt{wt}} \right),
\end{align}
where the last inequality follows from $w\leq \sqrt{w} =\frac{\sqrt{t}w}{\sqrt{wt}}< \frac{\sqrt{t} w+1}{\sqrt{wt}} $.
Therefore, when $y>\sqrt{w}$, we have
\begin{align}\label{step_8}
\left|v_\infty^{(t)}(r, y+w)- v_\infty^{(t)}(r,y)\right| \lesssim   \frac{\sqrt{t}w  + 1 }{\sqrt{t}(\sqrt{w}+w)} +\frac{\sqrt{t} w+1}{\sqrt{wt}} \leq 2\frac{\sqrt{t} w+1}{\sqrt{wt}} .
\end{align}
On the other hand, when $y\leq \sqrt{w}$, using the monotonicity of
$v_\infty^{(t)}(r,y)$ in $y$ and Lemma \ref{lemma2}(i),
\begin{align}\label{step_9}
	&\left|v_\infty^{(t)}(r, y+w)- v_\infty^{(t)}(r,y)\right| \leq v_\infty^{(t)}(r, y+w)\leq v_\infty^{(t)}(r, \sqrt{w}+w)\nonumber\\
	& \lesssim \frac{ \sqrt{t} ( \sqrt{w}+w)+1}{r^{\frac{1}{\alpha-1}+\frac{1}{2}} \sqrt{t}} \leq \frac{2}{(2r_0)^{\frac{1}{\alpha-1}+\frac{1}{2}}} \frac{\sqrt{tw}+1}{\sqrt{t}}\leq \frac{2}{(2r_0)^{\frac{1}{\alpha-1}+\frac{1}{2}}} \frac{\sqrt{t}w+1}{\sqrt{tw}}.
\end{align}
Together with \eqref{step_8} and \eqref{step_9}, we complete the proof of (i).

(ii)
Since $v_\infty^{(t)}(r,0)=0$, the assertion for $y=0$ is trivial.
So we assume $y>0$. By the monotonicity property of
$v_\infty^{(t)}(r,y)$ in  $r$,
\begin{align}\label{step_10}
\left|v^{(t)}_\infty(r, y)- v_\infty^{(t)}(r+q, y)\right|=v^{(t)}_\infty(r, y)- v_\infty^{(t)}(r+q, y).
\end{align}
Using Proposition \ref{Feynman-Kac} with $r$ replaced by $r+q$ and $w$ replaced by $q$,
and an argument  similar to  that for \eqref{step_6}, we have
\begin{align}
	& 	v_\infty^{(t)}(r+q,y) = \mathbf{E}_y\Big(\exp\Big\{-\int_0^q \psi^{(t)}\big(v_\infty^{(t)}(r+q-s, \xi_s^{(t)})\big)\mathrm{d}s \Big\}v_\infty^{(t)} (r, \xi_q^{(t)}); \tau_0^{(t),-} >q \Big)\nonumber\\
	& \geq e^{-\Gamma \int_0^q \frac{1}{r+q-s}\mathrm{d}s} \mathbf{E}_y\big(v_\infty^{(t)} (r, \xi_q^{(t)}); \tau_0^{(t),-} >q \big)
	\geq e^{-\frac{\Gamma}{2r_0} q} \mathbf{E}_y\Big(v_\infty^{(t)} (r, \xi_q^{(t)}); \tau_0^{(t),-} >q  \Big)\nonumber\\
	&\geq  e^{-\frac{\Gamma}{2r_0} q}
	\mathbf{E}_y\big(v_\infty^{(t)} (r, \xi_q^{(t)}); \tau_0^{(t),-} >q, |\xi_q^{(t)}- y|< q^{1/4} \big).
\end{align}
Plugging this  into \eqref{step_10}, we obtain
\begin{align}\label{step_11}
	& \left|v^{(t)}_\infty(r, y)- v_\infty^{(t)}(r+q, y)\right| \leq v^{(t)}_\infty(r, y) - e^{-\frac{\Gamma}{2r_0} q} \mathbf{E}_y\big(v_\infty^{(t)} (r, \xi_q^{(t)}); \tau_0^{(t),-} >q, |\xi_q^{(t)}- y|< q^{1/4} \big) \nonumber\\
	& \leq v_\infty^{(t)}(r,y)\big(1- e^{-\frac{\Gamma}{2r_0} q} \mathbf{P}_y\big( \tau_0^{(t),-} >q,  |\xi_q^{(t)}- y|< q^{1/4} \big) \big)\nonumber\\
	&\quad + e^{-\frac{\Gamma}{2r_0} q} \mathbf{E}_y\big(\big| v_\infty^{(t)} (r, y) - v_\infty^{(t)} (r, \xi_q^{(t)}) \big|; \tau_0^{(t),-} >q, |\xi_q^{(t)}- y|< q^{1/4} \big).
\end{align}
By part (i),
the last term of \eqref{step_11} is bounded above by $N_1(r_0)\frac{1+\sqrt{t}q^{1/4}}{\sqrt{t}q^{1/8}}$.
Similarly, by Lemma \ref{lemma2}(i),  Doob's maximal inequality and Markov's inequality, we have
\begin{align}\label{step_56}
	&v_\infty^{(t)}(r,y)\big(1- e^{-\frac{\Gamma}{2r_0} q} \mathbf{P}_y\big( \tau_0^{(t),-} >q,  |\xi_q^{(t)}- y|< q^{1/4} \big) \big)  \nonumber\\
	& \lesssim \frac{1}{r^{\frac{1}{\alpha-1}}}\big(1- e^{-\frac{\Gamma}{2r_0} q}  +\mathbf{P}_y\big( \tau_0^{(t),-} \leq q \big)  + \mathbf{P}_y\big( |\xi_q^{(t)}- y|> q^{1/4} \big)  \big)\nonumber\\
	& \leq  \frac{1}{(2r_0)^{\frac{1}{\alpha-1}}}\Big(\frac{\Gamma}{2r_0} q  +
	 \mathbf{P}_y\big(\inf_{s\leq q} \xi_s^{(t)} < 0 \big)
	+ \mathbf{P}_0\big( |\xi_q^{(t)}|> q^{1/4} \big)  \Big)\nonumber\\
	& \lesssim q + 	
	\mathbf{P}_0 \big(\sup_{s\leq q} (-\xi_s^{(t)})^+ > y\big)
	+ \frac{1}{q^{1/2}}\mathbf{E}_0\big(\big(\xi_q^{(t)}\big)^2\big)\nonumber\\
	& \leq q+ \frac{1}{y} \sqrt{\mathbf{E}_0\big(\big(\xi_q^{(t)}\big)^2\big)}+ \mathbf{E}_0(\xi_1^2) q^{1/2}\lesssim \sqrt{q} + \frac{\sqrt{q}}{y}.
\end{align}
Therefore, when $y>q^{1/4}$, combining \eqref{step_11} and \eqref{step_56},
\begin{align}\label{step_12}
	& \left|v^{(t)}_\infty(r, y)- v_\infty^{(t)}(r+q, y)\right| \lesssim \sqrt{q}+ \frac{\sqrt{q}}{y}+ \frac{1+\sqrt{t}q^{1/4}}{\sqrt{t}q^{1/8}} \lesssim \frac{1+\sqrt{t}q^{1/4}}{\sqrt{t}q^{1/8}}.
\end{align}
On the other hand, when $y\leq q^{1/4}$, by Lemma \ref{lemma2}(i) and the monotonicity of $V_\infty^{(t)}(r,y)$ in $r$,
\begin{align}\label{step_13}
 \left|v^{(t)}_\infty(r, y)- v_\infty^{(t)}(r+q, y)\right| \leq v^{(t)}_\infty(r, y) \lesssim
  \frac{ \sqrt{t}y +1}{r^{\frac{1}{\alpha-1}+\frac{1}{2}} \sqrt{t}}
 \leq  \frac{ \sqrt{t}q^{1/4} +1}{(2r_0)^{\frac{1}{\alpha-1}+\frac{1}{2}} \sqrt{t}}.
\end{align}
Combining \eqref{step_12} and \eqref{step_13}, we complete the proof of (ii).

\hfill$\Box$

\subsection{Preliminary estimates for the tail probability of $M^{(0,\infty)}$}

Recall that, for  $x,y>0$,  $K^{(x)}(y)$ is defined by \eqref{def-K(x)}.

\begin{lemma}\label{lemma7}
	Assume that {\bf (H1)},  {\bf(H2)}  and {\bf(H3)} hold.
	
	(i) For any $x>0$ and $y\in (0,1)$, it holds that
	\[
	K^{(x)}(y) \lesssim \frac{1}{(1-y)^{\frac{2}{\alpha-1}}}.
	\]
	(ii) There exists a constant $C_*>0$ such that for any $0<y<z<1$ and $x>0$,
	\[
	\int_0^{\tau_z^{(x^2),+}} \psi^{(x^2)}\big(K^{(x)}(\xi_s^{(x^2)})\big)\mathrm{d}s \leq C_* \frac{1}{(1-z)^2} \tau_z^{(x^2),+}, \quad
	\mathbf{P}_y\mbox{-a.s.}
	\]
\end{lemma}
\textbf{Proof: } (i) Since $M^{(0,\infty)}\leq M$, by \eqref{Tail-probability-M}, we have
\begin{align}
	K^{(x)}(y) \leq x^{\frac{2}{\alpha-1}} \P_{xy}(M\geq x) = x^{\frac{2}{\alpha-1}} \P(M\geq x(1-y))\lesssim x^{\frac{2}{\alpha-1}} \frac{1}{(x(1-y))^{\frac{2}{\alpha-1}}}= \frac{1}{(1-y)^{\frac{2}{\alpha-1}}}.
\end{align}
(ii)
Combining \eqref{step_35}, \eqref{step_34} and part (i), we get that
for any $0<y<z<1$ and $x>0$,
\begin{align}
	& 	\int_0^{\tau_z^{(x^2),+}} \psi^{(x^2)}\big(K^{(x)}(\xi_s^{(x^2)})\big)\mathrm{d}s   \lesssim \int_0^{\tau_z^{(x^2),+}} x^{\frac{2\alpha}{\alpha-1}}\frac{1}{ K^{(x)}(\xi_s^{(x^2)}) }\big(K^{(x)}(\xi_s^{(x^2)}) x^{-\frac{2}{\alpha-1}}\big)^\alpha   \mathrm{d}s\nonumber\\
	&= \int_0^{\tau_z^{(x^2),+}} \big(K^{(x)}(\xi_s^{(x^2)}) \big)^{\alpha-1}   \mathrm{d}s \lesssim \int_0^{\tau_z^{(x^2),+}} \left(\frac{1}{(1-\xi_s^{(x^2)})^{\frac{2}{\alpha-1}}} \right)^{\alpha-1}   \mathrm{d}s \leq \frac{1}{(1-z)^2} \tau_z^{(x^2),+}.
\end{align}

\hfill$\Box$

\begin{lemma}\label{lemma8}
	Assume that {\bf (H1)},  {\bf(H2)}  and {\bf(H3)} hold.
	For any $r_0\in (0,1/4)$, there exists a constant $N_3(r_0)$ such that for any $x>0$, any $y\in (0, 1-2r_0)$ and any $w\in (0, r_0^2)$, it holds that
	\[
	\left| K^{(x)} (y+w)- K^{(x)}(y) \right|\leq N_3(r_0)  \left(\frac{1+xw}{(y+w)x} + 1-\mathbf{E}_0\left(\exp\left\{-N_3(r_0) \tau_{w}^{(x^2),+} \right\}\right)\right).
	\]
\end{lemma}
\textbf{Proof: }
Let $r_0\in (0, 1/4)$ and $x>0$. For $z>y$, we have
\begin{align}
	& K^{(x)}(y)=  x^{\frac{2}{\alpha-1}} \P_{xy}\big(\exists\ t>0, \ u\in N(t): \ X_u(t)\geq x,\ \inf_{s\leq t} X_u(s) >0\big) \nonumber\\
	& \leq x^{\frac{2}{\alpha-1}} \P_{xy}\big(\exists\ t>0, \ u\in N(t): \ X_u(t)\geq x - x(z-y),\
	 \inf_{s\leq t} X_u(s) >-x(z-y)\big) \nonumber\\
	& = x^{\frac{2}{\alpha-1}} \P_{xz}\big(\exists\ t>0, \ u\in N(t): \ X_u(t)\geq x,\
	\inf_{s\leq t} X_u(s) >0\big) = K^{(x)}(z).
\end{align}
Therefore,
\begin{align}\label{step_39}
	\left| K^{(x)} (y+w)- K^{(x)}(y) \right| = K^{(x)} (y+w)- K^{(x)}(y).
\end{align}
Note that for $y\in (0, 1-2r_0)$ and $w\in (0, r_0^2)$, we have $y+w< 1-2r_0+r_0^2 <1-r_0$. Thus, combining Lemma \ref{lemma6} and Lemma \ref{lemma7} (ii), we obtain that for $y\in (0, 1-2r_0)$ and $w\in (0, r_0^2)$,
\begin{align}\label{step_51}
	& K^{(x)}(y) \geq \mathbf{E}_y\Big(\exp\Big\{-C_* \frac{1}{(1-y-w)^2} \tau_{y+w}^{(x^2),+} \Big\}  K^{(x)}\big(\xi_{\tau_{y+w}^{(x^2),+}}^{(x^2)} \big); \tau_{y+w}^{(x^2),+}< \tau_0^{(x^2),-} \Big)\nonumber\\
	& \geq K^{(x)}(y+w) \mathbf{E}_y\Big(\exp\Big\{-C_* \frac{1}{r_0^2} \tau_{y+w}^{(x^2),+}\Big\}; \tau_{y+w}^{(x^2),+}< \tau_0^{(x^2),-} \Big).
\end{align}
Together with Lemma \ref{lemma7}(i) and \eqref{step_51}, for all $y\in (0, 1-2r_0)$ and
$w\in (0, r_0^2)$,
\begin{align}\label{step_37}
	 & K^{(x)} (y+w)- K^{(x)}(y)\leq  K^{(x)} (y+w)\Big(1-\mathbf{E}_y\Big(\exp\big\{-C_* \frac{1}{r_0^2} \tau_{y+w}^{(x^2),+}\big\}; \tau_{y+w}^{(x^2),+}< \tau_0^{(x^2),-} \Big)\Big)\nonumber\\
	 & \lesssim \frac{1}{(1-y-w)^{\frac{2}{\alpha-1}}}
\Big( \mathbf{P}_y\big( \tau_{y+w}^{(x^2),+}\geq  \tau_0^{(x^2),-}  \big) + 1-\mathbf{E}_y\big(\exp\big\{-\frac{C_*}{r_0^2} \tau_{w+y}^{(x^2),+} \big\}\big)\Big)\nonumber\\
	 &\lesssim  \frac{1}{r_0^{\frac{2}{\alpha-1}}} \Big( \frac{xw+1}{x(y+w)}+ 1-\mathbf{E}_y\big(\exp\big\{-\frac{C_*}{r_0^2} \tau_{w+y}^{(x^2),+} \big\}\big)
	 \Big),
\end{align}
where in the last inequality we used Lemma \ref{lemma11}(ii). Therefore, there exists $N$ depending on $r_0$ such that
\[
K^{(x)} (y+w)- K^{(x)}(y)\leq  N\left( \frac{xw+1}{x(y+w)}+ 1-\mathbf{E}_y\left(\exp\left\{-\frac{C_*}{r_0^2} \tau_{w+y}^{(x^2),+} \right\}\right)\right),
\]
this completes the proof of the lemma with $N_3(r_0)=\max\{N, \frac{C_*}{r_0^2}\}.$

\hfill$\Box$

\section{Proofs of the main results}\label{Main}

Throughout this section we assume  {\bf (H1)}, {\bf(H2)} and  {\bf(H4)} hold.

 \subsection{Proof of Theorem \ref{thm1}}\label{Section-3.1}

By Lemma \ref{lemma2}(i), for any $r, y>0$, we have
\begin{align}
	\sup_{t>0} v_\infty^{(t)}(r,y)<\infty.
\end{align}
By a standard diagonalization argument, for any sequence of positive reals increasing to $\infty$,  we can find a subsequence $\{t_k: k\in \N\}$ such that $\lim_{k\to +\infty} t_k= \infty$ and that the following limit exists
\begin{align}\label{Subsequence-limit}
	\lim_{k\to\infty} v_\infty^{(t_k)}(r,y) =: v_\infty^X(r,y),\quad \mbox{for all }
	r, y\in (0,\infty)\cap \Q.
\end{align}
Since $v_\infty^{(t)}(r,y)$ is decreasing in $r$ and increasing in $y$, so is the limit $v_\infty^X(r,y)$ for rational number $r$ and $y$. Therefore, for any
$r, y>0$, we can define
\begin{align}\label{Subsequence-rational-V-infty}
	v_\infty^X(r,y):=
	\lim_{(0,\infty)\cap \Q\ni (r_k, y_k) \to (r,y)} v_\infty^X(r_k, y_k)
	= \sup_{w\in [r, \infty)\cap\Q, z\in (0, y]\cap\Q}v_\infty^X(w,z).
\end{align}
    We define $v_\infty^X(r,0)=0$ for all $r>0$.

\begin{lemma}\label{lemma4}
	The relation \eqref{Subsequence-limit} holds for all $r, y>0$.
\end{lemma}
\textbf{Proof: }
For any $r, y>0$, let $\{(r_m, y_m): m\in \N\}$ be a sequence in $((0,\infty)\cap \Q)\times ((0,\infty)\cap \Q)$ with $(r_m, y_m)\to (r, y)$.
Note that
\begin{align}\label{step_14}
	& \left|v_\infty^X(r,y)- v_\infty^{(t_k)}(r,y)\right|\leq  \left|v_\infty^X(r,y)- v_\infty^X(r_m,y_m)\right| \nonumber\\
	&\quad + \left|v_\infty^X(r_m,y_m)- v_\infty^{(t_k)}(r_m,y_m)\right|+\left|v_\infty^{(t_k)}(r_m,y_m)- v_\infty^{(t_k)}(r,y)\right|.
\end{align}
Fix  $r_0\in \left(0, \left(\frac{1}{2}\inf_{m} r_m\right)  \land 1 \right) $, then there exists $A>0$ such that
$|y_m-y|< r_0$ and $|r_m- r|< r_0^4$ for all $m>A$.
By Lemma \ref{lemma3}, we have that
\begin{align}\label{step_15}
	&\left|v_\infty^{(t_k)}(r_m,y_m)- v_\infty^{(t_k)}(r,y)\right|\leq \left|v_\infty^{(t_k)}(r_m,y_m)- v_\infty^{(t_k)}(r_m,y)\right|  + \left|v_\infty^{(t_k)}(r_m,y)- v_\infty^{(t_k)}(r,y)\right|\nonumber\\
	& \leq N_1(r_0)\frac{1+\sqrt{t_k}|y_m-y|}{\sqrt{t_k|y_m-y|}}1_{\{y_m\neq y\}} +
	N_2(r_0) \frac{1+\sqrt{t_k}|r_m-r|^{1/4}}{\sqrt{t_k}|r_m-r|^{1/8}} 1_{\{r_m\neq r\}}.
\end{align}
Combining \eqref{Subsequence-limit},  \eqref{step_14} and \eqref{step_15},
\begin{align}\label{step_16}
	\limsup_{k\to\infty} \left|v_\infty^X(r,y)- v_\infty^{(t_k)}(r,y)\right| \lesssim \left|v_\infty^X(r,y)- v_\infty^X(r_m,y_m)\right|+ \sqrt{\left|y_m-y\right|}+\left|r_m-r\right|^{1/8}.
\end{align}
By \eqref{Subsequence-rational-V-infty}, letting $m\to\infty$ in \eqref{step_16}, we complete the proof of lemma.
\hfill$\Box$

\bigskip
The next lemma shows that any subsequential limit $v_\infty^X(r,y)$ is
a solution to some initial-boundary problem. We postpone its proof to  Section \ref{Appendix2}.

\begin{lemma}\label{lemma5}
	The limit $v_\infty^X(r,y)$ solves the	following initial-boundary value problem
\begin{equation}
	\left\{\begin{array}{rl}
		&\frac{\partial}{\partial r}v_\infty^X(r,y)= \frac{\sigma^2}{2}\frac{\partial^2}{\partial y^2} v_\infty^X(r,y) -
				                \varphi \left(v_\infty^X(r,y)\right), \quad
\mbox {in } (0,\infty)\times (0,\infty), \label{inner-eq}
\\
&\lim_{r\to 0+}v_\infty^X(r,y)=\infty, \quad y\in(0,\infty), \\
&\lim_{y\to 0+}v_\infty^X(r,y)=0, \quad r\in(0,\infty), \\
	\end{array}\right.
	\end{equation}
and for each $r\in (0,\infty)$, $\sup_{y>0} v_\infty^X(r,y)<\infty$.
\end{lemma}

 The next proposition is on the uniqueness of the solution to the problem \eqref{inner-eq}.

\begin{prop}\label{prop1}
The solution to  the problem \eqref{inner-eq} is unique and  can be written as
\begin{equation}\label{Explicit-formula-for-V-infty}
 v_\infty^X(r,y)= -\log \P_{\delta_y}(X_r^{(0,\infty)}=0),
 \end{equation}
  where $X_r^{(0,\infty)}$ is the process defined in Subsection \ref{Intro-SBM}.
\end{prop}
\textbf{Proof: }
Suppose that $u$ solves  problem \eqref{inner-eq}. For any $\delta>0$, $v(r,y):=u(\delta+r,y)$ solves the following problem:
\begin{equation*}
	\left\{\begin{array}{rl}
		&\frac{\partial}{\partial r}v(r,y)= \frac{\sigma^2}{2}\frac{\partial^2}{\partial y^2} v(r,y) -
				                \varphi \left(v(r,y)\right), \quad
\mbox {in } (0,\infty)\times (0,\infty),
\\
	&\lim_{r\to 0+}v(r,y)=u(\delta, y), \quad y\in(0,\infty), \\
&\lim_{y\to 0+} v(r,y)=0, \quad r\in(0,\infty), \\
	\end{array}\right.
	\end{equation*}
which is equivalent to the integral equation \eqref{Evolution-cumulant-semigroup} with $f= u(\delta, \cdot)$.
By the uniqueness of the solution to  \eqref{Evolution-cumulant-semigroup}, we get
\[
u (r+\delta, y)= v_{u(\delta,\cdot)}^X(r,y)=-\log\E_{\delta_y} \big(\exp\big\{-\langle u(\delta,\cdot), X_t^{(0,\infty)}\rangle \big\}\big), \quad r>0, y>0.
\]
Now letting $\delta\to 0+$ in the above equation, by
Lemma \ref{lemma5}
and the continuity of $v_\infty^X(r,y)$ in $r$,
\begin{align}
	&u(r,y)= \lim_{\delta \to0+} u (r+\delta, y)= \lim_{\delta \to0+}-\log\E_{\delta_y} \big(\exp\big\{-\langle u(\delta,\cdot), X_t^{(0,\infty)}\rangle \big\}\big)\nonumber\\
	& = -\log \P_{\delta_y} (X_r^{(0,\infty)}=0).
\end{align}
\hfill$\Box$

To prove Theorem \ref{thm1}, we  need some results on the uniform convergence of $ v_\infty^{(t)} (s,y)$ to $v_\infty^X(s,y)$ as $t\to\infty$.
For each $0<w<r$, taking $t=t_k$  in Lemma \ref{lemma3} and letting $k\to\infty$, we see that $v_\infty^X(s, y)$ is jointly
continuous for all $s\in [r-w, r]$ and $y>0$.
Since $v_\infty^{(t)}(s,y)$ is increasing in $y$ and $\sup_{t>1,y\in\R_+}
v_\infty^{(t)}(s,y)<\infty$, we see that $v_\infty^X(s,y)$ is also increasing in $y$ and that
$\sup_{y\in \R_+}
v_\infty^X(s,y)<\infty$, which implies the existence of $v_\infty^X(s,\infty):= \lim_{y\to \infty} v_\infty^X(s,y).$
Letting
$t\to\infty$ first and then $y\to\infty$ in Lemma \ref{lemma3} (ii), we see that $v_\infty^X(s, \infty)$ is continuous in $s\in [r-w, r].$
Therefore, for any $\varepsilon>0$, there exist $J, L\in \N$ and $s_0=r-w <s_1<...< s_{J}= r$, $y_0=0< y_1<...< y_{L} <
y_{L+1}=\infty$
such that
\begin{align}\label{step_66}
	\max_{j \in [1, J], \ell \in [0, L +1]} \left(\left| v_\infty^X (s_j, y_\ell) - v_\infty^X(s_{j-1}, y_\ell) \right|\vee \left| v_\infty^X (s_j, y_\ell) - v_\infty^X(s_{j}, y_{\ell-1}) \right| \right)<\varepsilon
\end{align}
and that there exists $T_2>0$ such that for any $t>T_2$,
\begin{align}
	\max_{j \in [0, J], \ell \in [0, L ]} \left| v_\infty^{(t)} (s_j, y_\ell) - v_\infty^X(s_{j}, y_\ell) \right| <\varepsilon.
\end{align}
Therefore,
for all $s\in [s_{j-1}, s_j]$ and $y\in [y_{\ell-1}, y_\ell)$ with $j\in [1,J], \ell\in [1, L+1]$, by the monotonicity of $v_\infty^{(t)}$, we get
\begin{align}\label{Fact2}
	v_\infty^{(t)}(s,y)\geq v_\infty^{(t)}(s_{j}, y_{\ell-1})\geq -\varepsilon + v_\infty^X(s_{j}, y_{\ell-1})\geq  -3\varepsilon + v_\infty^X(s_{j-1},
	y_{\ell})
	\geq -3\varepsilon + v_\infty^X(s,y).
\end{align}
Similarly, we also have
\begin{align}\label{Fact3}
	v_\infty^{(t)}(s,y)\leq 3\varepsilon + v_\infty^X(s,y).
\end{align}

Now we are ready to prove Theorem \ref{thm1}.

\vspace{.1in}

\noindent
\textbf{Proof of Theorem \ref{thm1}: } (i)
Combining Lemmas \ref{lemma4}, \ref{lemma5} and Proposition \ref{prop1}, we get that, for any $y>0$,
\begin{align}\label{step_18}
	\lim_{t\to+\infty} v_\infty^{(t)}(1,y)= \lim_{t\to\infty} t^{\frac{1}{\alpha-1}}\P_{\sqrt{t}y}\left(\zeta^{(0,\infty)}>t\right) = -\log \P_{\delta_y} (X_1^{(0,\infty)}=0).
\end{align}
Taking $f= \theta 1_{(0,\infty)}$ in \eqref{N-measure-equation} and letting $\theta\to+\infty$, we have that
\begin{align}\label{step_19}
	-\log \P_{\delta_y} (X_1^{(0,\infty)}=0) = \lim_{\theta \to+\infty} \mathbb N_y\left(1- \exp\left\{-\langle \theta1_{(0,\infty)}(\cdot), w_1 \rangle \right\}\right)=\N_y\left(w_1((0,\infty))\neq 0\right).
\end{align}
Combining \eqref{step_18} and \eqref{step_19}, we arrive at assertion (i).

(ii) By \eqref{Def-V-infty}, we have
\[
v_\infty^{(t)}\left(1,\frac{y}{\sqrt{t}}\right)= t^{\frac{1}{\alpha-1}}\P_{y}
    \big(\zeta^{(0,\infty)}>t\big),\quad y>0.
\]
It suffices to show that there exists $C^{(0,\infty)}(\alpha)\in (0,\infty)$ such that
\begin{align}\label{limt-v}
	\lim_{t\to\infty} \sqrt{t}v_\infty^{(t)}\left(1,\frac{y}{\sqrt{t}}\right) = R(y) C^{(0,\infty)}(\alpha).
\end{align}
By Proposition \ref{Feynman-Kac} and \eqref{step_6}, there exists $\Gamma>0$ such that for any $r_0\in (0, 1/2)$ and $w\in (0,r_0)$,
\begin{align}\label{step_20}
	 &\sqrt{t} \mathbf{E}_{y/\sqrt{t}}\big(v_\infty^{(t)} (1-w, \xi_w^{(t)}); \tau_0^{(t),-} >w \big)
	\nonumber \\  & \geq  \sqrt{t}v_\infty^{(t)}\left(1,\frac{y}{\sqrt{t}}\right)
	 \geq e^{-\frac{\Gamma}{r_0}w} \sqrt{t} \mathbf{E}_{y/\sqrt{t}}\big(v_\infty^{(t)} (1-w, \xi_w^{(t)}); \tau_0^{(t),-} >w \big).
\end{align}
By \eqref{Fact3},
for any $\varepsilon>0$, when $t$ is large enough,
\[
\sqrt{t} \mathbf{E}_{y/\sqrt{t}}\big(v_\infty^{(t)} (1-w, \xi_w^{(t)}); \tau_0^{(t),-} >w \big)
	 \leq \sqrt{t} \mathbf{E}_{y/\sqrt{t}}\big((v_\infty^{X}(1-w, \xi_w^{(t)}) +3\varepsilon ); \tau_0^{(t),-} >w \big),
\]
which, by Lemma \ref{Technical-lemma-1},
tends to
\[
\frac{1}{\sqrt{w}}\frac{2R(y)}{\sqrt{2\pi \sigma^2}}\int_0^\infty ze^{-\frac{z^2}{2}} (v_\infty^{X}(1-w, z\sigma \sqrt{w}) +3\varepsilon )\mathrm{d}z
\]	
as $t\to\infty$.
Similarly, using \eqref{Fact2}, we have
\begin{align}
& \sqrt{t} \mathbf{E}_{y/\sqrt{t}}\big(v_\infty^{(t)} (1-w, \xi_w^{(t)}); \tau_0^{(t),-} >w  \big)
	\geq \sqrt{t} \mathbf{E}_{y/\sqrt{t}}\big((v_\infty^{X}(1-w, \xi_w^{(t)}) -3\varepsilon ); \tau_0^{(t),-} >w  \big) \nonumber\\
	&\stackrel{t\to\infty}{\longrightarrow} \frac{1}{\sqrt{w}}\frac{2R(y)}{\sqrt{2\pi \sigma^2}}\int_0^\infty ze^{-\frac{z^2}{2}} (v_\infty^{X}(1-w, z\sigma \sqrt{w}) -3\varepsilon )\mathrm{d}z.
\end{align}
Therefore, letting $\varepsilon\to 0$, we conclude that
\begin{align}\label{step_68}
	\lim_{t\to\infty}& \sqrt{t} \mathbf{E}_{y/\sqrt{t}}\big(v_\infty^{(t)} (1-w, \xi_w^{(t)}); \tau_0^{(t),-} >w  \big)\nonumber\\
	& = \frac{1}{\sqrt{w}}\frac{2R(y)}{\sqrt{2\pi \sigma^2}}\int_0^\infty ze^{-\frac{z^2}{2}} v_\infty^{X}(1-w, z\sigma \sqrt{w}) \mathrm{d}z=: R(y) G(1-w, w).
\end{align}
Plugging this limit into \eqref{step_20}, we get that
\begin{align}\label{Fact4}
	&R(y) G(1-w, w) \geq \limsup_{t\to\infty}\sqrt{t}v_\infty^{(t)}\left(1,\frac{y}{\sqrt{t}}\right) \nonumber\\
	&\geq \liminf_{t\to\infty} \sqrt{t}v_\infty^{(t)}\left(1,\frac{y}{\sqrt{t}}\right)  \geq e^{-\frac{\Gamma}{r_0}w} R(y) G(1-w, w).
\end{align}
Using \eqref{Fact4}, for any $w\in (0, 1)$, we easily see that $G(1-w, w)\in(0,\infty)$, which implies that
\[
\infty> \limsup_{t\to\infty}\sqrt{t}v_\infty^{(t)}\left(1,\frac{y}{\sqrt{t}}\right) \geq \liminf_{t\to\infty} \sqrt{t}v_\infty^{(t)}\left(1,\frac{y}{\sqrt{t}}\right)  >0.
\]
Therefore, letting $w\to 0+$ in \eqref{Fact4}, we finally conclude that
\begin{align}
	\lim_{t\to\infty} \sqrt{t}v_\infty^{(t)}\left(1,\frac{y}{\sqrt{t}}\right) = R(y)\lim_{w\to 0+} G(1-w, w)=:  R(y)
		C^{(0,\infty)}(\alpha),
\end{align}
which is  \eqref{limt-v}. The proof is complete.

\hfill$\Box$

 \subsection{Proof of Theorem  \ref{thm2}}\label{Section-3.2}

To prove Theorem \ref{thm2}, we need to show the convergence of $v_f^{(t)}(r,y)$ for every continuous function $f\in B_b^+((0,\infty))$. The next lemma shows that we can assume additionally that $f$ is Lipschitz.

\begin{lemma}\label{Technical-lemma-2}
	Let $\mu_n$ and $\mu$ be
	non-negative finite
	random measures on $\R$, then the following conditions are equivalent:
	
	(i) For any continuous $f\in B_b^+(\R)$, $\int f(x)\mu_n(\mathrm{d}x) \stackrel{\mathrm{d}}{\Longrightarrow} \int f(x)\mu(\mathrm{d}x)$;
	
	(ii) For any Lipschitz continuous $f\in B_b^+(\R)$, $\int f(x)\mu_n(\mathrm{d}x) \stackrel{\mathrm{d}}{\Longrightarrow} \int f(x)\mu(\mathrm{d}x)$;
\end{lemma}

\textbf{Proof:}
We only need to prove (ii) $\Rightarrow$ (i). First note that (ii) implies that
$(\mu_n)_{n\geq 1}$ is relatively compact in distribution.
In fact, by \cite[Lemma 16.15]{Kallenberg},  $(\mu_n)_{n\geq 1}$ is tight iff any  relatively compact Borel set $B$, $\mu_n(B)$ is tight.
Taking $f\equiv 1$ in (ii), we see that $\mu_n(\R)$ is tight. For any  relatively compact Borel set $B$, using the fact that $\mu_n(B)\le \mu_n(\R))$, we get $\mu_n(B)$ is tight.
Now it remains to show that the distribution of a random measure $\mu$ is determined by $\int f(x)\mu(\mathrm{d}x)$ for all Lipschitz continuous $f\in B_b^+(\R)$,
which can be shown via a routine argument.
We omit the details.

\hfill$\Box$

The next two results will be needed in the proof of Theorem \ref{thm2}. We postpone their proofs to Section \ref{Proof: Technical-Prop-1}.

\begin{prop}\label{Technical-Prop-2}
	Suppose that $f$ is a bounded Lipschitz function on $\R_+$ with $f(0)=0$ and that $T>0$.
	
	(i) For any $r\in [0, T]$ and any $w>0$, it holds that
	\begin{align}
		\sup_{y>0} \big|v_f^{(t)}(r,y)- v_f^{(t)}(r, y+w)\big|\lesssim \left(\frac{1}{\log t} +  w \right)(1+r^{-1/2}).
	\end{align}
	
	(ii) For any $r,q \geq 0$ with $r+q\leq T$, it holds that
	\begin{align}
		\sup_{y>0} \big|v_f^{(t)}(r,y)- v_f^{(t)}(r+q, y)\big|\lesssim  \left(\frac{1}{\log t}+ q^{1/4}\right)\left(1+r^{-1/2}\right).
	\end{align}
\end{prop}

\begin{prop}\label{Technical-Prop-1}
	For any continuous function $f\in B_b^+((0,\infty))$ and any $r, y>0$, it holds that
	\[
	\lim_{t\to\infty} v_f^{(t)}(r,y)= v_f^X(r,y),
	\]
	where $v_f^X(r,y)$ is the solution of \eqref{Evolution-cumulant-semigroup}.
\end{prop}

Now we are ready to prove Theorem \ref{thm2}.

\vspace{.1in}

\noindent
\textbf{Proof of Theorem \ref{thm2}:}  (i)
By the definition of $Z_{1}^{(0,\infty), t}$ in \eqref{def-Z_1},  for any
continuous function $f\in B_b^+((0,\infty))$,
\begin{align}
  &\E_{\sqrt{t}y} \Big(\exp\Big\{- \int_{(0,\infty)} f\left(y\right) Z_{1}^{(0,\infty), t}(\mathrm{d}y)\Big\} \big| \zeta^{(0,\infty)}> t\Big)\nonumber\\
 =& \E_{\sqrt{t}y} \Big(\exp\Big\{- \frac{1}{t^{\frac{1}{\alpha-1}}} \int_{(0,\infty)} f\left(\frac{x}{\sqrt{t}}\right) Z_{t}^{(0,\infty)}(\mathrm{d}x)\Big\} \big| \zeta^{(0,\infty)}> t\Big).
 \end{align}
 Note that
\begin{align}\label{step_26}
	& \E_{\sqrt{t}y} \Big(\exp\Big\{- \frac{1}{t^{\frac{1}{\alpha-1}}} \int_{(0,\infty)} f\left(\frac{y}{\sqrt{t}}\right)Z_{t}^{(0,\infty)}(\mathrm{d}x)\Big\} \big| \zeta^{(0,\infty)}> t\Big)\nonumber\\
	&= 1- \frac{1}{\P_{\sqrt{t}y}(\zeta^{(0,\infty)}>t) } \E_{\sqrt{t}y} \Big( 1- \exp\Big\{- \frac{1}{t^{\frac{1}{\alpha-1}}} \int_{(0,\infty)} f\left(\frac{x}{\sqrt{t}}\right)Z_{t}^{(0,\infty)}(\mathrm{d}x)\Big\}  \Big),
\end{align}
where in the equality we used the fact that $Z_{t}^{(0,\infty)}((0,\infty))=0$ on the set $\{\zeta^{(0,\infty)}\leq t\}$.
    Recall the definitions of $v_f^{(t)}$ and $v_\infty^{(t)}$
in \eqref{Def-V-phi-psi} and \eqref{Def-V-infty},
\eqref{step_26} is equivalent to
\begin{align}\label{step_27}
	& \E_{\sqrt{t}y} \Big(\exp\Big\{- \frac{1}{t^{\frac{1}{\alpha-1}}} \int_{(0,\infty)} f\left(\frac{x}{\sqrt{t}}\right) Z_{t}^{(0,\infty)}(\mathrm{d}x)\Big\} \big| \zeta^{(0,\infty)}> t\Big) = 1- \frac{v_f^{(t)}(1,y)}{v_\infty^{(t)}(1,y)}.
\end{align}
Combining Proposition \ref{Technical-Prop-1} and \eqref{N-measure-equation}, we get
\begin{align}\label{step_28}
\lim_{t\to\infty} v_f^{(t)}(1,y)= v_f^X(1,y)= -\log \P_{\delta_y}\left(-\langle f, X_1 \rangle \right)= \N_y\left(1- \exp\left\{-\langle f, w_1\rangle\right\}\right).
\end{align}
Plugging \eqref{step_18}, \eqref{step_19} and \eqref{step_28} into \eqref{step_27}, we conclude that
\begin{align}
& \lim_{t\to\infty} \E_{\sqrt{t}y} \Big(\exp\Big\{- \frac{1}{t^{\frac{1}{\alpha-1}}} \int_{(0,\infty)} f\left(\frac{x}{\sqrt{t}}\right) Z_{t}^{(0,\infty)}(\mathrm{d}x)\Big\} \big| \zeta^{(0,\infty)}> t\Big)\nonumber\\
	& = 1- \frac{ \N_y\left(1- \exp\left\{-\langle f, w_1\rangle\right\}\right)}{\N_y\left(w_1((0,\infty))\neq 0\right)} = 1- \frac{ \N_y\left(1- \exp\left\{-\langle f, w_1|_{(0,\infty)}\rangle\right\}\right)}{\N_y\left(w_1((0,\infty))\neq 0\right)} \nonumber\\
	& = 1- \N_y\left(1- \exp\left\{-\langle f, w_1|_{(0,\infty)}\rangle\right\} \big| w_1((0,\infty))\neq 0 \right)\nonumber\\
	& = \N_y\left( \exp\left\{-\langle f, w_1|_{(0,\infty)}\rangle\right\} \big| w_1((0,\infty))\neq 0\right).
\end{align}
The completes the proof of (i).

(ii)
Let $f$ be an arbitrary non-negative bounded Lipschitiz function on $(0,\infty)$.
By \eqref{step_27}, we see that
\begin{align}\label{step_31}
	& \E_{y} \Big(\exp\Big\{- \frac{1}{t^{\frac{1}{\alpha-1}}} \int_{(0,\infty)} f\big(\frac{x}{\sqrt{t}}\big) Z_{t}^{(0,\infty)}(\mathrm{d}x)\Big\} \big| \zeta^{(0,\infty)}> t\Big) = 1- \frac{v_f^{(t)}\big(1,y t^{-\frac{1}{2}}\big)}{v_\infty^{(t)}\big(1,y t^{-\frac{1}{2}}\big)}.
\end{align}
By using an argument similar to that leading to \eqref{step_20}, we get that exists $\Gamma>0$ such that for any
$r_0\in (0, 1/2)$ and $w\in (0, r_0)$,
\begin{align}\label{step_67}
&\sqrt{t} \mathbf{E}_{y/\sqrt{t}}\big(v_f^{(t)} (1-w, \xi_w^{(t)}); \tau_0^{(t),-} >w \big)
	\nonumber \\  & \geq  \sqrt{t}v_f^{(t)}\big(1,\frac{y}{\sqrt{t}}\big)
	\geq e^{-\frac{\Gamma}{r_0}w} \sqrt{t} \mathbf{E}_{y/\sqrt{t}}\big(v_f^{(t)} (1-w, \xi_w^{(t)}); \tau_0^{(t),-} >w \big).
\end{align}
Proposition \ref{Technical-Prop-2} implies that, for any  $T>r_0$,
\begin{align}
\big|v_f^{(t)} (r,y)- v_f^{(t)}(s,z)\big|\lesssim \frac{1}{\log t}+ \left|y-z\right|+ \left| r-s\right|^{1/4}, \quad \mbox{ for all } r, s\in (r_0, T) \mbox{ and } y, z>0.
\end{align}
Therefore, for any  large $N>0$
and any $\varepsilon>0$, we can find $s_0= r_0<..<s_J=T$ and
$y_0=0<...<y_{L+1}= N$
 such that \eqref{step_66} holds, which in turn implies that \eqref{Fact2} and \eqref{Fact3} hold for all $s\in (r_0, T)$ and  $y\in (0,N)$  when $t$ is large enough.  Therefore, using an argument
 similar to that leading to \eqref{step_68}  and the following  consequence of Lemma \ref{Solve-Technical-lemma-1}
\begin{align}
	& \lim_{N\to\infty} \limsup_{t\to\infty} \sqrt{t} \mathbf{E}_{y/\sqrt{t}}\big(v_\infty^{(t)} (1-w, \xi_w^{(t)}); \tau_0^{(t),-} >w, \xi_w^{(t)}> N \big)\nonumber\\
	&\lesssim \lim_{N\to\infty} \limsup_{t\to\infty} \sqrt{t} \mathbf{P}_{y/\sqrt{t}}\big(\tau_0^{(t),-} >w,  \xi_w^{(t)}> N\big) = \lim_{N\to\infty}\frac{2}{\sqrt{2\pi \sigma^2}} R(y)\int_{N/\sigma} ze^{-\frac{z^2}{2}}\mathrm{d}z=0,
\end{align}
 we get the following result analogous to \eqref{step_68}:
 \begin{align}
 \lim_{t\to\infty}& \sqrt{t} \mathbf{E}_{y/\sqrt{t}}\big(v_f^{(t)} (1-w, \xi_w^{(t)}); \tau_0^{(t),-} >w \big)\nonumber\\
 	& = \frac{1}{\sqrt{w}}\frac{2R(y)}{\sqrt{2\pi \sigma^2}}\int_0^\infty ze^{-\frac{z^2}{2}} v_f^X(1-w, z\sigma \sqrt{w}) \mathrm{d}z.
 \end{align}
Therefore,
\begin{align}
	\lim_{t\to\infty} \sqrt{t}v_f^{(t)}\big(1,\frac{y}{\sqrt{t}}\big) = \lim_{w\to 0}\frac{1}{\sqrt{w}}\frac{2R(y)}{\sqrt{2\pi \sigma^2}}\int_0^\infty ze^{-\frac{z^2}{2}} v_f^X(1-w, z\sigma \sqrt{w}) \mathrm{d}z,
\end{align}
Together with Theorem \ref{thm1} (i), we conclude that
\begin{align}\label{step_69}
& \lim_{t\to\infty} \E_{y} \Big(\exp\Big\{- \frac{1}{t^{\frac{1}{\alpha-1}}} \int_{(0,\infty)} f\left(\frac{x}{\sqrt{t}}\right) Z_{t}^{(0,\infty)}(\mathrm{d}x)\Big\} \big| \zeta^{(0,\infty)}> t\Big) \nonumber\\
	&= 1- \frac{1}{C^{(0,\infty)}(\alpha)} \frac{2}{\sqrt{2\pi \sigma^2}} \lim_{w\to 0} \frac{1}{\sqrt{w}} \int_0^\infty ze^{-\frac{z^2}{2}}v_f^X(1-w, z\sigma \sqrt{w})\mathrm{d}z.
\end{align}
If we could show that
\begin{align}\label{e:t1.3}
	& \lim_{t\to\infty} \E_{y} \Big(\exp\Big\{-  \frac{1}{t^{\frac{1}{\alpha-1}}} \int_{(0,\infty)} \varepsilon f\left(\frac{y}{\sqrt{t}}\right) Z_{t}^{(0,\infty)}(\mathrm{d}y)\Big\} \big| \zeta^{(0,\infty)}> t\Big)\stackrel{\varepsilon \to 0+}{\longrightarrow}1,
\end{align}
we would get that there exists a random measure $\eta_1$ such that the right -hand  side of \eqref{step_69} is equal to
    $\E\big(\exp\big\{-\langle f, \eta_1\rangle\big\}\big)$.
Combining this with Lemma \ref{Technical-lemma-2}, we arrive at the assertion (ii).
Now we prove \eqref{e:t1.3}. By \eqref{Evolution-cumulant-semigroup},
\begin{align}
	v_{\varepsilon f}^X(1,y) \leq \varepsilon \mathbf{E}_y\left(f(W_1^0) \right) \leq \varepsilon \sup_{x>0} |f(x)| \mathbf{P}_y\big(\tau_0^{W,-} >1 \big)\lesssim \varepsilon y,
\end{align}
which implies that
\begin{align}
	\lim_{w\to0} \frac{1}{\sqrt{w}} \int_0^\infty ze^{-\frac{z^2}{2}}v_f^X(1-w, z\sigma \sqrt{w})\mathrm{d}z
	\lesssim 	 \varepsilon \int_0^\infty z^2e^{-\frac{z^2}{2}}  \mathrm{d}z \stackrel{\varepsilon \to 0+}{\longrightarrow}0.
\end{align}
Thus \eqref{e:t1.3} is valid.

\hfill$\Box$

\subsection{Proof of Theorem \ref{thm3}}

By Lemma \ref{lemma7},  for any $y\in (0,1)$,
\[
\sup_{x>0} K^{(x)}(y)<\infty.
\]
Therefore, using a diagonalization argument, we can find, for any sequence of positive reals increasing to $\infty$, a subsubsequence $\{x_k: k\in \N\}$ such that the following limit exists
\begin{align}\label{Subsequence-rational-K-infty}
	\lim_{k\to\infty} K^{(x_k)}(y) = : K^X(y),\quad \mbox{for all } y\in (0,1)\cap \Q.
\end{align}
Since $K^{(x)}(y)$ is monotone in $y\in (0,1)$, we can define
\[
K^X(y):=  \lim_{(0,1)\cap \Q\ni y_m\to y} K^X(y_m) = \inf_{z\in [y,1)\cap \Q} K^X(z), \quad y\in (0, 1).
\]
Using an argument similar to that used in the proof of  Lemma \ref{lemma4}, with Lemma \ref{lemma3} replaced by Lemma \ref{lemma8}, we can easily get
the following lemma, whose proof is omitted:

\begin{lemma}\label{lemma10}
	The relation  \eqref{Subsequence-rational-K-infty} holds for all $y\in (0,1)$.
\end{lemma}

The following Lemma \ref{lemma9} says that the limit $K^X(y)$ solves the boundary value problem \eqref{PDEin(0,1)} below.  Proposition \ref{prop2} is about the uniqueness and probabilistic representation to problem \eqref{PDEin(0,1)}.
Since the main idea of the proof of Lemma \ref{lemma9} is similar to that of Lemma \ref{lemma5}, and since we need to introduce exit measures of superprocesses in the proof of Proposition \ref{prop2},
we postpone the proofs to Section \ref{Appendix}.

\begin{lemma}\label{lemma9}
	The limit $K^X(y)$ solves the following problem
\begin{equation}\label{PDEin(0,1)}
	\left\{\begin{array}{rl}
		&\frac{\sigma^2}{2} (K^X)''(y) = \varphi(K^X(y)),\quad y\in (0,1),\\
&\lim_{y\to 0+}K^X(y)=0,\quad \lim_{y\to 1-} K^X(y)=\infty.
	\end{array}\right.
\end{equation}
\end{lemma}

\begin{prop}\label{prop2}
The problem in \eqref{PDEin(0,1)}
has a unique solution and the unique solution admits the representation $K^X(y)= -\log \P_{\delta_y}\left(M^{(0,\infty),X}<1\right)$.
\end{prop}

Now we are ready to prove Theorem \ref{thm3}.

\noindent
\textbf{Proof of Theorem \ref{thm3}: } Combining Lemma \ref{lemma10} and Proposition \ref{prop2}, we get that
\begin{align}
	\lim_{x\to\infty} K^{(x)}(y)= \lim_{x\to\infty} x^{\frac{2}{\alpha-1}}\P_{xy}(M^{(0,\infty)}\geq x)= -\log \P_{\delta_y}\big(M^{(0,\infty),X}<1\big),
\end{align}
which proves (i).
For (ii),
by the definition \eqref{def-K(x)} of $K^{(x)}(y)$ and Lemma \ref{lemma6}, for any
fixed small $z<\frac{1}{2}$, when $y< xz$, we have
\begin{align}
	& x^{\frac{2}{\alpha-1}+1} \P_y(M^{(0,\infty)}\geq x)= xK^{(x)}(yx^{-1})\nonumber\\
	&=x \mathbf{E}_{yx^{-1}}\Big(\exp\Big\{-\int_0^{\tau_z^{(x^2),+}} \psi^{(x^2)}\big( K^{(x)}\big(\xi_s^{(x^2)} \big)  \big)\mathrm{d}s \Big\}  K^{(x)}\big(\xi_{\tau_z^{(x^2),+}}^{(x^2)} \big); \tau_z^{(x^2),+}< \tau_0^{(x^2),-} \Big).
\end{align}
Therefore,
by Lemma \ref{lemma7}(ii) and the fact that
    $\xi_{\tau_z^{(x^2),+}}^{(x^2)} \geq z$, we have
\begin{align}\label{e:thrm1.3}
	& xK^{(x)} (z) \mathbf{E}_{yx^{-1}}\big(\exp\big\{-C_* \frac{1}{(1-z)^2} \tau_z^{(x^2),+} \big\}; \tau_z^{(x^2),+}< \tau_0^{(x^2),-}\big)\leq xK^{(x)}(yx^{-1})\nonumber\\
  & \leq  x\mathbf{E}_{yx^{-1}}\big( K^{(x)}\big(\xi_{\tau_z^{(x^2),+}}^{(x^2)} \big); \tau_z^{(x^2),+}< \tau_0^{(x^2),-} \big).
\end{align}
It follows from Lemma \ref{lemma12} that, for any $\delta>0$,
\begin{align}\label{Technical-step}
	 & x\mathbf{E}_{yx^{-1}}\big( K^{(x)}\big(\xi_{\tau_z^{(x^2),+}}^{(x^2)} \big); \tau_z^{(x^2),+}< \tau_0^{(x^2),-} \big)\nonumber\\
	&\leq K^{(x)}(z+\delta) x\mathbf{P}_{yx^{-1}}\big(  \tau_z^{(x^2),+}< \tau_0^{(x^2),-}\big)+ x^{1+\frac{2}{\alpha -1}}\mathbf{P}_{yx^{-1}}\big( \xi_{\tau_z^{(x^2),+}}^{(x^2)} > z+\delta\big) \nonumber\\
	&= K^{(x)}(z+\delta) x\mathbf{P}_{yx^{-1}}\big(  \tau_z^{(x^2),+}< \tau_0^{(x^2),-}\big)+x^{1+\frac{2}{\alpha -1}} \mathbf{P}_{-xz+y}\big(\xi_{\tau_0^+}> x\delta \big)\nonumber\\
	&\leq K^{(x)}(z+\delta) x\mathbf{P}_{yx^{-1}}\big(  \tau_z^{(x^2),+}< \tau_0^{(x^2),-}\big) +  \frac{1}{\delta^{r-2}} \frac{1}{x^{r_0-1-\frac{2\alpha}{\alpha -1}}}\sup_{y>0} \mathbf{E}_{-y}(|\xi_{\tau_0^+}|^{r_0-2}).
\end{align}
Therefore, combining \eqref{e:thrm1.3},  \eqref{Technical-step} and Lemma \ref{Technical-lemma-3}(i),  letting $x\to\infty$ first and then $\delta \to0$, we get that
\begin{align}\label{eq:5}
	\limsup_{x\to\infty}xK^{(x)}(yx^{-1}) \leq R(y)\frac{K^X(z)}{z}.
\end{align}
On the other hand,
by Lemma \ref{Technical-lemma-3} (ii), there exists a constant $C$ such that
$$
C_*x\mathbf{E}_{yx^{-1}} \big( \tau_{z}^{(x^2),+}; \tau_{z}^{(x^2),+}<\tau_0^{(x^2),-}  \big)
		\le  Cx z^2 \mathbf{P}_y\left( \tau_{xz}^+<\tau_0^-  \right)+  \frac{C}{x}.
$$
Thus,  by \eqref{e:thrm1.3}, using the inequality $e^{-x}\geq 1-x$ and Lemma \ref{Technical-lemma-3} (i), we have
\begin{align}\label{eq:6}
	 &\liminf_{x\to\infty} xK^{(x)}(yx^{-1})  \geq  K^X(z) \liminf_{x\to\infty} x\mathbf{E}_{yx^{-1}}\Big(\Big(1- C_* \frac{1}{(1-z)^2} \tau_z^{(x^2),+}\Big); \tau_z^{(x^2),+}< \tau_0^{(x^2),-}\Big) \nonumber\\
	&\geq
	K^X(z ) \Big(1- \frac{Cz^2}{(1-z)^2}\Big)\lim_{x\to\infty} x\mathbf{P}_{yx^{-1}}\big(\tau_z^{(x^2),+}< \tau_0^{(x^2),-}\big) \nonumber\\
& =R(y)\frac{K^X(z)}{z} \Big(1- \frac{Cz^2}{(1-z)^2}\Big).
\end{align}
Letting $z\to 0$, we conclude from \eqref{eq:5} and \eqref{eq:6} that
 $\lim_{z\to 0+} \frac{K^X(z)}{z}$
exists. Define
 $$\theta^{(0,\infty)}(\alpha):=\lim_{z\to 0+} \frac{K^X(z)}{z}.$$
 Then we have
\begin{align}\label{expression-theta}
	&\lim_{x\to\infty} xK^{(x)}(yx^{-1})  = 
\theta^{(0,\infty)}(\alpha)R(y).
\end{align}
Choose $z_0\in (0,1)$ such that $Cz_0^2/(1-z_0)^2 <1$. Then taking 
$z=z_0$
 in \eqref{eq:5} and \eqref{eq:6}, we get
\[
0< \frac{K^X(z_0)}{z_0} \Big(1- \frac{Cz_0^2}{(1-z_0)^2}\Big) \leq \theta^{(0,\infty)}(\alpha)\leq \frac{K^X(z_0)}{z_0}<\infty,
\]
which implies that $\theta^{(0,\infty)}(\alpha)\in(0,1)$. We complete the proof of the theorem.

\hfill$\Box$

\subsection{Proof of Theorem \ref{thm4}}

For $t, r, z>0$, define $Q_z^{(t)}(r,y):= t^{\frac{1}{\alpha-1}} \P_{\sqrt{t}y}\big(M_{tr}^{(0,\infty)}> \sqrt{t}z\big)$.

\begin{lemma}\label{e:newlemma}
Let  $z>0$
 and $\varepsilon \in (0, (z/2)\wedge 1)$. There exists a constant $L=L(\varepsilon)>0$ such that
for any $\delta>0$,
there exists $T=T(z, \varepsilon, \delta)$ such that when $t>T$,
\begin{align}\label{Upper-bound-of-Q}
	&Q_z^{(t)}(\delta,z-2\varepsilon) \leq L(\varepsilon)\delta.
\end{align}
\end{lemma}
\textbf{Proof: } Note that for $t, r, y, z>0$,
 \begin{align}
	\mathbb{P}_{\sqrt{t}y}\big(M_{tr}^{(0,\infty)}>\sqrt{t}z\big) = \lim_{\theta \to+\infty} \Big(1- \mathbb{E}_{\sqrt{t}y}\big(\exp\big\{- \theta Z_{tr}^{(0,\infty)}\left((\sqrt{t}z,\infty)\right) \big\}\big)\Big).
\end{align}
Taking $f= \theta 1_{(z,\infty)}$ in Proposition \ref{Feynman-Kac}
first and then letting $\theta \to +\infty$, we see that for any $w \in (0,  r]$, $Q_z^{(t)}(r,y)$ solves the equation
\begin{align}\label{eq:13}
	Q_z^{(t)}(r,y) = \mathbf{E}_y\Big(\exp\Big\{-\int_0^w \psi^{(t)}\big(Q_z^{(t)}(r-s, \xi_s^{(t)})\big)\mathrm{d}s \Big\}Q_z^{(t)} (r-w, \xi_w^{(t)}); \tau_0^{(t),-} >w \Big).
\end{align}
Taking $r=\delta, y= z-2\varepsilon$ and 
using the argument of getting \eqref{Feynman-Kac-1} with $T= \tau_{z-\varepsilon}^{(t),+}$, 
we get
\begin{align}\label{eq:17}
	&Q_z^{(t)}(\delta,z-2\varepsilon) \leq  \mathbf{E}_{z-2\varepsilon}\big(Q_z^{(t)} (\delta -\delta\land \tau_{z-\varepsilon}^{(t),+} ,  \xi_{\delta\land \tau_{z-\varepsilon}^{(t),+}}^{(t)}) \big)\nonumber\\
	& = \mathbf{E}_{z-2\varepsilon}\big(Q_z^{(t)} \big(\delta - \tau_{z-\varepsilon}^{(t),+} ,  \xi_{ \tau_{z-\varepsilon}^{(t),+}}^{(t)}\big); \delta > \tau_{z-\varepsilon}^{(t),+}  \big),
\end{align}
where in the equality we used the fact that $Q_z^{(t)}(0, y)=0$. By \eqref{Tail-probability-M}
and the fact that $M_{tr}^{(0,\infty)}\leq M$, we see that for any $0<y<z$ and $r>0$,
\[
Q_z^{(t)}(r,y) \leq t^{\frac{1}{\alpha-1}} \P_{\sqrt{t}y}\big(M> \sqrt{t}z\big)= t^{\frac{1}{\alpha-1}} \P_{0}\big(M> \sqrt{t}(z-y)\big)\lesssim\frac{1}{(z-y)^{\frac{2}{\alpha-1}}}.
\]
Combining the inequality above with the monotonicity of $Q_z^{(t)}(r,y)$ in $y$, we get that
\begin{align}\label{eq:16}
	& \mathbf{E}_{z-2\varepsilon}\big(Q_z^{(t)} \big(\delta - \tau_{z-\varepsilon}^{(t),+} ,  \xi_{ \tau_{z-\varepsilon}^{(t),+}}^{(t)}\big); \delta > \tau_{z-\varepsilon}^{(t),+} \big)\nonumber\\
	& = \mathbf{E}_{z-2\varepsilon}\big(Q_z^{(t)} \big(\delta - \tau_{z-\varepsilon}^{(t),+} ,  \xi_{ \tau_{z-\varepsilon}^{(t),+}}^{(t)}\big);\xi_{ \tau_{z-\varepsilon}^{(t),+}}^{(t)}> z -2^{-1}\varepsilon, \delta > \tau_{z-\varepsilon}^{(t),+} \big)\nonumber\\
	&\quad +  \mathbf{E}_{z-2\varepsilon}\big(Q_z^{(t)} \big(\delta - \tau_{z-\varepsilon}^{(t),+} ,  \xi_{ \tau_{z-\varepsilon}^{(t),+}}^{(t)}\big); \xi_{ \tau_{z-\varepsilon}^{(t),+}}^{(t)}\leq  z -2^{-1}\varepsilon, \delta > \tau_{z-\varepsilon}^{(t),+}\big)\nonumber\\
	& \leq t^{\frac{1}{\alpha-1}} \mathbf{P}_{z-2\varepsilon}\big(\xi_{ \tau_{z-\varepsilon}^{(t),+}}^{(t)}> z -2^{-1}\varepsilon \big)+ \mathbf{E}_{z-2\varepsilon}\big(Q_z^{(t)} \big(\delta - \tau_{z-\varepsilon}^{(t),+} ,  z -2^{-1}\varepsilon\big); \delta > \tau_{z-\varepsilon}^{(t),+} \big)\nonumber\\
	&\lesssim t^{\frac{1}{\alpha-1}} \mathbf{P}_{z-2\varepsilon}\big(\xi_{ \tau_{z-\varepsilon}^{(t),+}}^{(t)}> z -2^{-1}\varepsilon \big) + \frac{1}{\varepsilon^{\frac{2}{\alpha-1}}}\mathbf{P}_{z-2\varepsilon}\big(\delta> \tau_{z-\varepsilon}^{(t),+}\big).
\end{align}
Since $(r_0-2)/2>1/(\alpha-1)$, by Markov's inequality and Lemma \ref{lemma12}, we have
\begin{align}\label{eq:14}
	&t^{\frac{1}{\alpha-1}} \mathbf{P}_{z-2\varepsilon}\Big(\xi_{ \tau_{z-\varepsilon}^{(t),+}}^{(t)}> z -2^{-1}\varepsilon \Big) = t^{\frac{1}{\alpha-1}} \mathbf{P}_{-\sqrt{t}\varepsilon}\Big(\xi_{ \tau_{0}^{+}}> \frac{\varepsilon\sqrt{t}}{2} \Big)\nonumber\\
	    &\leq t^{\frac{1}{\alpha-1}} \left(\frac{2}{\varepsilon \sqrt{t}}\right)^{r_0-2} \mathbf{E}_{-\sqrt{t}\varepsilon}\left(\xi_{\tau_0^+}^{r_0-2}\right)\lesssim \frac{1}{\log t}\frac{1}{\varepsilon^{r_0-2}}.
\end{align}
By Doob's inequality,
\begin{align}\label{eq:15}
	\mathbf{P}_{z-2\varepsilon}\big(\delta> \tau_{z-\varepsilon}^{(t),+}\big)= 	
	     \mathbf{P}_{0}\big(\sup_{s\leq \delta}\xi_s^{(t)}\geq \varepsilon \big)\leq
	 \frac{\mathbf{E}_0\big(|\xi_\delta^{(t)}|^2\big)}{\varepsilon^2}\lesssim \frac{\delta}{\varepsilon^2}.
\end{align}
Plugging \eqref{eq:14} and \eqref{eq:15} into \eqref{eq:16}, we see that  for $t>e^{1/\delta}$, we have
\begin{align}\label{eq:18}
\mathbf{E}_{z-2\varepsilon}\Big(Q_z^{(t)} \Big(\delta - \tau_{z-\varepsilon}^{(t),+} ,  \xi_{ \tau_{z-\varepsilon}^{(t),+}}^{(t)}\Big); \delta > \tau_{z-\varepsilon}^{(t),+}  \Big)
  \lesssim \frac{\delta}{\varepsilon^{r_0-2}} + \frac{\delta }{\varepsilon^{\frac{2\alpha}{\alpha-1}}} \lesssim \frac{\delta}{\varepsilon^{r_0}}.
\end{align}
Combining \eqref{eq:17} and \eqref{eq:18}, we get the assertion of the lemma.

\hfill$\Box$

\noindent
\textbf{Proof of Theorem \ref{thm4}:}
It suffices to study the limits of the conditional probabilities of the
$\{M_t^{(0,\infty)}> \sqrt{t}z \}$ for $z>0$.

(i)
We first prove the lower bound. For any fixed $z>0$, define $g_1(x):= \min\{1, (x-z)^+\}$. then by Theorem \ref{thm2}(i),  for any $\theta>0$,
\begin{align}
	&\liminf_{t\to\infty} \P_{\sqrt{t} y}\big(M_{t}^{(0,\infty)}> \sqrt{t}z\big| \zeta^{(0,\infty)}> t\big)\nonumber\\
&=\liminf_{t\to\infty}\Big(1- \P_{\sqrt{t} y}\big(M_{t}^{(0,\infty)}\leq\sqrt{t}z\big| \zeta^{(0,\infty)}> t\big)\Big)\nonumber\\
	& \geq  \lim_{t\to\infty}\E_{\sqrt{t}y}\Big( 1- \exp\Big\{ - \frac{\theta }{t^{\frac{1}{\alpha-1}} } \int_{(0,\infty)} g_1\big(\frac{a}{\sqrt{t}}\big) Z_t^{(0,\infty)}(\mathrm{d}a) \Big\} \big| \zeta^{(0,\infty)}>t\Big)\nonumber\\
	& = \N_y\Big( 1-\exp\Big\{-\theta \int_{(0,\infty)} g_1(a)w_1(\mathrm{d}a)\Big\} \big| w_1((0,\infty))\neq 0\Big).
\end{align}
Letting $\theta\to +\infty$, we conclude that
\begin{align}\label{Lower-Bound-3}
	&\liminf_{t\to\infty} \P_{\sqrt{t}y}\big(M_{t}^{(0,\infty)}> \sqrt{t}z\big| \zeta^{(0,\infty)}> t\big)\nonumber\\
	&\geq \N_y\big( M_1^{(0,\infty), X}>z \big| w_1((0,\infty))\neq 0\big).
\end{align}
For the upper bound, we fix an arbitrary $z>0$. Let $\varepsilon\in (0, (z/2)\wedge 1)$ and $\delta\in (0, 1)$. 
We note that for any $r>\delta$,
\begin{align}
	&\P_{\sqrt{t} y} \big(M_{tr}^{(0,\infty)}> \sqrt{t}z\big) = 	\mathbb{E}_{\sqrt{t}y}	\Big(1-\exp\Big\{\int_{(0,\infty)} \log \P_a\big(M_{t\delta}^{(0,\infty)}\leq \sqrt{t}z\big) Z_{t(r-\delta)}^{(0,\infty)}(\mathrm{d}a) \Big\}\Big).
\end{align}
Note that for all $a>0$,
\[
\P_a\big(M_{t\delta}^{(0,\infty)}\leq \sqrt{t}z\big) \geq \P_a \big(\zeta^{(0,\infty)}\leq t\delta \big)\geq
\P(\zeta\le t\delta)
\stackrel{t\to\infty}{\longrightarrow}1.
\]
Using the fact that $\log x\sim x-1$ as $x\to 1$, we get that there exists $t_0=t_0(\delta)>0$ such that for all $t\ge t_0$,
\begin{align}\label{eq:21}
	&\P_{\sqrt{t} y} \big(M_{tr}^{(0,\infty)}> \sqrt{t}z\big) \leq	\E_{\sqrt{t}y}	\Big(1-\exp\Big\{-\frac{1}{2}\int_{(0,\infty)}  \P_a\big(M_{t\delta}^{(0,\infty)}> \sqrt{t}z\big) Z_{t(r-\delta)}^{(0,\infty)}(\mathrm{d}a) \Big\}\Big).
\end{align}
    When $a<\sqrt{t}(z-\varepsilon)$, by Lemma \ref{e:newlemma} with $\varepsilon$ replaced by $\frac{\varepsilon}{2}$ and using the monotonicity of $Q_z^{(t)}$,
we get that when $t\geq T(z,\frac{\varepsilon}{2},\delta)$,
\begin{align}\label{eq:19}
	t^{\frac{1}{\alpha-1}} \P_a\big(M_{t\delta}^{(0,\infty)}> \sqrt{t}z\big) \leq 	t^{\frac{1}{\alpha-1}} \P_{\sqrt{t}(z-\varepsilon)}\big(M_{t\delta}^{(0,\infty)}> \sqrt{t}z\big) = Q_z^{(t)}(\delta, z-\varepsilon)\leq L\left(\frac{\varepsilon}{2}\right)\delta.
\end{align}
When $a\geq \sqrt{t}(z-\varepsilon)$, by \eqref{Survival-prob-zeta}, there exists a constant $L_1$ such that
\begin{align}\label{eq:20}
	t^{\frac{1}{\alpha-1}} \P_a\big(M_{t\delta}^{(0,\infty)}> \sqrt{t}z\big) \leq t^{\frac{1}{\alpha-1}} \P\big(\zeta > \sqrt{t}\delta\big)\leq \frac{L_1}{\delta^{\frac{1}{\alpha-1}}}.
\end{align}
For any fixed $\varepsilon\in (0, (z/2)\wedge 1)$, let
$\delta_*>0$
 be small enough so that
\[
L\left(\frac{\varepsilon}{2}\right)\delta_*< \frac{L_1}{\delta_*^{\frac{1}{\alpha-1}}}.
\]
 Define another  non-negative bounded continuous function
\begin{align}
 &g_2(x):= L\left(\frac{\varepsilon}{2}\right)\delta_* 1_{\{x\leq z-2\varepsilon\}} +  \frac{L_1}{\delta_*^{\frac{1}{\alpha-1}}} 1_{\{x\geq z-\varepsilon\}}  \nonumber\\ &\quad\quad+ \bigg(\bigg(\frac{L_1}{\delta_*^{\frac{1}{\alpha-1}}}- L\left(\frac{\varepsilon}{2}\right)\delta_*\bigg) \frac{x-(z-2\varepsilon)}{\varepsilon}+ L\left(\frac{\varepsilon}{2}\right)\delta_* \bigg)1_{\{x\in (z-2\varepsilon, z-\varepsilon)\}}.
\end{align}
Then \eqref{eq:19} and \eqref{eq:20} imply that for all $a\in (0,\infty)$,
\begin{align}\label{Upper-bound-3}
	t^{\frac{1}{\alpha-1}} \P_a\big(M_{t\delta_*}^{(0,\infty)}> \sqrt{t}z\big) \leq g_2(a).
\end{align}
Plugging this upper bound into \eqref{eq:21}, we conclude that
\begin{align}\label{Upper-bound-4}
	&\P_{\sqrt{t} y} \big(M_{tr}^{(0,\infty)}> \sqrt{t}z\big) \leq \mathbb{E}_{\sqrt{t}y} \Big(1-\exp\Big\{- \frac{1}{t^{\frac{1}{\alpha-1}}}\int_{(0,\infty)} g_2\big(\frac{a}{\sqrt{t}} \big) Z_{t(r-\delta_*)}^{(0,\infty)}(\mathrm{d}a) \Big\}\Big)\nonumber\\
	& = t^{-\frac{1}{\alpha-1}}v_{g_2}^{(t)}(r-\delta_*, y),
\quad r>\delta_*.
\end{align}
Combining Proposition \ref{Technical-Prop-1} (applied to $g_2$)  and
 Theorem \ref{thm1},
we get that
\begin{align}\label{Upper-bound-1}
	&\limsup_{t\to\infty}\P_{\sqrt{t} y} \big(M_{t}^{(0,\infty)}> \sqrt{t}z\big| \zeta^{(0,\infty)}> t\big) = \limsup_{t\to\infty} \frac{\P_{\sqrt{t} y} \big(M_{t}^{(0,\infty)}> \sqrt{t}z\big) }{\P_{\sqrt{t} y} \big(\zeta^{(0,\infty)}> t\big) }  \nonumber\\
	& = \frac{ \limsup_{t\to\infty} (t(1+\delta_*))^{\frac{1}{\alpha-1}}\P_{\sqrt{t} \sqrt{1+\delta_*}y} \big(M_{t(1+\delta_*)}^{(0,\infty)}> \sqrt{t(1+\delta_*)}z\big) }{  \lim_{t\to\infty}t^{\frac{1}{\alpha-1}}\P_{\sqrt{t} y} \big(\zeta^{(0,\infty)}>t\big) } \nonumber\\
	& \leq  (1+\delta_*)^{\frac{1}{\alpha-1}}\frac{ \limsup_{t\to\infty} t^{\frac{1}{\alpha-1}}\P_{\sqrt{t} \sqrt{1+\delta_*}y} \big(M_{t(1+\delta_*)}^{(0,\infty)}> \sqrt{t}z\big) }{  \N_y\big(w_1((0,\infty))\neq 0\big)} \nonumber\\
 &\leq (1+\delta_*)^{\frac{1}{\alpha-1}}\frac{ \lim_{t\to\infty} v_{g_2}^{(t)}(1, y\sqrt{1+\delta_*})}{  \N_y\left(w_1((0,\infty))\neq 0\right)}
	\nonumber\\	&
	= (1+\delta_*)^{\frac{1}{\alpha-1}}\frac{-\log \E_{\delta_{\sqrt{1+\delta_*}y}} \big(\exp\big\{- \int_{(0,\infty)} g_2(a)X_{1}(\mathrm{d}a)\big\}\big)}{  \N_y\big(w_1((0,\infty))\neq 0\big)},
\end{align}
where in the second inequality we used \eqref{Upper-bound-4} with $r=1+ \delta_*$.
By the inequality $e^{-x}\geq 1-x$ and the fact that $\{X_1((z,\infty))=0\}= \{X_1^{(0,\infty)}((z,\infty))=0\}=\{M_1^{(0,\infty)}\leq z\}$, we have
\begin{align}
	& \E_{\delta_{\sqrt{1+\delta_*}y}} \Big(\exp\Big\{- \int_{(0,\infty)} g_2(a)X_{1}(\mathrm{d}a)\Big\}\Big)\nonumber\\	& \geq  \E_{\delta_{\sqrt{1+\delta_*}y}} \left(\exp\left\{- L\left(\frac{\varepsilon}{2}\right)\delta_* X_1 ((-\infty, z-2\varepsilon])\right\}; X_1((z-2\varepsilon, \infty))=0\right)\nonumber\\	& \geq  \E_{\delta_{\sqrt{1+\delta_*}y}} \Big(\Big(1-L\left(\frac{\varepsilon}{2}\right)\delta_* X_1 ((-\infty, z-2\varepsilon])\Big) ; X_1((z-2\varepsilon, \infty))=0 \Big)\nonumber\\	& \geq \P_{\delta_{\sqrt{1+\delta_*}y}} \big(M_1^{(0,\infty)} \leq z-2\varepsilon\big) - L\left(\frac{\varepsilon}{2}\right)\delta_*\E_{\delta_{\sqrt{1+\delta_*}y}} X_1(\R) .
\end{align}
Note that $X_1(\R)$ under $\P_{\delta_{\sqrt{1+\delta_*}y}}$ is a critical continuous-state branching process, the last term in the above inequality is equal to $-L\left(\frac{\varepsilon}{2}\right)\delta_*$. Therefore, letting $\delta_*\to 0$, we conclude that
\begin{align}
	& \limsup_{\delta_*\to0}-\log \E_{\delta_{\sqrt{1+\delta_*}y}} \Big(\exp\Big\{- \int_{(0,\infty)} g_2(a)X_{1}(\mathrm{d}a)Big\}\Big)\nonumber\\
	&\leq -\log \P_{\delta_{y}} \big(M_1^{(0,\infty)} \leq z-2\varepsilon\big).
\end{align}
Now letting $\varepsilon\to0$, we conclude that
\begin{align}\label{Upper-bound-2}
& \limsup_{\varepsilon\to 0}\limsup_{\delta_* \to 0} (1+\delta_*)^{\frac{1}{\alpha-1}}\frac{-\log \E_{\delta_{\sqrt{1+\delta_*}y}} \Big(\exp\Big\{- \int_{(0,\infty)} g_2(a)X_{1}(\mathrm{d}a)\Big\}\Big)}{  \N_y\left(w_1((0,\infty))\neq 0\right)} \nonumber\\& \leq  \frac{-\log \P_{\delta_{y}} \big(M_1^{(0,\infty)} \leq z\big) }{\N_y\left(w_1((0,\infty))\neq 0\right) }
= \frac{\N_y\big(M_1^{(0,\infty)}>z\big)}{\N_y\left(w_1((0,\infty))\neq 0\right)}\nonumber\\	& = \N_y\left(M_1^{(0,\infty)}>z\big| w_1((0,\infty))\neq 0 \right).
\end{align}
Combining \eqref{Lower-Bound-3}, \eqref{Upper-bound-1} and \eqref{Upper-bound-2},
we get the assertion of (i).

(ii) The proof of (ii) is similar.
In \eqref{Lower-Bound-3}, by replacing $\sqrt{t}y$ by $y$ and applying Theorem \ref{thm2}(ii), we get that
\begin{align}
	&\liminf_{t\to\infty} \P_{y}\big(M_{\sqrt{t}}^{(0,\infty)}> \sqrt{t}z\big| \zeta^{(0,\infty)}> t\big)\nonumber\\
	&\geq \lim_{\theta \to\infty} \P\Big( 1-\exp\Big\{-\theta \int_{(0,\infty)} g_1(a)\eta_1(\mathrm{d}a)\Big\} \Big) = \P(M^{\eta_1}>z).
\end{align}
For the upper bound, the argument in \eqref{Upper-bound-3} still holds in this case. Therefore,
by \eqref{Upper-bound-4} with $r=1+\delta_*$,
we get that
\begin{align}\label{Upper-bound-5}
	&t^{\frac{1}{\alpha-1}+\frac{1}{2}}\P_{ y} \big(M_{t(1+\delta_*)}^{(0,\infty)}> \sqrt{t}z\big) \nonumber\\
	& \leq t^{\frac{1}{\alpha-1}+\frac{1}{2}}\mathbf{E}_{y} \Big(1-\exp\Big\{- \frac{1}{t^{\frac{1}{\alpha-1}}}\int_{(0,\infty)} g_2\left(\frac{a}{\sqrt{t}} \right) Z_{t}^{(0,\infty)}(\mathrm{d}a) \Big\}\Big).
\end{align}
Combining \eqref{Upper-bound-5}, Theorem \ref{thm1}(i) and Theorem \ref{thm2}(ii), we see that
\begin{align}
	& \limsup_{t\to\infty} t^{\frac{1}{\alpha-1}+\frac{1}{2}}\P_{ y} \big(M_{t}^{(0,\infty)}> \sqrt{t}z\big) \nonumber\\
	&=  (1+\delta_*)^{ \frac{1}{\alpha-1}+\frac{1}{2}}\limsup_{t\to\infty} t^{\frac{1}{\alpha-1}+\frac{1}{2}}\P_{ y} \big(M_{t(1+\delta_*)}^{(0,\infty)}> \sqrt{t(1+\delta_*)}z\big)\nonumber\\
	&\leq  (1+\delta_*)^{ \frac{1}{\alpha-1}+\frac{1}{2}} \limsup_{t\to\infty} t^{\frac{1}{\alpha-1}+\frac{1}{2}}\P_{ y} \left(M_{t(1+\delta_*)}^{(0,\infty)}> \sqrt{t}z\right) \nonumber\\
	&\leq
	(1+\delta_*)^{ \frac{1}{\alpha-1}+\frac{1}{2}}
C^{(0,\infty)}(\alpha)
	 R(y) \mathbb{E}\Big(1-\exp\Big\{-\int_{(0,\infty)} g_2(a)\eta_1(\mathrm{d}a)\Big\}\Big).
\end{align}
Now letting $\delta_*\to 0$ first and then $\varepsilon\to 0$, we see that
\begin{align}
	& \limsup_{t\to\infty} \P_{ y} \big(M_{t}^{(0,\infty)}> \sqrt{t}z \big| \zeta^{(0,\infty)} >t\big) \nonumber\\
	& \leq \frac{
C^{(0,\infty)}(\alpha)
R(y) \P(M^{\eta_1}>z)}{\lim_{t\to\infty}t^{\frac{1}{\alpha-1}+\frac{1}{2}}  \P_{ y} \left(\zeta^{(0,\infty)} >t\right)   } = \P(M^{\eta_1}>z),
\end{align}
which completes the proof of (ii).

\hfill$\Box$

\section{Proofs of the auxiliary results}\label{aux}

\subsection{Proof of Lemma \ref{Technical-lemma-1}}\label{Proof: Technical-lemma-1}

In this subsection,
we assume that the L\'{e}vy process $\xi$ satisfies {\bf(H2)} and \eqref{eq:7}.
We use
\[
\Phi^+(t):=\int_0^t ze^{-\frac{z^2}{2}}\mathrm{d}z, \quad t\ge 0
\]
to denote the  Rayleigh distribution function.

\begin{lemma}\label{Solve-Technical-lemma-1}
	For any $y>0$ and $a\in (0,\infty]$,  it holds that
	\begin{align}
		\lim_{t\to\infty} \sqrt{t} \mathbf{P}_y\big(\xi_t\leq a \sqrt{t} , \tau_0^->t\big)=  \frac{2}{\sqrt{2\pi \sigma^2}} R(y)\Phi^+\left(\frac{a}{\sigma}\right).
	\end{align}
\end{lemma}
\textbf{Proof: } For any $r>0$ and
    $\varepsilon\in (0, \delta/(2(5+2\delta)))$,
where $\delta$ is the constant in \eqref{eq:7}, we define
\begin{align}
	A_r:= \Big\{\sup_{0\leq s\leq 1} \left|\xi_{sr} -W_{sr} \right| \leq r^{\frac{1}{2}-2\varepsilon} \Big\}.
\end{align}
Recall that the random walk $S_n$ is given by $S_n =\xi_n$. For any $b\in \R$, define
\[
\tau_b^{S,+}:=\inf\left\{j\in \N: S_j > b\right\}.
\]
Then we have the following decomposition:
\begin{align}
	 \sqrt{t} \mathbf{P}_y\left(\xi_t\leq a \sqrt{t} , \tau_0^->t\right)= \sum_{k=1}^4 I_k,
\end{align}
where $I_k$ are defined by
\begin{align}
	I_1 &:= \sqrt{t}\mathbf{P}_y \Big( \xi_t\leq a \sqrt{t} , \tau_0^->t, \tau_{t^{1/2-\varepsilon}}^{S,+}> [t^{1-\varepsilon}]\Big),\nonumber\\
	I_2 &:=\sqrt{t} \sum_{k=1}^{[t^{1-\varepsilon}]} \mathbf{E}_y\big(\mathbf{P}_{\xi_k}\big(\xi_{t-k}\leq a\sqrt{t}, \tau_0^->t-k, A_{t-k}^c\big); \tau_0^->k, \tau_{t^{1/2-\varepsilon}}^{S,+} =k \big),\nonumber\\
	I_3&:= \sqrt{t} \sum_{k=1}^{[t^{1-\varepsilon}]} \mathbf{E}_y\big(\mathbf{P}_{\xi_k}\big(\xi_{t-k}\leq a\sqrt{t}, \tau_0^->t-k, A_{t-k}\big); \tau_0^->k, \xi_k> t^{(1-\varepsilon)/2},  \tau_{t^{1/2-\varepsilon}}^{S,+} =k  \big),\nonumber\\
	I_4 &:= \sqrt{t} \sum_{k=1}^{[t^{1-\varepsilon}]} \mathbf{E}_y\big(\mathbf{P}_{\xi_k}\big(\xi_{t-k}\leq a\sqrt{t}, \tau_0^->t-k, A_{t-k}\big); \tau_0^->k, \xi_k\leq  t^{(1-\varepsilon)/2},  \tau_{t^{1/2-\varepsilon}}^{S,+} =k \big).
\end{align}

(i) In this part, we show that $\lim_{t\to\infty} I_1=0$.
Since $I_1\leq \sqrt{t}\mathbf{P}_y\big( \tau_{t^{1/2-\varepsilon}}^{S,+}> [t^{1-\varepsilon}]\big)$, it suffices to prove that
\begin{align}\label{Goal-1}
	\lim_{t\to\infty} \sqrt{t}\mathbf{P}_y\big( \tau_{t^{1/2-\varepsilon}}^{S,+}> [t^{1-\varepsilon}]\big)=0.
\end{align}
Since $[t^{1-\varepsilon}]\geq [t^\varepsilon -1] [t^{1-2\varepsilon}]=: L_1 \cdot L_2,$ we have
\begin{align}\label{step_61}
	\mathbf{P}_y\Big( \tau_{t^{1/2-\varepsilon}}^{S,+}> [t^{1-\varepsilon}]\Big)= \mathbf{P}_y\Big( \max_{j\leq [t^{1-\varepsilon}]} |S_j|\leq t^{1/2-\varepsilon} \Big) \leq  \mathbf{P}_y\Big( \max_{j\leq L_1} |S_{L_2j}|\leq t^{1/2-\varepsilon} \Big).
\end{align}
Applying the Markov property repeatedly, we get that
\begin{align}\label{step_60}
	&\mathbf{P}_y\Big( \max_{j\leq L_1} |S_{L_2j}|\leq t^{1/2-\varepsilon} \Big)\leq \sup_{x\in \R} \mathbf{P}_x\Big(|S_{L_2}| \leq ^{1/2-\varepsilon} \Big) \mathbf{P}_y\Big( \max_{j\leq L_1-1} |S_{L_2j}|\leq t^{1/2-\varepsilon} \Big)\nonumber\\
	&\leq\cdots\leq \left(\sup_{x\in \R} \mathbf{P}_x\Big(|S_{L_2}| \leq t^{1/2-\varepsilon} \Big) \right)^{L_1}.
\end{align}
The classical central limit theorem implies that when $x>2\sqrt{L_2}$
\begin{align}\label{CLT-1}
	& \mathbf{P}_x\big(|S_{L_2}| \leq  t^{1/2-\varepsilon} \big) \leq \mathbf{P}_x\big(|S_{L_2}| \leq  2\sqrt{L_2}\big) \nonumber\\
	&=\mathbf{P}_0\big( -x-2\sqrt{L_2}\leq S_{L_2} \leq 2\sqrt{L_2}-x \big)\leq \mathbf{P}_0\left(  S_{L_2}\leq 0 \right) \stackrel{t\to\infty}{\longrightarrow}
	\frac{1}{2}.
\end{align}
Similarly, when $x<-2\sqrt{L_2}$, we have
\begin{align}\label{CLT-2}
	& \mathbf{P}_x\big(|S_{L_2}| \leq  t^{1/2-\varepsilon} \big) \leq \mathbf{P}_0\left(  S_{L_2}\geq 0 \right) \stackrel{t\to\infty}{\longrightarrow}\frac{1}{2}
\end{align}
and that for $|x|\leq 2\sqrt{L_2}$,
\begin{align}\label{CLT-3}
	\mathbf{P}_x\big(|S_{L_2}| \leq  t^{1/2-\varepsilon} \big) \leq \mathbf{P}_0\big(  -4\sqrt{L_2}\leq S_{L_2}
	    \leq 4\sqrt{L_2} \big)
	\stackrel{t\to\infty}{\longrightarrow} \int_{-4}^4 \frac{1}{\sqrt{2\pi \sigma^2}}e^{-\frac{z^2}{2\sigma^2}}\mathrm{d}z.
\end{align}
Combining \eqref{CLT-1}, \eqref{CLT-2} and \eqref{CLT-3}, we see that
there exist $c\in (0,1)$ and $t_0>0$ such that
\[
\mathbf{P}_x\big(|S_{L_2}| \leq  t^{1/2-\varepsilon} \big) <c, \quad x\in \R, t\ge t_0.
\]
Plugging this into \eqref{step_60} and combining the conclusion with \eqref{step_61}, we get \eqref{Goal-1}.

(ii)
In this part, we show that $\lim_{t\to\infty} I_2=0$.
By Lemma \ref{Coupling-Levy-BM} and the definition of $A_r$, we have
\begin{align}\label{Upp-bound-of-I-2}
	I_2& \leq \sqrt{t} \sum_{k=1}^{[t^{1-\varepsilon}]} \mathbf{E}_y\Big(\mathbf{P}_{\xi_k}\left( A_{t-k}^c\right); \tau_0^->k,  \tau_{t^{1/2-\varepsilon}}^{S,+} =k \Big), \nonumber\\
	&\lesssim \sqrt{t}\sum_{k=1}^{[t^{1-\varepsilon}]}
	    \frac{N_*(2\varepsilon)}{(t-k)^{(2+\delta)(\frac{1}{2}-2\varepsilon)-1}}
	\mathbf{P}_y\big(\tau_0^->k, \tau_{t^{1/2-\varepsilon}}^{S,+} =k \big) \nonumber\\
	 & \lesssim \frac{  \sqrt{t}N_*(2\varepsilon)}{t^{(2+\delta)(\frac{1}{2}-2\varepsilon)-1}}\sum_{k=1}^{[t^{1-\varepsilon}]}  \mathbf{P}_y\big(\tau_0^{S,-}>k, \tau_{t^{1/2-\varepsilon}}^{S,+} =k \big).
\end{align}
Since $S_k$ is also a martingale under $\mathbf{P}_y$, using the fact that $S_k\geq t^{1/2 -\varepsilon}$ on the event $\{\tau_{t^{1/2-\varepsilon}}^{S,+} =k\}$, by \eqref{Upp-bound-of-I-2},
\begin{align}
	& I_2 \lesssim \frac{  \sqrt{t}N_*(2\varepsilon)}{t^{(2+\delta)(\frac{1}{2}-2\varepsilon)-1} t^{1/2-\varepsilon}}\sum_{k=1}^{[t^{1-\varepsilon}]}  \mathbf{P}_y\big( S_k;  \tau_0^{S,-}>k, \tau_{t^{1/2-\varepsilon}}^{S,+} =k \big)\nonumber\\
	& \leq \frac{  \sqrt{t}N_*(2\varepsilon)}{t^{(2+\delta)(\frac{1}{2}-2\varepsilon)-1} t^{1/2-\varepsilon}}\mathbf{E}_y\big( S_{[t^{1-\varepsilon}]\land \tau_0^{S,-}\land \tau_{t^{1/2}-\varepsilon}^{S,+}} \big) = \frac{N_*(2\varepsilon) y}{t^{\delta/2 - (5+2\delta)\varepsilon}},
\end{align}
we see that when $\varepsilon$ is sufficiently small  so that
$\varepsilon<\delta /(2(5+2\delta))$,
we have $\lim_{t\to\infty} I_2=0$.

(iii) In this part, we show that $\lim_{t\to\infty} I_3=0$.
Set $x'= \xi_k> t^{(1-\varepsilon)/2}$, by Lemma \ref{lemma11}(i) with $t=1$, we have that
\begin{align}
	&\mathbf{P}_{x'}\big(\xi_{t-k}\leq a\sqrt{t}, \tau_0^->t-k, A_{t-k}\big) \leq	\mathbf{P}_{x'}\left( \tau_0^->t-k\right)\lesssim \frac{x'+1}{\sqrt{t-k}}\lesssim\frac{x'}{\sqrt{t}},
\end{align}
where in the last inequality we used the fact that $t\lesssim t-k$ for all $k\leq [t^{1-\varepsilon}]$.
Therefore, we have
\begin{align}\label{eq:12}
	& I_3\lesssim  \sum_{k=1}^{[t^{1-\varepsilon}]} \mathbf{E}_y\big(\xi_k; \tau_0^->k, \xi_k> t^{(1-\varepsilon)/2}, \tau_{t^{1/2-\varepsilon}}^{S,+} =k  \big)\nonumber\\
	& \leq  \sum_{k=1}^{[t^{1-\varepsilon}]} \mathbf{E}_y\big(S_k; \tau_0^{S,-}>k, S_k> t^{(1-\varepsilon)/2},  \tau_{t^{1/2-\varepsilon}}^{S,+} =k \big).
\end{align}
Now we deal with the random walk $S_k$. Set $\Delta_k:= S_k- S_{k-1}$. We have that
\begin{align}\label{eq:10}
	  &\sum_{k=1}^{[t^{1-\varepsilon}]} \mathbf{E}_y\big(S_k; \tau_0^{S,-}>k, S_k> t^{(1-\varepsilon)/2},  \tau_{t^{1/2-\varepsilon}}^{S,+} =k  \big)\nonumber\\
	  &\leq \sum_{k=1}^{[t^{1-\varepsilon}]} \mathbf{E}_y\big((S_{k-1}+\Delta_k); \tau_0^{S,-}>k-1, S_k> t^{(1-\varepsilon)/2},  \Delta_k > t^{(1-\varepsilon)/2} - t^{1/2-\varepsilon} \big) \nonumber\\
	  & \leq  \sum_{k=1}^{[t^{1-\varepsilon}]} \mathbf{E}_y\big(S_{k-1} ; \tau_0^{S,-}>k-1 \big) \mathbf{P}_y\big(  \Delta_k > t^{(1-\varepsilon)/2} - t^{1/2-\varepsilon}  \big)\nonumber\\
	     & \quad + \sum_{k=1}^{[t^{1-\varepsilon}]} \mathbf{P}_y\big(\tau_0^{S,-}>k-1 \big) \mathbf{E}_y\big(\Delta_k ;  \Delta_k > t^{(1-\varepsilon)/2} - t^{1/2-\varepsilon}  \big).
\end{align}
Noting  that, for any fixed $y$, $\mathbf{E}_y\big(S_{k-1}; \tau_0^{S,-}>k-1 \big)\lesssim1$,
$\mathbf{P}_y\big(\tau_0^{S,-}>k-1 \big) \lesssim \frac{1}{\sqrt{k}}$
and that $(\Delta_k, \mathbf{E}_y)\stackrel{\mathrm{d}}{=}(\xi_1, \mathbf{E}_0)$, we can continue
the estimates in the display above to get
\begin{align}\label{eq:11}
	&
\sum_{k=1}^{[t^{1-\varepsilon}]} \mathbf{E}_y\big(S_k; \tau_0^{S,-}>k, S_k> t^{(1-\varepsilon)/2},  \tau_{t^{1/2-\varepsilon}}^{S,+} =k  \big)\nonumber\\
	&\lesssim \sum_{k=1}^{[t^{1-\varepsilon}]}  \mathbf{P}_0\big(  \xi_1  > t^{(1-\varepsilon)/2} - t^{1/2-\varepsilon}  \big)+ \sum_{k=1}^{[t^{1-\varepsilon}]} \frac{1}{\sqrt{k}} \mathbf{E}_0\big(\xi_1 ;  \xi_1  > t^{(1-\varepsilon)/2} - t^{1/2-\varepsilon}  \big)\nonumber\\
	&\leq  [t^{1-\varepsilon}] \frac{ \mathbf{E}_0 \left( \xi_1^2 ;\xi_1 > t^{(1-\varepsilon)/2} - t^{1/2-\varepsilon}\right)}{\left(t^{(1-\varepsilon)/2} - t^{1/2-\varepsilon}\right)^2}+ \int_0^{[t^{1-\varepsilon}]}\frac{1}{\sqrt{x}}\mathrm{d}x  \frac{ \mathbf{E}_0 \left( \xi_1^2 ;\xi_1 > t^{(1-\varepsilon)/2} - t^{1/2-\varepsilon}\right)}{t^{(1-\varepsilon)/2} - t^{1/2-\varepsilon}}\nonumber\\
	& \lesssim \mathbf{E}_0 \left( \xi_1^2 ;\xi_1 > t^{(1-\varepsilon)/2} - t^{1/2-\varepsilon}\right).
\end{align}
Combining \eqref{eq:12} and \eqref{eq:11}, we get that
\begin{align}
	I_3 \lesssim \mathbf{E}_0 \big( \xi_1^2;\xi_1 > t^{(1-\varepsilon)/2} - t^{1/2-\varepsilon}\big) \stackrel{t\to\infty}{\longrightarrow}0.
\end{align}

(iv) In this part, we deal with $I_4$. We allow $a$ to be $\infty$.
For $k\le [t^{1-\varepsilon}]$ and $x'>0$, define
\[
K(k, x'):
= \mathbf{P}_{x'}\left(\xi_{t-k}\leq a\sqrt{t}, \tau_0^->t-k, A_{t-k}\right).
\]
By the definition of $A_r$, we see that
\begin{align}\label{e:upper-K}
	K(k, x')
	\leq \mathbf{P}_{x'}\left(W_{t-k}\leq a\sqrt{t}+ (t-k)^{\frac{1}{2}-2\varepsilon}, \min_{s\leq t-k} W_s> -(t-k)^{\frac{1}{2}-2\varepsilon} \right).
\end{align}
Since $B_t:=W_t/\sigma$ is a
standard Brownian motion, we see that $(\sigma B_t, \mathbf{P}_y)\stackrel{\mathrm{d}}{=}(W_t, \mathbf{P}_{\sigma y})$. For any $z>0$, define
\begin{align}\label{change-of-measure}
	\frac{\mathrm{d} \mathbf{P}_z^\uparrow}{\mathrm{d} \mathbf{P}_z}\bigg|_{\sigma(B_s,s\leq t)}:= \frac{B_t}{z}1_{\{\min_{s\leq t} B_s>0\}}.
\end{align}
It is well-known that under $\mathbf{P}_z^\uparrow$, $B_s$ is a Bessel-3 process with transition density
\begin{align}\label{transition-density-Bessel}
	p_t^\uparrow(x,y)= \frac{y}{x\sqrt{2\pi t}}e^{-\frac{(y-x)^2}{2t}}\left(1-e^{-\frac{2xy}{t}}\right)1_{\{y>0\}}.
\end{align}
Set
\begin{align}
	x^*:=\frac{x'+(t-k)^{\frac{1}{2}-2\varepsilon}}{\sigma} \quad\mbox{and }\quad a^*:=  \frac{a\sqrt{t}}{\sigma \sqrt{t-k}}+ \frac{2}{\sigma(t-k)^{2\varepsilon}}.
\end{align}
Combining \eqref{e:upper-K}, \eqref{change-of-measure},  \eqref{transition-density-Bessel} and the inequality $1-e^{-x}\leq x$,  we obtain that
\begin{align}\label{e:upper-K-2}
	&K(k, x')
	\leq \mathbf{P}_{x^*}\big(B_{t-k}\leq a^* \sqrt{t-k}, \min_{s\leq t-k} B_s> 0 \big) = x^* \mathbf{E}_{x^*}^\uparrow \Big(\frac{1}{B_{t-k}}; B_{t-k}\leq a^* \sqrt{t-k}\Big)\nonumber\\
	& = \frac{1}{\sqrt{2\pi (t-k)}}\int_0^{a^* \sqrt{t-k}} e^{-\frac{(y-x^*)^2}{2(t-k)}}\big(1-e^{-\frac{2x^*y}{t-k}}\big)\mathrm{d}y\nonumber\\
	& \leq \frac{2x^*}{\sqrt{2\pi (t-k)^3}}\int_0^{a^* \sqrt{t-k}} ye^{-\frac{(y-x^*)^2}{2(t-k)}}\mathrm{d}y= \frac{2x^*}{\sqrt{2\pi (t-k)}}\int_0^{a^*} ye^{-\frac{y^2}{2} + \frac{yx^*}{\sqrt{t-k}} -\frac{(x^*)^2}{2(t-k)}}\mathrm{d}y\nonumber\\
	& \leq \frac{2}{\sqrt{2\pi (t-k)}}\Big( \frac{x'}{\sigma}+ \frac{t^{\frac{1}{2}-2\varepsilon}}{\sigma}\Big) \int_0^{a^*} ye^{-\frac{y^2}{2} + \frac{yx^*}{\sqrt{t-k}} }\mathrm{d}y\nonumber\\
	& = \frac{2}{\sqrt{2\pi (t-k)}}\Big( \frac{x'}{\sigma}+ \frac{t^{\frac{1}{2}-2\varepsilon}}{\sigma}\Big) \Big(\int_0^{a/\sigma } ye^{-\frac{y^2}{2} + \frac{yx^*}{\sqrt{t-k}} }\mathrm{d}y+\int_{a/\sigma}^{a^*} ye^{-\frac{y^2}{2} + \frac{yx^*}{\sqrt{t-k}} }\mathrm{d}y\Big),
\end{align}
where the last term on the right-hand side of the inequality above is $0$ when $a=\infty$. Note that for all $x'\leq t^{(1-\varepsilon)/2}$ and all $k\leq [t^{1-\varepsilon}]$,
\[
\frac{x^*}{\sqrt{t-k}} \lesssim \frac{t^{(1-\varepsilon)/2} + t^{\frac{1}{2}-2\varepsilon}}{\sqrt{t}}\lesssim t^{-\varepsilon/2}
\to 0 \quad\mbox{ as }t\to\infty.
\]
Note also that, for $a\in (0, \infty)$, $a^*-\sigma^{-1}a\to 0$ as $t\to\infty$.
Therefore, we conclude from \eqref{e:upper-K-2} that for
any $a\in (0,\infty]$ and $\delta_0>0$,
there exists $T>0$ such that when $t>T$, for all $x'\leq t^{(1-\varepsilon)/2}$,
\begin{align}
	K(k, x')
	\leq \frac{2(1+\delta_0)}{\sqrt{2\pi t}}
		     \Big(\frac{x'}{\sigma}+\frac{t^{\frac{1}{2}-2\varepsilon}}{\sigma} \Big)
	\Phi^+\left(\frac{a}{\sigma}\right).
\end{align}
Plugging this into the definition of $I_4$, we see that when $t$ is large enough,
\begin{align}
	I_4\leq  \frac{2(1+\delta_0)}{\sqrt{2\pi}} \Phi^+\left(\frac{a}{\sigma}\right)\sum_{k=1}^{[t^{1-\varepsilon}]} \mathbf{E}_y\Big( \Big(\frac{\xi_k}{\sigma}+
		\frac{t^{\frac{1}{2}-2\varepsilon}}{\sigma}\Big); \tau_0^->k, \xi_k\leq  t^{(1-\varepsilon)/2}, \tau_{t^{1/2-\varepsilon}}^{S,+} =k  \Big).
\end{align}
Note that
 on $\{ \tau_{t^{1/2-\varepsilon}}^{S,+} =k \}$, we have
$\xi_k=S_k \geq t^{1/2-\varepsilon}$,
which implies that
$t^{\frac{1}{2}-2\varepsilon}\leq \delta_0 \xi_k$ for $t$ large enough.
Hence, for $t$ large enough, by the inequality $R(x)\geq x$, we have
\begin{align}\label{Upper-bound}
		& I_4\leq  \frac{2(1+\delta_0)^2}{\sqrt{2\pi \sigma^2}}\Phi^+\left(\frac{a}{\sigma}\right) \sum_{k=1}^{[t^{1-\varepsilon}]} \mathbf{E}_y\big( \xi_k; \tau_0^->k, \xi_k\leq  t^{(1-\varepsilon)/2}, \tau_{t^{1/2-\varepsilon}}^{S,+} =k  \big)\nonumber\\
		& = \frac{2(1+\delta_0)^2}{\sqrt{2\pi \sigma^2}} \Phi^+\left(\frac{a}{\sigma}\right)\mathbf{E}_y\Big(R \big(\xi_{\tau_{t^{1/2-\varepsilon}}^{S,+}} \big);\tau_0^->\tau_{t^{1/2-\varepsilon}}^{S,+},\tau_{t^{1/2-\varepsilon}}^{S,+}\leq [t^{1-\varepsilon}]  \Big)\nonumber\\
		& \leq \frac{2(1+\delta_0)^2}{\sqrt{2\pi \sigma^2}} \Phi^+\left(\frac{a}{\sigma}\right)\mathbf{E}_y\Big(R \big(\xi_{\tau_{t^{1/2-\varepsilon}}^{S,+}\land [t^{1-\varepsilon}}]\big); \tau_0^->\tau_{t^{1/2-\varepsilon}}^{S,+}\land[t^{1-\varepsilon}]  \Big)\nonumber\\
		&=  \frac{2(1+\delta_0)^2}{\sqrt{2\pi \sigma^2}} R(y)\Phi^+\left(\frac{a}{\sigma}
		\right),
\end{align}
where in the last equality we used Lemma \ref{Prop-Renewal-Function}(ii).

For the lower bound, we have, similarly, for
$t^{1/2-\varepsilon}\leq x'\leq t^{(1-\varepsilon)/2}$,
\begin{align}
	K(k, x')
	\geq  \mathbf{P}_{x'}\Big(W_{t-k}\leq a\sqrt{t}- (t-k)^{\frac{1}{2}-2\varepsilon}, \min_{s\leq t-k} W_s> (t-k)^{\frac{1}{2}-2\varepsilon} \Big).
\end{align}
In this case, we define
\begin{align}
	x_*:=\frac{x'-(t-k)^{\frac{1}{2}-2\varepsilon}}{\sigma} \quad\mbox{and }\quad a_*:=  \frac{a\sqrt{t}}{\sigma \sqrt{t-k}}- \frac{2}{\sigma(t-k)^{2\varepsilon}}.
\end{align}
Then combining the inequalities $(y-x)^2 \leq y^2+x^2$,
$1-e^{-x}\geq x(1-x)$ for all $x,y>0$  and
an argument similar to that used in \eqref{e:upper-K-2}, we get
\begin{align}
	&K(k, x')
	\geq \mathbf{P}_{x_*}\Big(B_{t-k}\leq a_* \sqrt{t-k}, \min_{s\leq t-k} B_s> 0 \Big)  = \frac{1}{\sqrt{2\pi (t-k)}}\int_0^{a_* \sqrt{t-k}} e^{-\frac{(y-x_*)^2}{2(t-k)}}\big(1-e^{-\frac{2x_*y}{t-k}}\big)\mathrm{d}y\nonumber\\
	& \geq \frac{2x_*}{\sqrt{2\pi (t-k)^3}}\int_0^{a_* \sqrt{t-k}} y\Big(1-\frac{2x_*y}{t-k}\Big)e^{-\frac{y^2+ x_*^2}{2(t-k)}}\mathrm{d}y\nonumber\\
	&= \frac{2x_* e^{-\frac{x_*^2}{2(t-k)}}}{\sqrt{2\pi (t-k)}}\int_0^{a_* } y\Big(1-\frac{2x_*y}{\sqrt{t-k}}\Big)e^{-\frac{y^2}{2}}\mathrm{d}y\nonumber\\
	& \geq \frac{2e^{-\frac{x_*^2}{2(t-k)}} }{\sqrt{2\pi t}} \Big( \frac{x'}{\sigma}-\frac{t^{\frac{1}{2}-2\varepsilon}}{\sigma}\Big)  \Big( \int_0^{a_* } ye^{-\frac{y^2}{2}}\mathrm{d}y- \frac{2x_*}{\sqrt{t-k}}\int_0^{\infty} y^2e^{-\frac{y^2}{2}}\mathrm{d}y \Big).
\end{align}
Noting that for all $t^{1/2-\varepsilon}\leq x'\leq t^{(1-\varepsilon)/2}$ and $k\leq [t^{1-\varepsilon}]$,
\[
\frac{x_*}{\sqrt{t-k}}\lesssim \frac{t^{(1-\varepsilon)/2}+t^{\frac{1}{2}-2\varepsilon}}{\sqrt{t}}\lesssim t^{-\varepsilon/2} \quad \mbox{and } \quad \frac{t^{\frac{1}{2}-2\varepsilon}}{x'}\leq \frac{t^{\frac{1}{2}-2\varepsilon}}{ t^{1/2-\varepsilon}}=t^{-\varepsilon},
\]
and  for any  $a\in (0,\infty)$, $|\sigma^{-1}a-a_*|\lesssim t^{-\varepsilon}$.
Therefore,
 for any $\delta_0\in (0, 1)$,
when $t$ is large enough,
we have for all $k\leq [t^{1-\varepsilon}]$ and $t^{1/2-\varepsilon}\leq x'\leq t^{(1-\varepsilon)/2}$,
\begin{align}
	K(k, x')
	\geq
	     \frac{2(1-\delta_0)}{\sqrt{2\pi t \sigma^2}}
	x'\Phi^+\left(\frac{a}{\sigma}\right).
\end{align}
Therefore, for $t$ large enough,
\begin{align}
	& I_4\geq \frac{2(1-\delta_0)}{\sqrt{2\pi \sigma^2}}\Phi^+\left(\frac{a}{\sigma}\right) \sum_{k=1}^{[t^{1-\varepsilon}]} \mathbf{E}_y\big( \xi_k; \tau_0^->k,\xi_k\leq  t^{(1-\varepsilon)/2}, \tau_{t^{1/2-\varepsilon}}^{S,+} =k \big).
\end{align}
It follows from Lemma \ref{Prop-Renewal-Function} that $x\geq (1-\delta_0) R(x)$ for $x\geq t^{1/2-\varepsilon}$ with $t$ large enough.
Thus, when $t$ is large enough,
\begin{align}\label{Lower-Bound-1}
	& I_4\geq  \frac{2(1-\delta_0)^2}{\sqrt{2\pi \sigma^2}}\Phi^+\left(\frac{a}{\sigma}\right) \mathbf{E}_y\Big(R \big(\xi_{\tau_{t^{1/2-\varepsilon}}^{S,+}} \big); \tau_0^->\tau_{t^{1/2-\varepsilon}}^{S,+}, \tau_{t^{1/2-\varepsilon}}^{S,+}\leq [t^{1-\varepsilon}] \Big)\nonumber\\
	& =\frac{2(1-\delta_0)^2}{\sqrt{2\pi \sigma^2}}\Phi^+\left(\frac{a}{\sigma}\right)\Big(R(y)-\mathbf{E}_y\big(R \left(\xi_{[t^{1-\varepsilon}]} \right); \tau_0^->[t^{1-\varepsilon}], \tau_{t^{1/2-\varepsilon}}^{S,+}> [t^{1-\varepsilon}] \big) \Big),
\end{align}
where in the last inequality we used the following fact:
\begin{align*}
&R(y)\\=
&\mathbf{E}_y\Big(R \big(\xi_{\tau_{t^{1/2-\varepsilon}}^{S,+}\wedge[t^{1-\varepsilon}]} \big); \tau_0>\tau_{t^{1/2-\varepsilon}}^{S,+}\wedge[t^{1-\varepsilon}]\Big)\\
=& \mathbf{E}_y\Big(R \big(\xi_{\tau_{t^{1/2-\varepsilon}}^{S,+}} \big); \tau_0^->\tau_{t^{1/2-\varepsilon}}^{S,+}, \tau_{t^{1/2-\varepsilon}}^{S,+}\leq [t^{1-\varepsilon}] \Big)+\mathbf{E}_y\big(R \left(\xi_{[t^{1-\varepsilon}]} \right); \tau_0^->[t^{1-\varepsilon}], \tau_{t^{1/2-\varepsilon}}^{S,+}> [t^{1-\varepsilon}] \big).
\end{align*}
Noting that $\xi_{[t^{1-\varepsilon}]}\leq t^{1/2-\varepsilon}$ on $\{ \tau_{t^{1/2-\varepsilon}}^{S,+}> [t^{1-\varepsilon}]\}$
and applying Lemma \ref{lemma11} with $t=1$, we get
\begin{align}\label{Lower-Bound-2}
	&\mathbf{E}_y\Big(R \big(\xi_{[t^{1-\varepsilon}]} \big); \tau_0^->[t^{1-\varepsilon}], \tau_{t^{1/2-\varepsilon}}^{S,+}> [t^{1-\varepsilon}]  \Big) \nonumber\\
	& \leq R(t^{1/2-\varepsilon})\mathbf{P}_y\left(\tau_0^-> [t^{1-\varepsilon}]\right) \lesssim t^{1/2-\varepsilon}  \frac{y +1}{\sqrt{ [t^{1-\varepsilon}]}}\lesssim (y+1)t^{-\varepsilon/2}.
\end{align}
Combining \eqref{Upper-bound}, \eqref{Lower-Bound-1} and \eqref{Lower-Bound-2},
we conclude that
\begin{align}
	\frac{2(1+\delta_0)^2}{\sqrt{2\pi \sigma^2}} R(y)\Phi^+\left(\frac{a}{\sigma}\right)\geq \limsup_{t\to\infty} I_4\geq \liminf_{t\to\infty} I_4\geq \frac{2(1-\delta_0)^2}{\sqrt{2\pi \sigma^2}} R(y)\Phi^+\left(\frac{a}{\sigma}\right).
\end{align}
Letting $\delta_0\to 0$, we arrive at the assertion of the lemma.

\hfill$\Box$

\noindent
\textbf{Proof of Lemma \ref{Technical-lemma-1}: }
Define a sequence of measures
\[
\mu^{(t)}(D):= \sqrt{t} \mathbf{P}_y\Big(\frac{\xi_t}{\sigma\sqrt{t}}\in D, \tau_0^->t\Big)\quad \mbox{and}\quad \mu(D):= \frac{2 R(y)}{\sqrt{2\pi \sigma^2}} \int_D ze^{-\frac{z^2}{2}}\mathrm{d}z,
\quad D\in {\cal B}((0,\infty)).
\]
Lemma \ref{Solve-Technical-lemma-1} implies that for any $y>0$ and any $a\in (0,\infty)$,
\begin{align}
\lim_{t\to\infty} \mu^{(t)}((0, a])= \mu((0,a])
\end{align}
and that $\lim_{t\to\infty} \mu^{(t)}((0,\infty))= \mu((0,\infty))$.
Therefore, $\mu^{(t)}$  weakly converge to $\mu$ and this completes the proof of the lemma.

\hfill$\Box$

\subsection{Proof of Lemma \ref{lemma5} }\label{Appendix2}
\noindent
\textbf{Proof of Lemma \ref{lemma5}: }
First, it follows from  Lemma \ref{lemma2}(i)  that $\sup_{y>0} v_\infty^X(r,y)<\infty$ for any $ r>0$.
Next we prove that for any $y\in(0,\infty)$, $\lim_{r\to 0+}v_\infty^X(r,y)=\infty.$
By the definition of $v_\infty^{(t)}(r,y) $ in \eqref{Def-V-infty}, we have
\begin{align}\label{step_17}
	&v_\infty^{(t)}(r,y) \geq  t^{\frac{1}{\alpha-1}}\P_{\sqrt{t}y} \big(\zeta^{(0,\infty)}> tr, \inf_{s>0} \inf_{u\in N(s)} X_u(s)>0 \big)\nonumber\\
	& = t^{\frac{1}{\alpha-1}}\P_{\sqrt{t}y} \big(\zeta> tr, \inf_{s>0} \inf_{u\in N(s)} X_u(s)>0 \big)\nonumber\\
	& \geq t^{\frac{1}{\alpha-1}}\P_{\sqrt{t}y} \left(\zeta> tr \right)- t^{\frac{1}{\alpha-1}}\P_{\sqrt{t}y} \big(\inf_{s>0} \inf_{u\in N(s)} X_u(s)\leq 0 \big)\nonumber\\
	& = t^{\frac{1}{\alpha-1}}\P_{\sqrt{t}y} \left(\zeta> tr \right)- t^{\frac{1}{\alpha-1}}\P_{0} \big(\widetilde{M} \geq \sqrt{t}y \big),
\end{align}
where in the last equality, $\widetilde{M}$ is the maximal displacement of the critical branching L\'{e}vy process with spatial motion $-\xi$, branching rate $\beta$ and offspring distribution $\{p_k\}$. Combining \eqref{Survival-prob-zeta}, \eqref{Tail-probability-M} with $M$ replaced by $\widetilde{M}$ and \eqref{step_17}, we see that
\begin{align}
	v_\infty^X(r,y)\geq \frac{C(\alpha)}{r^{\frac{1}{\alpha-1}}}- \frac{\widetilde{\theta}(\alpha)}{y^{\frac{2}{\alpha-1}}} \stackrel{r\to0+}{\longrightarrow}
	\infty,
\end{align}
which implies $\lim_{r\to 0+}v_\infty^X(r,y)=\infty.$

Letting $t\to\infty$ first and then $y\to 0+$ in Lemma \ref{lemma2} (i), we easily see that $\lim_{y\to 0+}v_\infty^X(r,y)=0$ for any $r>0$.

Now we prove that $v^X_\infty$ satisfies the partial differential differential equation in \eqref{inner-eq}.
For  any  $0<w<r$ and $y>0$, by \eqref{Evolution-eq-V2},
\begin{align}\label{Evolution-eq-V2'}
	v_\infty^{(t)}(r,y)=\mathbf{E}_{y}\big(v_\infty^{(t)}(r-w,\xi_{w\land \tau_0^{(t),-}}^{(t)} )\big) - \mathbf{E}_{y}\big(\int_0^{w} \phi^{(t)}(v_\infty^{(t)}(r-s, \xi_{s\land \tau_0^{(t),-}}^{(t)} ))\mathrm{d} s\big).
\end{align}
By Lemma \ref{lemma2}(i),
$K:=\sup_{t>0, s\leq w, y\in \R} v_\infty^{(t)}(r-s, y)<\infty$.
By Lemma \ref{lemma2}(ii) and the definition of $\phi^{(t)}$ in \eqref{Def-V-phi-psi}, $\lim_{t\to\infty}\frac{\phi^{(t)}(v)}{v^{\alpha}}={\cal C}(\alpha)$ uniformly for $v\in[0,K]$. Note that $\varphi(\lambda)= \mathcal{C}(\alpha) \lambda^\alpha$ by definition \eqref{Stable-Branching-mechanism}. Therefore, for any $\varepsilon>0$, there exists $T_1>0$ such that when $t>T_1$,
\begin{align}\label{Fact1}
	\phi^{(t)}(v_\infty^{(t)}(r-s, \xi_{s\land \tau_0^{(t),-}}^{(t)} )) \geq (1-\varepsilon) 	\varphi
	\left( v_\infty^{(t)}(r-s, \xi_{s\land \tau_0^{(t),-}}^{(t)} )\right).
\end{align}
 Plugging \eqref{Fact2}, \eqref{Fact3} and \eqref{Fact1} into
\eqref{Evolution-eq-V2'}
we see that when $t> \max\{T_1, T_2\}$,
\begin{align}
	v_\infty^{(t)}(r,y) & \leq 3\varepsilon + \mathbf{E}_{y}\big(v_\infty^X(r-w,\xi_{w\land \tau_0^{(t),-}}^{(t)} )\big) \nonumber\\
	&\quad\quad -(1-\varepsilon) \mathbf{E}_{y}\Big(\int_0^{w\land \tau_0^{(t),-}} \varphi\big(\big(v_\infty^X(r-s, \xi_{s}^{(t)} ) -3\varepsilon\big)^+\big) \mathrm{d} s\Big).
\end{align}
Taking $t=t_k$ in the inequality above and letting $k\to\infty$, applying Lemma \ref{Coupling-Levy-BM}, we get that
\begin{align}
	v_\infty^X(r,y) & \leq 3\varepsilon + \mathbf{E}_{y}\left(v_\infty^X(r-w,
	W^0_w)\right) \nonumber\\
	&\quad\quad -(1-\varepsilon) \mathbf{E}_{y}\Big(\int_0^{w} \varphi\big(\big(v_\infty^X(r-s, 	W^0_s ) -3\varepsilon\big)^+\big) \mathrm{d} s\Big),
\end{align}
where $W^0_s$ is a Brownian motion with variance $\sigma^2t$ stopped upon exiting $(0, \infty)$. Now letting $\varepsilon \to0$, we finally conclude that
\begin{align}
	v_\infty^X(r,y) & \leq  \mathbf{E}_{y}\left(v_\infty^X(r-w, W^0_w	)\right)- \mathbf{E}_{y}\Big(\int_0^{w} \varphi\big(v_\infty^X(r-s,
	W^0_s) \big) \mathrm{d} s\Big).
\end{align}
A very similar argument for the lower bound implies that,  for any $0<w<r$ and $y>0$,   $v_\infty^X(r,y)$ solves the equation
\begin{align}\label{Evolution-eq-V-2}
	v_\infty^X (r,y)& =\mathbf{E}_{y}\left(v_\infty^X(r-w, 	W^0_w
	)\right)
	-\mathbf{E}_{y}\Big(\int_0^{w} \varphi(v_\infty^X(r-s, W^0_s	))\mathrm{d} s\Big).
\end{align}
This implies the desired result.
Indeed, we may rewrite \eqref{Evolution-eq-V-2} as
\begin{align}\label{Evolution-eq-V-2'}
	v_\infty^X (r+w,y)& =\mathbf{E}_{y}\left(v_\infty^X(r, 	W^0_w
	)\right)
	-\mathbf{E}_{y}\left(\int_0^{w} \varphi(v_\infty^X(r+w-s, W^0_s	))\mathrm{d} s\right),\quad
	y>0, r, w\geq 0.
\end{align}
For each  fixed $r>0$, set $f(y)=v_\infty^X(r, y)$. Then the function $u(w,y):=v_\infty^X(r+w, y)$ is the solution of the integral equation
\begin{align*}
	u(w,y)& =\mathbf{E}_{y}\left(f(	W^0_w)\right)
	-\mathbf{E}_{y}\big(\int_0^{w} \varphi(u(w-s, W^0_s))\mathrm{d} s\big),\quad
	y>0, w\geq 0.
\end{align*}
Recall that $W_s^0$ is the Brownian motion $W$ (with diffusion coefficient $\sigma^2$) stopped at $\tau_0^{W,-}$ and the generator of $W_s^0$ in the domain $(0,\infty)$ is $\frac{\sigma^2}{2}\frac{\mathrm{d}^2}{\mathrm{d} x^2}$.
Thus we have
$$\frac{\partial}{\partial w}v_\infty^X(r+w, y)= \frac{\sigma^2}{2}\frac{\partial^2}{\partial y^2} v_\infty^X(r+w, y) -
\varphi \left(v_\infty^X(r+w, y)\right), \quad r, w, y>0.$$
Since $r>0$ is arbitrary, we get
$$\frac{\partial}{\partial w}v_\infty^X(w,y)= \frac{\sigma^2}{2}\frac{\partial^2}{\partial y^2} v_\infty^X(w,y) -
\varphi \left(v_\infty^X(w,y)\right), \quad
\mbox {in } (0,\infty)\times (0,\infty).$$
The proof is now complete.
\hfill$\Box$

\subsection{Proofs of Propositions \ref{Technical-Prop-2} and \ref{Technical-Prop-1}}\label{Proof: Technical-Prop-1}

In this subsection,
we assume that {\bf (H1)} {\bf(H2)} and \eqref{eq:7} hold.
When $\xi$ is a standard Brownian motion, $\big(\xi_r^{(t)}, \mathbf{P}_y\big)\stackrel{\mathrm{d}}{=} \left(\xi_r, \mathbf{P}_y\right)$, and Proposition \ref{Technical-Prop-1} follows immediately from \cite[Proposition 4.5]{LZ}.

\begin{lemma}\label{Solve-Technical-Prop-1}
	Let $f$ be a bounded Lipschitz function on $\R_+$ with $f(0)=0$.
	
	(i) For any $r,y,w>0$ and $t>1$, it holds that
	\begin{align}
	& t^{\frac{1}{\alpha-1}}	\Big|\mathbf{E}_y \Big(1-\exp\Big\{-\frac{1}{t^{\frac{1}{\alpha-1}}} f\big(\xi_{r\land \tau_0^{(t),-}}^{(t)}\big)\Big\} \Big) -\mathbf{E}_{y+w} \Big(1-\exp\Big\{-\frac{1}{t^{\frac{1}{\alpha-1}}} f\big(\xi_{r\land \tau_0^{(t),-}}^{(t)}\big)\Big\} \Big) \Big|\nonumber\\
	& \lesssim \Big(\frac{1}{\log t} +  w \Big)(1+r^{-1/2})=: G_1^{(t)}(r, w).
	\end{align}
	
	(ii) For any $r,y,q>0$ and $t>1$, it holds that
	\begin{align}
		& t^{\frac{1}{\alpha-1}}	\Big|\mathbf{E}_y \Big(1-\exp\Big\{-\frac{1}{t^{\frac{1}{\alpha-1}}} f\big(\xi_{r\land \tau_0^{(t),-}}^{(t)}\big)\Big\} \Big) -\mathbf{E}_{y} \Big(1-\exp\Big\{-\frac{1}{t^{\frac{1}{\alpha-1}}} f\big(\xi_{(r+q)\land \tau_0^{(t),-}}^{(t)}\big)\Big\} \Big) \Big|\nonumber\\
		& \lesssim \Big(\frac{1}{\log t}+ q^{1/4}\Big)\big(1+r^{-1/2}\big)=: G_2^{(t)}(r,w).
	\end{align}
\end{lemma}
\textbf{Proof: }
(i) By the  inequality $x-(1-e^{-x})\lesssim x^2$ for all $x>0$, we have that
\begin{align}
	& t^{\frac{1}{\alpha-1}}	\Big|\mathbf{E}_y \Big(1-\exp\Big\{-\frac{1}{t^{\frac{1}{\alpha-1}}} f\big(\xi_{r\land \tau_0^{(t),-}}^{(t)}\big)\Big\} \Big) -\mathbf{E}_{y+w} \Big(1-\exp\Big\{-\frac{1}{t^{\frac{1}{\alpha-1}}} f\big(\xi_{r\land \tau_0^{(t),-}}^{(t)}\big)\Big\} \Big) \Big|\nonumber\\
	& \lesssim \frac{(\sup_{x\in \R}|f(x)| )^2}{t^{\frac{1}{\alpha-1}}}+ 	\Big|\mathbf{E}_y \Big(f\big(\xi_{r\land \tau_0^{(t),-}}^{(t)}\big) \Big) -\mathbf{E}_{y+w} \Big(f\big(\xi_{r\land \tau_0^{(t),-}}^{(t)}\big) \Big) \Big|\nonumber\\
	& \lesssim \frac{1}{t^{\frac{1}{\alpha-1}}}+ \Big|  \mathbf{E}_y
	 \Big(f\big(\xi_{r}^{(t)}\big); \tau_0^{(t),-}>r  \Big) -  \mathbf{E}_{y+w} \Big(f\big(\xi_{r}^{(t)}\big); \tau_0^{(t),-}>r \Big)\Big|.
\end{align}
Since $f$ is a bounded Lipschitz function, we have
\begin{align}\label{step_62}
	&\Big|  \mathbf{E}_y \Big(f\big(\xi_{r}^{(t)}\big); \tau_0^{(t),-}>r \Big) -  \mathbf{E}_{y+w} \Big(f\big(\xi_{r}^{(t)}\big); \tau_0^{(t),-}>r \Big)\Big|\nonumber\\
& = \Big|  \mathbf{E}_y \Big(f\big(\xi_{r}^{(t)}\big); \inf_{s\leq r } \xi_s^{(t)}> 0  \Big) -  \mathbf{E}_{y} \Big(f\big(\xi_{r}^{(t)} +w\big); \inf_{s\leq r } \xi_s^{(t)}> -w \Big)\Big| \nonumber\\	& \lesssim w+  \mathbf{P}_0\Big( \inf_{s\leq r } \xi_s^{(t)} \in (-w-y, -y]\Big).
\end{align}
 Recall the coupling in \eqref{Coupling-Levy-BM} and set $W_s^{(t)}= W_{st}/\sqrt{t}$,  we see that,  for any $\gamma\in (0,\frac{\delta}{2(2+\delta)})$,
 \begin{align}\label{step_64}
  & \mathbf{P}_0\Big( \inf_{s\leq r } \xi_s^{(t)} \in (-w-y, -y]\Big)\lesssim
   \frac{1}{t^{(\frac{1}{2}-\gamma)(\delta+2)-1}}
 + \mathbf{P}_0\Big( \min_{s\leq r} W_s^{(t)} \in (-w-y-t^{-\gamma}, -y+t^{-\gamma}]\Big)\nonumber\\
 	 &\lesssim \frac{1}{\log t} + \mathbf{P}_0\left(W_r\in [y-t^{-\gamma}, y+w+t^{-\gamma}) \right).
 \end{align}
Here the last inequality holds by the reflection principle. Therefore, by estimating the density of Brownian motion, we obtain that
\begin{align}
& \mathbf{P}_0\Big( \inf_{s\leq r } \xi_s^{(t)} \in (-w-y, -y]\Big)
	\lesssim  \frac{1}{\log t} + \frac{w+t^{-\gamma}}{\sqrt{r}}\lesssim \frac{1}{\log t} + \frac{w+(\log t)^{-1}}{\sqrt{r}},
\end{align}
which  gives the assertion (i).

(ii)
Similar to the beginning of the proof of (i), we also have
\begin{align}
		& t^{\frac{1}{\alpha-1}}	\Big|\mathbf{E}_y \Big(1-\exp\Big\{-\frac{1}{t^{\frac{1}{\alpha-1}}} f\big(\xi_{r\land \tau_0^{(t),-}}^{(t)}\big)\Big\} \Big) -\mathbf{E}_{y} \Big(1-\exp\Big\{-\frac{1}{t^{\frac{1}{\alpha-1}}} f\big(\xi_{(r+q)\land \tau_0^{(t),-}}^{(t)}\big)\Big\} \Big) \Big|\nonumber\\
	& \lesssim \frac{1}{\log t}+ \Big|  \mathbf{E}_y \Big(f\big(\xi_{r}^{(t)}\big); \tau_0^{(t),-}>r \Big) -  \mathbf{E}_{y} \Big(f\big(\xi_{r+q}^{(t)}\big); \tau_0^{(t),-}>r+q  \Big)\Big|.
\end{align}
Again using the fact that $f$ is bounded Lipschitz, we have
\begin{align}
 &\Big|  \mathbf{E}_y \Big(f\big(\xi_{r}^{(t)}\big); \tau_0^{(t),-}>r  \Big) -  \mathbf{E}_{y} \Big(f\big(\xi_{r+q}^{(t)}\big); \tau_0^{(t),-}>r+q  \Big)\Big|\nonumber\\
	 &\lesssim \sup_{x>0} |f(x)| \mathbf{P}_y\big( \tau_0^{(t),-} \in (r, r+q] \big) + \mathbf{E}_y\big(\big| \xi_{r}^{(t)}-\xi_{r+q}^{(t)} \big|\big)\nonumber\\
	 &\lesssim
 \mathbf{P}_0\Big( \inf_{s\leq r}\xi_s^{(t)} \geq -y,   \inf_{s\leq r+q}\xi_s^{(t)} < -y \Big) +\sqrt{q}.
\end{align}
Here in the last inequality we used the fact that $\mathbf{E}_y\big(\big| \xi_{r}^{(t)}-\xi_{r+q}^{(t)} \big|\big) \leq \sqrt{\mathbf{E}_0 (|\xi_1|^2) q }$.  Therefore, using a coupling argument similar to that leading to \eqref{step_64}, we get
\begin{align}
&\mathbf{P}_0\Big( \inf_{s\leq r}\xi_s^{(t)} \geq -y,   \inf_{s\leq r+q}\xi_s^{(t)} < -y \Big)\nonumber\\
	& \lesssim \frac{1}{\log t}+ \mathbf{P}_0\Big( \min_{s\leq r} W_s >-y -t^{-\gamma},   \min_{s\leq r+q}W_s \leq -y +t^{-\gamma} \Big).
\end{align}
Using
\begin{align}
	&\mathbf{P}_0\Big( \min_{s\leq r} W_s >-y -t^{-\gamma},   \min_{s\leq r+q}W_s \leq -y +t^{-\gamma} \Big)\nonumber\\
	& \leq \mathbf{P}_0\big( W_r\in (-y-t^{-\gamma}, -y+t^{-\gamma } +q^{1/4}) \big)+ \mathbf{P}_0\Big( W_r> -y+t^{-\gamma } +q^{1/4},   \min_{s\leq r+q}W_s \leq -y +t^{-\gamma}\Big)\nonumber\\
	&\lesssim \frac{t^{-\gamma}+q^{1/4}}{\sqrt{r}}+ \mathbf{P}_0 \big(\min_{s\leq q} W_s < -q^{1/4}\big) \lesssim \frac{t^{-\gamma}+q^{1/4}}{\sqrt{r}}+ q^{1/4}\lesssim  \Big(\frac{1}{\log t}+ q^{1/4}\Big)\big(1+r^{-1/2}\big),
\end{align}
we easily get the assertion of (ii).

\hfill$\Box$

The following lemma is a generalized Gronwall inequality. We omit the proof here since the proof is standard.

\begin{lemma}\label{Generlized-Grownwall}
	Suppose that $F$ and $G$ are two bounded non-negative measurable function on $[0,T]$.  If for any $r\in [0,T]$,
	\[
	F(r)\leq G(r)+ C\int_0^r F(s)\mathrm{d}s,
	\]
	then we have for all $r\in [0,T]$,
	\[
	F(r)\leq G(r)+ C\int_0^r e^{C(r-s)}G(s)\mathrm{d}s
	\]
\end{lemma}

\noindent
\textbf{Proof of Proposition \ref{Technical-Prop-2}: }
(i) By Lemma \ref{Solve-Technical-Prop-1} (i), we have that
\begin{align}\label{step_63}
	  & \big|v_f^{(t)}(r,y)- v_f^{(t)}(r, y+w)\big|  \lesssim G_1^{(t)}(r,w) + \Big| \mathbf{E}_{y}\Big(\int_0^{r} \phi^{(t)}(v_f^{(t)}(r-s, \xi_{s}^{(t)} ))
	  1_{\{\tau_0^{(t),-}>s\}}\mathrm{d} s\Big) \nonumber\\
	  &\quad\quad- \mathbf{E}_{y+w}\Big(\int_0^{r} \phi^{(t)}(v_f^{(t)}(r-s, \xi_{s}^{(t)} ))1_{\{\tau_0^{(t),-}>s\}}\mathrm{d} s\Big) \Big|.
\end{align}
Using the fact that $\left| \phi^{(t)}(u) -\phi^{(t)}(v)\right| \lesssim |u-v|$ for all $u,v\in [0,K]$ and $t>K$, and an argument similar to that leading to  \eqref{step_62}, we get that, for $t$ large enough so that $t> \sup_{x\in\R} |f(x)|$,
the second term on the right-hand side of \eqref{step_63} is bounded above by a constant multiple of
\begin{align}
	&	\int_0^r \sup_{y\in\R}\big| v_f^{(t)}(r-s, y)-   v_f^{(t)}(r-s, y+w)\big|\mathrm{d}s + \int_0^r
 \mathbf{P}_0\big( \inf_{\ell\leq s } \xi_\ell^{(t)} \in (-w-y, -y]\big)\mathrm{d}s\nonumber\\
	&\lesssim \int_0^r \sup_{y\in\R}\big| v_f^{(t)}(r-s, y)-   v_f^{(t)}(r-s, y+w)\big|\mathrm{d}s  +  \int_0^r  \frac{1}{\log t} + \frac{w+(\log t)^{-1}}{\sqrt{s}}\mathrm{d}s \nonumber\\
	&\lesssim  \int_0^r \sup_{y\in\R}\big| v_f^{(t)}(r-s, y)-   v_f^{(t)}(r-s, y+w)\big|\mathrm{d}s + \frac{1}{\log t} + w\nonumber\\
	&\leq  \int_0^r \sup_{y\in\R}\big| v_f^{(t)}(r-s, y)-   v_f^{(t)}(r-s, y+w)\big|\mathrm{d}s +G_1^{(t)}(r,w).
\end{align}
Plugging this into \eqref{step_63},
we conclude that there exists a constant $L$ independent of $r$ and $t$ such that for all $r\in [0,T]$ and $t>1$,
\begin{align}
	& \sup_{y>0} \big|v_f^{(t)}(r,y)- v_f^{(t)}(r, y+w)\big|  \leq L G_1^{(t)}(r,w) \nonumber\\
	&\quad + L\int_0^r \sup_{y\in\R}\big| v_f^{(t)}(r-s, y)-   v_f^{(t)}(r-s, y+w)\big|\mathrm{d}s.
\end{align}
Applying Lemma \ref{Generlized-Grownwall}, we obtain that for all $r\in [0, T]$, we have that
\begin{align}
	&\sup_{y>0} \big|v_f^{(t)}(r,y)- v_f^{(t)}(r, y+w)\big|  \leq L G_1^{(t)}(r,w)+ L\int_0^r e^{C(r-s)} G_1^{(t)}(s,w)\mathrm{d}s\nonumber\\
	&\lesssim \Big(\frac{1}{\log t} +  w \Big)(1+r^{-1/2}) + \int_0^r \Big(\frac{1}{\log t} +  w \Big)(1+s^{-1/2})\mathrm{d}s\nonumber\\
	&\lesssim \Big(\frac{1}{\log t} +  w \Big)(1+r^{-1/2}).
\end{align}
This completes the proof of (i).

(ii) By Lemma \ref{Solve-Technical-Prop-1} (ii), we see that
\begin{align}
	& \big|v_f^{(t)}(r,y)- v_f^{(t)}(r+q, y)\big| \lesssim G_2^{(t)}(r,w) +\Big| \mathbf{E}_{y}\Big(\int_0^{r} \phi^{(t)}(v_f^{(t)}(r-s, \xi_{s}^{(t)} ))1_{\{\tau_0^{(t),-}>s\}}\mathrm{d} s\Big) \nonumber\\
	&\quad\quad- \mathbf{E}_{y}\Big(\int_0^{r+q} \phi^{(t)}(v_f^{(t)}(r+q-s, \xi_{s}^{(t)} ))1_{\{\tau_0^{(t),-}>s\}}\mathrm{d} s\Big) \Big|\nonumber\\
	&\lesssim G_2^{(t)}(r,w)+q+  \mathbf{E}_{y}\Big(\int_0^{r} \Big| \phi^{(t)}(v_f^{(t)}(r-s, \xi_{s}^{(t)} )) -  \phi^{(t)}(v_f^{(t)}(r+q-s, \xi_{s}^{(t)} ))  \Big| 1_{\{\tau_0^{(t),-}>s\}}\mathrm{d} s\Big).
\end{align}
Again by the inequality $|\phi^{(t)}(u)-\phi^{(t)}(v)|\lesssim |u-v|$, we get that the last term on the right-hand side of the inequality above is bounded from above
by a constant multiple of
\begin{align}
	\int_0^{r} \sup_{y>0} \big| v_f^{(t)}(r-s, y) - v_f^{(t)}(r+q-s, y )  \big| \mathrm{d} s.
\end{align}
Therefore, there exists a constant $L$ independent of $t,q$ and $r$ such that for all $r+q\leq T$,
\begin{align}
	& \sup_{y>0} \big|v_f^{(t)}(r,y)- v_f^{(t)}(r+q, y)\big| \nonumber\\
	& \leq L G_2^{(t)}(r,w)+Lq+ L \int_0^{r} \sup_{y>0} \big| v_f^{(t)}(r-s, y) - v_f^{(t)}(r+q-s, y )  \big| \mathrm{d} s.
\end{align}
Applying Lemma \ref{Generlized-Grownwall} for any fixed $q$ yields that
\begin{align}
	& \sup_{y>0} \big|v_f^{(t)}(r,y)- v_f^{(t)}(r+q, y)\big| \leq L G_2^{(t)}(r,w)+Lq + L\int_0^r e^{L(r-s)} \big(L G_2^{(t)}(s,w)+Lq \big)\mathrm{d}s\nonumber\\
	&\lesssim G_2^{(t)}(r,w)+q + \int_0^r \Big(\frac{1}{\log t}+ q^{1/4}\Big)\big(1+s^{-1/2}\big)  \mathrm{d}s\nonumber\\
	&\lesssim \Big(\frac{1}{\log t}+ q^{1/4}\Big)\big(1+r^{-1/2}\big),
\end{align}
which completes the proof of (ii).

\hfill$\Box$

\noindent
\textbf{Proof of  Proposition \ref{Technical-Prop-1}: } Fix a continuous function $f\in B_b^+((0,\infty))$ and $T>0$.  By Lemma \ref{Technical-lemma-2}, without loss of generality, we assume that $f$ is Lipschitz continuous. Since $v_f^{(t)}(r,y)$ is uniformly bounded for all $r\in [0,T], y>0$ and $t>1$, we can find a sequence $\{t_k\}$ and a limit
$v_f^X(r,y)$ such that
\begin{align}
	v_f^X(r,y)
	= \lim_{k\to\infty} v_f^{(t_k)}(r,y),\quad \mbox{for all }\ r\in[0,T]\cap \Q, y\in (0,\infty)\cap \Q.
\end{align}
Proposition \ref{Technical-Prop-2} implies that for any $r\in (0,T)$, $y>0$ and any  $((0, T)\cap \Q)\times ((0, \infty)\times \Q)\ni (r_m, y_m)\to (r,y)$,
we have that
$v_f^X(r_m, y_m)$ is a Cauchy sequence. Thus we define, for any $r\in (0,T)$ and $y>0$,
\begin{align}
	v_f^X(r,y):	=\lim_{((0, T)\cap \Q)\times ((0, \infty)\times \Q)\ni (r_m, y_m)\to (r,y)}
	v_f^X(r_m, y_m).
\end{align}
Using an argument similar to that leading to Lemma \ref{lemma4}, we can get
\begin{align}\label{step_65}
v_f^X(r,y)
	= \lim_{k\to\infty} v_f^{(t_k)}(r,y),\quad \mbox{for all } r\in (0,T), y\in (0,\infty).
\end{align}
Note that $v_f^{(t)}$ solves the equation
\begin{align}
		v_f^{(t)}(r,y)=  t^{\frac{1}{\alpha-1}}	\mathbf{E}_y \Big(1-\exp\Big\{-\frac{1}{t^{\frac{1}{\alpha-1}}} f\big(\xi_{r\land \tau_0^{(t),-}}^{(t)}\big)\Big\} \Big)
		- \mathbf{E}_{y}\Big(\int_0^{r} \phi^{(t)}(v_f^{(t)}(r-s,  \xi_{s\land \tau_0^{(t),-}}^{(t)}))\mathrm{d}	s\Big).
\end{align}
Using the invariance principle and an argument similar to that
leading to \eqref{Evolution-eq-V-2}, we arrive at the desired result.

\hfill$\Box$

\subsection{Proof of Lemma \ref{lemma9} and Proposition \ref{prop2}}\label{Appendix}

\noindent
\textbf{Proof of Lemma \ref{lemma9}:}
We first show that $\lim_{y\to0+}K^X(y)=0$.
Taking $z=\frac{1}{2}$ in Lemma \ref{lemma6} and applying Lemma \ref{lemma7} (i), we get that for $y<\frac{1}{2}$,
\begin{align}\label{step_70}
	&K^{(x)}(y)\leq
	\mathbf{E}_y\Big(K^{(x)}\Big(\xi_{\tau_{1/2}^{(x^2),+}}^{(x^2)} \Big)	;\tau_{1/2}^{(x^2),+}< \tau_0^{(x^2),-} \Big) \nonumber\\
	& \leq K^{(x)}\big(\frac{2}{3}\big)
	\mathbf{P}_y\Big(  \tau_{1/2}^{(x^2),+}< \tau_0^{(x^2),-}\Big)+ x^{\frac{2}{\alpha -1}} \mathbf{P}_y\Big(\xi_{\tau_{1/2}^{(x^2),+}}^{(x^2)} > \frac{2}{3}\Big)\nonumber\\
	& = K^{(x)}\big(\frac{2}{3}\big)
	\mathbf{P}_y\Big(  \tau_{1/2}^{(x^2),+}< \tau_0^{(x^2),-}\Big)
	+ x^{\frac{2}{\alpha -1}}\mathbf{P}_{ -\frac{1}{2}x +xy}\Big( \xi_{\tau_0^+} > \frac{1}{6}x\Big)\nonumber\\
	& \lesssim
	\mathbf{P}_y\Big(  \tau_{1/2}^{(x^2),+}< \tau_0^{(x^2),-}\Big)+ x^{\frac{2}{\alpha -1}}
	\frac{6^{r_0-2}}{x^{r_0-2}}\mathbf{E}_{ -\frac{1}{2}x +xy}\Big( \left| \xi_{\tau_0^+} \right| ^{r_0-2}	\Big),
\end{align}
where in the last inequality we used Markov's inequality.
It follows from  Lemma \ref{lemma12} that
$$
\mathbf{E}_{ -\frac{1}{2}x +xy}\Big( \left| \xi_{\tau_0^+} \right| ^{r_0-2}	\Big)\le C
$$
for some constant $C>0$.  Thus, since $r_0-2> \frac{2}{\alpha-1}$, taking $x=x_k$ and letting $k\to\infty$ in \eqref{step_70}, we get that
\begin{align}
	K^X(y)\lesssim
	\mathbf{P}_y\big(\tau_{1/2}^{W,+}< \tau_0^{W,-}\big)
	\stackrel{y\to 0+}{\longrightarrow}0.
\end{align}

Next we show that $\lim_{y\to1-}K^X(y)=\infty$. Note that
\begin{align}\label{step_43}
	& K^{(x)}(y) \geq x^{\frac{2}{\alpha-1}}\P_{xy}\Big(M^{(0,\infty)}\geq x, \ \inf_{t>0}\inf_{u\in N(t)}X_u(t)>0\Big)\nonumber\\
	& =  x^{\frac{2}{\alpha-1}}\P_{xy}\Big(M\geq x, \ \inf_{t>0}\inf_{u\in N(t)}X_u(t)>0\Big) \nonumber\\
	& \geq x^{\frac{2}{\alpha-1}}\P_{xy}\left(M\geq x\right) -x^{\frac{2}{\alpha-1}}\P_{xy}\Big( \inf_{t>0}\inf_{u\in N(t)}X_u(t)\leq 0\Big) \nonumber\\
	& = x^{\frac{2}{\alpha-1}}\P\left(M\geq x(1-y)\right)- x^{\frac{2}{\alpha-1}}\P(\widetilde{M}\geq xy),
\end{align}
where $\widetilde{M}$ is the maximal displacement of the critical branching L\'{e}vy process with branching rate $\beta$, offspring distribution $\{p_k\}$ and spatial motion $-\xi$. Applying \eqref{Tail-probability-M}  to $\widetilde{M}$ and $M$  we see that under {\bf(H4)},
\[
K^X(y)= \lim_{k\to\infty} K^{(x_k)}(y)
\geq \frac{\theta(\alpha)}{(1-y)^{\frac{2}{\alpha-1}}} -\frac{\widetilde{\theta(}\alpha)}{y^{\frac{2}{\alpha-1}}} \stackrel{y\to 1-}{\longrightarrow} +\infty.
\]

Finally, we show that $K^X(\cdot)$ satisfies the differential equation in    \eqref{PDEin(0,1)}. We fix an arbitrary  $z\in (0,1)$ in the remainder of this proof. By Lemma \ref{lemma7},
\[
\sup_{s\in (0,\tau_z^{(x^2_k),+})} K^{(x_k)}\big(\xi_s^{(x^2_k)}\big) \leq K^{(x_k)}(z) \lesssim \frac{1}{(1-z)^{\frac{2}{\alpha-1}}}.
\]
Therefore, by Lemma \ref{lemma2}(ii), for any $\varepsilon>0$, there exists $N>0$ such that
for any $k>N$ and $s\in (0,\tau_z^{(x^2_k),+})$,
\[
\mathcal{C}(\alpha) (1-\varepsilon)\leq
\frac{\psi^{(x^2_k)}\big( K^{(x_k)}\big(\xi_s^{(x^2_k)}\big)\big) }{\big(  K^{(x_k)}\big(\xi_s^{(x^2_k)}\big)\big)^{\alpha -1}}
\leq  \mathcal{C}(\alpha) (1+\varepsilon).
\]
Recall that $\varphi(\lambda)= \mathcal{C}(\alpha) \lambda^\alpha$ defined in  \eqref{Stable-Branching-mechanism}.
Set $\psi^X(v):= \varphi(v)/v$. For simplicity, we will use $x_k$ as $x$ in the remainder of this proof.
Applying the display above to Lemma \ref{lemma6}, we see that for $k>N$,
\begin{align}\label{step_44}
	&\mathbf{E}_y\Big(\exp\Big\{- (1-\varepsilon) \int_0^{\tau_z^{(x^2),+}} \psi^X\big(  K^{(x)}(\xi_s^{(x^2)})\big) \mathrm{d}s \Big\}  K^{(x)}\big(\xi_{\tau_z^{(x^2),+}}^{(x^2)} \big); \tau_z^{(x^2),+}< \tau_0^{(x^2),-} \Big)\geq K^{(x)}(y)\nonumber \\
	&\geq \mathbf{E}_y\Big(\exp\Big\{- (1+\varepsilon) \int_0^{\tau_z^{(x^2),+}} \psi^X\big(  K^{(x)}(\xi_s^{(x^2)})\big) \mathrm{d}s \Big\}  K^{(x)}\big(\xi_{\tau_z^{(x^2),+}}^{(x^2)} \big); \tau_z^{(x^2),+}< \tau_0^{(x^2),-} \Big).
\end{align}
Now we will let $k\to+\infty$ in \eqref{step_44}.  For the upper bound,
we note that for any $\delta\in (0,1-z)$,
$K^{(x)}\big(\xi_{\tau_z^{(x^2),+}}^{(x^2)} \big)\leq K^{(x)}(z+\delta)$ on the event $\big\{\xi_{\tau_z^{(x^2),+}}^{(x^2)} \leq z+\delta\big\}$ and $K^{(x)}\big(\xi_{\tau_z^{(x^2),+}}^{(x^2)} \big)\leq x^{\frac{2}{\alpha-1}}$ on the event $\big\{\xi_{\tau_z^{(x^2),+}}^{(x^2)} > z+\delta\big\}$. Thus,
\begin{align}\label{step_44-2}	
	&\mathbf{E}_y\Big(\exp\Big\{- (1-\varepsilon) \int_0^{\tau_z^{(x^2),+}} \psi^X\big(  K^{(x)}(\xi_s^{(x^2)})\big) \mathrm{d}s \Big\}  K^{(x)}\big(\xi_{\tau_z^{(x^2),+}}^{(x^2)} \big); \tau_z^{(x^2),+}< \tau_0^{(x^2),-} \Big)\nonumber\\
	&\leq K^{(x)}(z+\delta) \mathbf{E}_y\Big(\exp\Big\{- (1-\varepsilon) \int_0^{\tau_z^{(x^2),+}} \psi^X\big(  K^{(x)}(\xi_s^{(x^2)})\big) \mathrm{d}s \Big\}; \tau_z^{(x^2),+}< \tau_0^{(x^2),-}\Big)\nonumber\\
	&\quad + x^{\frac{2}{\alpha-1}} \mathbf{P}_y\big( \xi_{\tau_z^{(x^2),+}}^{(x^2)}>z+\delta\big).
\end{align}
The last term of the upper bound converges to $0$ as
$k\to\infty$. Indeed, since $r-2> 2/(\alpha-1)$, by Lemma \ref{lemma12}, we have
\begin{align}
	x^{\frac{2}{\alpha-1}} \mathbf{P}_y\Big( \xi_{\tau_z^{(x^2),+}}^{(x^2)}>z+\delta\Big)=  x^{\frac{2}{\alpha-1}} \mathbf{P}_{-x(z-y)}\big( \xi_{\tau_0^+}>x\delta\big)\leq \frac{x^{\frac{2}{\alpha-1}}}{(\delta x)^{r-2}}\sup_{w>0}
	    \mathbf{E}_{-w}\big(\xi_{\tau_0^+}^{r-2}\big)
	\stackrel{x\to\infty}{\longrightarrow}0.
\end{align}
Therefore, combining \eqref{step_44} and \eqref{step_44-2}, letting $k\to\infty$, we get
\begin{align}\label{e:step_44-rs}
	&K^X(y) \nonumber\\
	&\leq K^X(z+\delta) \limsup_{k\to\infty} \mathbf{E}_y\Big(\exp\Big\{- (1-\varepsilon) \int_0^{\tau_z^{(x^2),+}} \psi^X\big(  K^{(x)}(\xi_s^{(x^2)})\big) \mathrm{d}s \Big\}; \tau_z^{(x^2),+}< \tau_0^{(x^2),-} \Big).
\end{align}
Using the continuity of $K^X(\cdot )$ and the fact that $\lim_{y\to0+}K^X(y)=0$, we get that, for any $\varepsilon>0$, there exist $L\in \N$ and $0=w_0<w_1<...< w_L=z$ such that
\[
\max_{j\in \{1,...,L\}} \Big|K^X(w_j) - K^X(w_{j-1}) \Big|<\varepsilon.
\]
Let $T=T(L,\varepsilon)$ be large enough so that for all $k\geq T$,
\[
\max_{j\in \{0,..., L\}} \Big|K^{(x_k)}(w_j) - K^X(w_j)\Big|<\varepsilon.
\]
For $w\in [0, z)$, we have $w\in [w_{j-1}, w_j)$ for some $j\in \{1,...,L\}$. Using the fact that both $K^{(x_k)}(w)$ and $K^X(w)$ are increasing in $w$, we get
\begin{align}
	K^{(x_k)}(w) \geq K^{(x_k)}(w_{j-1})\geq K^X(w_{j-1})-\varepsilon \geq K^X(w_{j})- 2\varepsilon\geq 	K^X(w)-2\varepsilon.
\end{align}
Therefore,  when $k$ is sufficiently large,
\begin{align}\label{step_44-3}
	&\mathbf{E}_y\Big(\exp\Big\{- (1-\varepsilon) \int_0^{\tau_z^{(x^2),+}} \psi^X\big(  K^{(x)}(\xi_s^{(x^2)})\big) \mathrm{d}s \Big\}; \tau_z^{(x^2),+}< \tau_0^{(x^2),-} \Big)\nonumber\\
	&\leq \mathbf{E}_y\Big(\exp\Big\{- (1-\varepsilon) \int_0^{\tau_z^{(x^2),+}} \psi^X\big(  \big(K^X(\xi_s^{(x^2)}) -2\varepsilon \big)^+\big) \mathrm{d}s \Big\}; \tau_z^{(x^2),+}< \tau_0^{(x^2),-} \Big).
\end{align}
Plugging \eqref{step_44-3} into \eqref{e:step_44-rs}, we obtain
\begin{align}\label{step_44-4}
	&\frac{K^X(y)}{K^X(z+\delta)} \nonumber\\
	&\leq  \limsup_{k\to\infty} \mathbf{E}_y\Big(\exp\Big\{- (1-\varepsilon) \int_0^{\tau_z^{(x^2),+}} \psi^X\Big(  \big(K^X(\xi_s^{(x^2)}) -2\varepsilon \big)^+\Big) \mathrm{d}s \Big\}; \tau_z^{(x^2),+}< \tau_0^{(x^2),-} \Big).
\end{align}
Fix a large real number $A$ and an integer $N$,
and set $t_i= \frac{A}{N}i$ for $i\in \{0,...,N\}$.
Then we have
\begin{align}
	&\int_0^{\tau_z^{(x^2),+} \land A} \psi^X\Big(  \big(K^X(\xi_s^{(x^2)}) -2\varepsilon \big)^+\Big) \mathrm{d}s = \sum_{i=1}^N \int_{t_{i-1}}^{t_i} 1_{\{s<\tau_z^{(x^2),+} \}}  \psi^X\Big(  \big(K^X(\xi_s^{(x^2)}) -2\varepsilon \big)^+\Big) \mathrm{d}s \nonumber\\
	& \geq \sum_{i=1}^N \frac{A}{N}1_{\{t_{i}<\tau_z^{(x^2),+} \}}  \psi^X\Big(  \big(K^X\left(\inf_{s\in[t_{i-1}, t_i]}\xi_s^{(x^2)}\right) -2\varepsilon \big)^+\Big) .
\end{align}
Using an argument similar to that in \cite[Step 1 in Lemma 3.3]{HJRS}, with \cite[Lemma 2.4]{HJRS} there replaced by Lemma  \ref{Coupling-Levy-BM}, we see that
\begin{align}\label{step_44-5}
	&\limsup_{k\to\infty} \mathbf{E}_y\Big(\exp\Big\{- (1-\varepsilon) \int_0^{\tau_z^{(x^2),+}} \psi^X\big(  \big(K^X(\xi_s^{(x^2)}) -2\varepsilon \big)^+\big) \mathrm{d}s \Big\}; \tau_z^{(x^2),+}< \tau_0^{(x^2),-} \Big)\nonumber\\
	& \leq  \limsup_{k\to\infty} \mathbf{P}_y\big(\tau_0^{(x^2),-}> A\big) + \mathbf{E}_y\Big(\exp\Big\{- (1-\varepsilon) \sum_{i=1}^N \frac{A}{N}1_{\{t_{i}<\tau_z^{W,+} \}} \nonumber\\
	&\qquad\qquad\qquad\times  \psi^X\big(  \big(K^X\big(\inf_{s\in[t_{i-1}, t_i]}W_s \big) -2\varepsilon \big)^+\big) \Big\}; \tau_z^{W,+}< \tau_0^{W,-}<A \Big),
\end{align}
where $\tau_z^{W,+}$ is the exit time of the  process $W$ on $(-\infty, z)$.
Combining \eqref{step_44-4} and \eqref{step_44-5}, taking $N\to\infty$ first and then $A\to +\infty$, we get
\begin{align}
	&\frac{K^X(y)}{K^X(z+\delta)}  \nonumber\\
	& \leq  \limsup_{A\to\infty} \mathbf{P}_y\big(\tau_0^{W,-}> A\big) + \limsup_{A\to\infty}  \limsup_{N\to\infty} \mathbf{E}_y\Big(\exp\Big\{- (1-\varepsilon) \sum_{i=1}^N \frac{A}{N}1_{\{t_{i}<\tau_z^{W,+} \}} \nonumber\\
	&\qquad\qquad\qquad\times  \psi^X\big(  \big(K^X\big(\inf_{s\in[t_{i-1}, t_i]}W_s \big) -2\varepsilon \big)^+\big) \Big\}; \tau_z^{W,+}< \tau_0^{W,-}<A \Big)\nonumber\\
	& = \limsup_{A\to\infty}  \mathbf{E}_y\Big(\exp\Big\{- (1-\varepsilon) \int_0^{A\land \tau_z^{W,+}}   \psi^X\big(  \big(K^X\big(W_s \big) -2\varepsilon \big)^+\big)\mathrm{d}s \Big\}; \tau_z^{W,+}< \tau_0^{W,-}<A \Big)\nonumber\\
	& =  \mathbf{E}_y\Big(\exp\Big\{- (1-\varepsilon) \int_0^{\tau_z^{W,+}}   \psi^X\big(  \big(K^X\left(W_s \right) -2\varepsilon \big)^+\big)\mathrm{d}s \Big\}; \tau_z^{W,+}< \tau_0^{W,-}\Big).
\end{align}
Since $\varepsilon$ and $\delta$ are independent, letting $\varepsilon, \delta\to 0$ in the above inequality, we conclude that
\begin{align}
	K^X(y) \leq K^X(z)\mathbf{E}_y\Big(\exp\Big\{-  \int_0^{\tau_z^{W,+}} \psi^X\big(  K^X(W_s)\big) \mathrm{d}s \Big\}; \tau_z^{W,+}< \tau_0^{W, -} \Big).
\end{align}
Using a similar argument, we can prove that
\begin{align}
	K^X(y) \geq K^X(z)\mathbf{E}_y\Big(\exp\Big\{-  \int_0^{\tau_z^{W,+}} \psi^X\big(  K^X(W_s)\big) \mathrm{d}s \Big\}; \tau_z^{W,+}< \tau_0^{W, -} \Big).
\end{align}
Therefore,
\begin{align}\label{Feynman-Kac-of-K}
	&K^X(z)\mathbf{E}_y\Big(\exp\Big\{-  \int_0^{\tau_z^{W,+}} \psi^X\big(  K^X(W_s)\big) \mathrm{d}s \Big\}; \tau_z^{W,+}< \tau_0^{W, -} \Big)=
	K^X(y).
\end{align}
Note that $z$ is fixed.
The display above implies that $K^X(y)$ satisfies the differential equation in \eqref{PDEin(0,1)}.
The proof is now complete.

\hfill$\Box$

\bigskip
To prove Proposition \ref{prop2},
we first  recall some basics on exit measures of superprocesses. Let $S:= \R_+\times \R_+$ and we consider the evolution of the supper process in $S$. Let $\mathbb{O}\subset \mathcal{B}(\mathcal{S})$ the class of open subsets of $S$. Roughly speaking, we obtain the exit measures $\{X_O; O\in \mathbb{O}\}$ by freezing ``particles'' once they exit $O$.
For supercritical branching Brownian motion, similar ideas but with
a different terminology ``stopping line"
are used in \cite{Ky}.   For applications of exit measures of supercritical super Brownian motion, one can see \cite{KLMR}.
Now we formally introduce the exit measures.
For any $r>0$ and $x\geq 0$, we use $\mathbf{P}_{r,x}$ to denote the law $\mathbf{P}\left(\cdot\ | W_{r}=x\right)$.
Let $\mathcal{B}(S)$ be the Borel $\sigma$-field on $S$,
and $\mathcal{M}_F(S)$ the space of finite Borel measures on $S$. A measure $\mu \in \mathcal{M}_F(\R_+)$ is identified with the corresponding measure on $S$ concentrated on $\{0\} \times \R_+ $.
According to Dynkin \cite{E.B1.}, there exists a family of random measures
$\{(X_Q, \mathbb{P}_\mu); Q\in \mathbb{O}, \mu\in \mathcal{M}_F(S)\}$ such that for any $Q\in \mathbb{O}$,
$\mu\in \mathcal{M}_F(S)$ with $\textup{supp}\  \mu \subset Q$, and bounded non-negative Borel function $f(t,x)$ on $S$,
$$
\mathbb{E}_\mu \left(\exp\left\{-\langle f, X_Q \rangle\right\}\right) = \exp\big\{-\langle v_f^{X,Q} , \mu \rangle \big\},
$$
where
$v^{X,Q}_f(s,x)$
is the unique positive solution of the equation
\begin{align}\label{Equation-exit-measure}
	v_f^{X,Q} (s,x) & = \mathbf{E}_{s,x} \left(f(\tau, W_\tau^0)\right)- \mathbf{E}_{s,x} \int_s^\tau \varphi\big(v_f^{X,Q}(r, W_r^0)\big)\mathrm{d} r \nonumber\\
	&= \mathbf{E}_{s,x} \big(f(\tau, W_{\tau\land \tau_0^{W,-} })\big)- \mathbf{E}_{s,x} \int_s^{\tau\land \tau_0^{W,-}} \varphi \big(v_f^{X,Q}(r, W_r)\big)	\mathrm{d} r,
\end{align}
with $\tau:= \inf\left\{ r: (r, W_r) \notin Q\right\}$.
For $Q= D_z:=(0,\infty) \times (0, z)$,
$\tau= \tau_0^{W,-}\land \tau_z^{W,+}$.
Taking $f(t,x)= \theta 1_{\{x>0\}}$ in \eqref{Equation-exit-measure} and using the time-homogeneity of $W$,
we get that $v_f^{X,D_z}(s,x)=: v_\theta^{X,D_z}(x)$ is independent of $s$ and
is the unique positive solution of the equation of
\begin{align}\label{Equation-exit-measure-2}
	v_\theta^{X,D_z} (x) = \theta \mathbf{P}_{x} \big(\tau_z^{W, +} <\tau_0^{W,-} \big)- \mathbf{E}_{x} \int_0^{\tau_z^{W,+}\land \tau_0^{W,-}} \varphi \big(v_\theta^{X, D_z}( W_r)\big)\mathrm{d} r.
\end{align}
Moreover, by \eqref{Coupling},
\begin{align}\label{Solution-Equation-exit-measure-2}
	v_\theta^{X,D_z} (x) &= -\log \E_{\delta_x} \left(\exp\left\{- \theta X_{D_z} ([0,\infty)\times \{z\} ) \right\}\right) \nonumber\\
	& = -\log \E_{\delta_x} \big(\exp\big\{- \theta X_{D_z}^{(0,\infty)} ([0,\infty)\times \{z\} ) \big\}\big).
\end{align}
Letting $\theta\to +\infty$ in the display above, by the definition of $X_{D_z}$, we see that
\begin{align}
	v_\infty^{X, D_z}(x) & = -\log \P_{\delta_x} \big( X_{D_z}^{(0,\infty)} ([0,\infty)\times \{z\} )=0\big)\nonumber\\
	& = -\log \P_{\delta_x} \big( M^{(0,\infty),X} <z\big).
\end{align}

\noindent
\textbf{Proof of Proposition \ref{prop2}: }
Note that if $K^X$ is a solution to the problem in \eqref{PDEin(0,1)}, then for any $0<y<z<1$,
\begin{align}\label{step_46}
	K^X(y)+\mathbf{E}_y\Big(\int_0^{\tau_z^{W,+}\land \tau_0^{W,-}} \varphi(K^X((W_s))\mathrm{d}s \Big)=\mathbf{E}_y\big(K^X(W_{\tau_z^{W, +}\land \tau_0^{W,-}})\big).
\end{align}
Thus, for each fixed $z\in (0,1)$, $K^X(y)$ is a solution to the equation
\begin{align}\label{step_45}
	& f(y)+\mathbf{E}_y\Big(\int_0^{\tau_z^{W,+}\land \tau_0^{W,-}} \varphi(f(W_s))\mathrm{d}s \Big)\nonumber\\
	&=\mathbf{E}_y\big(K^X(W_{\tau_z^{W,+}\land \tau_0^{W,-}})\big)= K^X(z)\mathbf{P}_y\big(\tau_z^{W,+}<\tau_0^{W,-}\big),\quad y\in (0,z),
\end{align}
where the last inequality holds since $K^X(0)=0$.
By
\eqref{Equation-exit-measure-2} and \eqref{Solution-Equation-exit-measure-2},
\eqref{step_45} has a unique solution
given by
\begin{align}\label{step_47}
	v_{K^X(z)}^{X, D_z}(y)= -\log \E_{\delta_y}\big(\exp\big\{- K^X(z)X_{D_z}^{(0,\infty)}([0,\infty)\times \{z\}) \big\}\big).
\end{align}
Since $K^X$ is a solution to \eqref{step_45}, we have
\begin{align}\label{step_48}
	K^X(y)= -\log \E_{\delta_y}\big(\exp\big\{- K^X(z) X_{D_z}^{(0,\infty)}([0,\infty)\times \{z\} )\big\}\big), \quad y\in (0,z).
\end{align}
On one hand,
\begin{align}\label{step_49}
	K^X(y)\leq -\log \P_{\delta_y}\big(X_{D_z}^{(0,\infty)}([0,\infty)\times \{z\} )=0\big)= -\log \P_{\delta_y}\big(M^{(0,\infty),X}<z\big).
\end{align}
On the other hand, for any fixed $z_0\in (y,1)$, we choose $z\in (z_0, 1)$ so that $K^X(z)>K^X(z_0)$. Then
\begin{align}\label{step_50}
	K^X(y)\geq -\log \E_{\delta_y}\big(\exp\big\{- K^X(z_0) X_{D_z}^{(0,\infty)}([0,\infty)\times \{z\} )\big\}\big)=: K^X_{z_0}(y;z).
\end{align}
Note that $K^X_{z_0}(\cdot;z)$ is the unique bounded solution to
\[
K^X_{z_0}(y;z)+\mathbf{E}_y\Big(\int_0^{\tau_z^{W,+}\land \tau_0^{W,-}} \varphi(K^X_{z_0}(W_s;z))\mathrm{d}s \Big)=K^X(z_0)\mathbf{P}_y\big(\tau_z^{W,+}<\tau_0^{W,-}\big).
\]
Define $\widehat{K}^X_{z_0}(y):= z^{\frac{2}{\alpha-1}}K^X_{z_0}(yz;z)$, then the above equation is equivalent to
\begin{align}\label{step_57}
	&\widehat{K}^X_{z_0}\left(\frac{y}{z}\right)+ z^{\frac{2}{\alpha-1}}\mathbf{E}_y\Big(\int_0^{\tau_z^{W,+}\land \tau_0^{W,-}} \varphi\big(z^{-\frac{2}{\alpha-1}}\widehat{K}^X_{z_0}(z^{-1}W_s;z)\big)\mathrm{d}s \Big)\nonumber\\
	&=z^{\frac{2}{\alpha-1}}K^X(z_0)\mathbf{P}_y\big(\tau_z^{W,+}<\tau_0^{W,-}\big).
\end{align}
Using the scaling property of Brownian motion and the fact that $x^{\frac{2}{\alpha-1}}\varphi\big(x^{-\frac{2}{\alpha-1}} v\big) x^2 = \varphi(v)$, \eqref{step_57} is equivalent to
\begin{align}
	&\widehat{K}^X_{z_0}\left(\frac{y}{z}\right)+ \mathbf{E}_{y/z}\Big(\int_0^{\tau_1^{W,+}\land \tau_0^{W,-M}} \varphi\big(\widehat{K}^X_{z_0}(z^{-1}W_s;z)\big)\mathrm{d}s \Big)\nonumber\\
	&=z^{\frac{2}{\alpha-1}}K^X(z_0)\mathbf{P}_{y/z}\big(\tau_1^{W,+}<\tau_0^{W,-}\big).
\end{align}
Again using the uniqueness of the solution to \eqref{Equation-exit-measure-2}, we conclude that
\begin{align}\label{step_58}
	& -\log \E_{\delta_y}\big(\exp\big\{- K^X(z_0) X_{D_z}^{(0,\infty)}([0,\infty)\times \{z\} )\big\}\big)= K^X_{z_0}(y;z)\nonumber\\
	& =z^{-\frac{2}{\alpha-1}}
	\widehat{K}^X_{z_0}\left(\frac{y}{z}\right)= \Big(-\log \E_{\delta_{y/z}}\big(\exp\big\{- z^{\frac{2}{\alpha-1}}K^X(z_0) X_{D_1}^{(0,\infty)}([0,\infty)\times \{1\} )\big\}\big)\Big)\cdot z^{-\frac{2}{\alpha-1}}\nonumber\\
	& \geq \Big(-\log \E_{\delta_{y/z}}\big(\exp\big\{-
	z_0^{\frac{2}{\alpha-1}}
	K^X(z_0) X_{D_1}^{(0,\infty)}([0,\infty)\times \{1\} )\big\}\big)\Big)\cdot z^{-\frac{2}{\alpha-1}}.
\end{align}
Therefore, plugging \eqref{step_58} into \eqref{step_50} and then letting $z\to 1-$ in both \eqref{step_49} and \eqref{step_50}, we conclude that
\begin{align}
	& -\log \P_{\delta_y}\big(M^{(0,\infty),X}<1\big)\geq K^X(y)\geq -\log \E_{\delta_y}\big(\exp\big\{-
	z_0^{\frac{2}{\alpha-1}}
	K^X(z_0) X_{D_1}^{(0,\infty)}([0,\infty)\times \{1\} )\big\}\big).
\end{align}
Since $K^X(z_0)\to +\infty$ as $z_0 \to 1-$,  letting $z_0\to 1-$ in the above inequality yields that
\begin{align}
	& -\log \P_{\delta_y}\big(M^{(0,\infty),X}<1\big)\geq K^X(y)\geq -\log \P_{\delta_y}\big(X_{D_1}^{(0,\infty)}([0,\infty)\times \{1\} )=0\big)\nonumber\\
	& = -\log \P_{\delta_y}\big(M^{(0,\infty),X}<1\big),
\end{align}
which implies that $K^X(y) = -\log \P_{\delta_y}\left(M^{(0,\infty),X}<1\right)$.  This completes the proof of the proposition.

\hfill$\Box$

\noindent

\begin{singlespace}

\end{singlespace}

\vskip 0.2truein
\vskip 0.2truein

\noindent{\bf Haojie Hou:}  School of Mathematical Sciences, Peking
University,   Beijing, 100871, P.R. China. Email: {\texttt
houhaojie@pku.edu.cn}

\smallskip

\noindent{\bf Yan-Xia Ren:} LMAM School of Mathematical Sciences \& Center for
Statistical Science, Peking
University,  Beijing, 100871, P.R. China. Email: {\texttt
yxren@math.pku.edu.cn}

\smallskip
\noindent {\bf Renming Song:} Department of Mathematics,
University of Illinois,
Urbana, IL 61801, U.S.A.
Email: {\texttt rsong@math.uiuc.edu}

\end{document}